\def\IH{{\mathbb H}} 
 
\def\IR{{\mathbb R}} 
\def\ax{{\rm axis}} 
\def\tr{{\rm tr}} 
\def\IQ{{\mathbb Q}}

\def\IS{{\mathbb S}} 
 
\def\IZ{{\mathbb Z}} 
\def\s{\sigma} 
 
\def\IC{\mathbb C} 
\def\ID{{\mathbb D}}
\def\oC{\hat{\IC}} 
\def\G{\Gamma} 

\def\arccosh{{\rm arccosh}}

\def\i{{\bf i}} 
\def\j{{\bf j}} 
\def\k{{\bf k}}

\newcommand{\HS}[3]{\left({ #1, #2 \over #3} \right)}

\documentclass[twoside]{irmaems}
\usepackage{graphicx}
\usepackage{amssymb}  
\usepackage{amsmath}  
\usepackage{latexsym} 
\usepackage{epsfig,multicol}

\renewcommand\theequation{\arabic{section}.\arabic{equation}}

\theoremstyle{definition}  

 \newtheorem{definition}{Definition}[section]

\theoremstyle{plain}

 \newtheorem{theorem}[definition]{Theorem}
 \newtheorem{corollary}[definition]{Corollary}
 \newtheorem{lemma}[definition]{Lemma}

\begin{document}

\title{ {\large The Geometry and Arithmetic of Kleinian Groups } \\  { $\;\;$} \\{\small To the memory of  Fred Gehring and Colin Maclachlan}
\\
} 
 \author{Gaven J. Martin\thanks{
Work partially supported by the New Zealand Marsden Fund}
\address{
Institute for Advanced Study\\
Massey University,  Auckland\\
New Zealand\\
email:\,\tt{g.j.martin@massey.ac.nz}
\\
}}
\maketitle 
 \begin{abstract} In this article we survey and describe various aspects of the geometry and arithmetic of Kleinian groups - discrete nonelementary groups of isometries of hyperbolic $3$-space. In particular we make a detailed study of two-generator groups and discuss the classification of the arithmetic generalised triangle groups (and their near relatives).  This work is mainly based around my collaborations over the last two decades with Fred Gehring and Colin Maclachlan,  both of whom  passed away in 2012. There are many others involved as well.  

Over the last few decades the theory of Kleinian groups has flourished because of its intimate connections with low dimensional topology and geometry.  We give  little of the general theory and its connections with $3$-manifold theory here, but focus on two main problems: Siegel's problem of identifying the minimal covolume hyperbolic lattice and the Margulis constant problem.  These are both ``universal constraints'' on Kleinian groups -- a feature of discrete isometry groups in negative curvature and include results such as J\o rgensen's inequality,  the higher dimensional version of Hurwitz's $84g-84$ theorem and a number of other things.   We will see  that big part of the work necessary to obtain these results is in getting concrete descriptions of  various analytic spaces of two-generator Kleinian groups,  somewhat akin to the Riley slice.
  \end{abstract}

\newpage

\tableofcontents

\section{Introduction} 

Over the last few decades the theory of Kleinian groups has flourished because of its intimate connections with low dimensional topology and geometry.  The culmination must certainly be Perelman's proof of Thurston's geometrisation conjecture which states that compact 3-manifolds can be decomposed canonically into submanifolds that have geometric structures, \cite{MT1,MT2}. This is an analogue for 3-manifolds of the uniformization theorem for surfaces. This conjecture  implied, for instance,  the Poincar\'e conjecture.  There have been many other recent advances.  These include the density conjecture \cite{Agol,CG},  the ending lamination conjecture \cite{BCM},  the surface subgroup conjecture \cite{KMar} and the virtual Haken conjecture \cite{Agol2}. While we will not discuss these results here (nor offer statements of precise theorems) together they give a remarkably complete picture of the structure of hyperbolic group actions and their quotient spaces in three-dimensions.  

\medskip

Our focus here will be on another aspect of the theory of Kleinian groups. Namely their geometry and the geometric constrains placed on them in order to be discrete - early examples of these are \cite{Ap, Jorgensen, Jorgensen2, Mard}.  Here we will focus on identification of  minimal covolume actions,  collaring theorems,  singular set structure,  geometric decomposition theorems (thick and thin),   the Margulis constant and so forth.  Precise definitions and statements about what is known can be found in the body of this article.  In particular,  in 1943 Siegel \cite{Siegel,Siegel2} posed the problem of identifying the minimal covolume lattices of isometries of hyperbolic $n$--space,  or more
generally rank--$1$  symmetric spaces.  He solved the problem in two dimensions identifying the $(2,3,7)$--triangle group as the unique lattice of minimal
coarea.  This group is the orientation-preserving index-two subgroup of the group generated by reflections across the sides of a hyperbolic triangle with interior angles $\pi/2,\pi/3$ and $\pi/7$. This tessellation of the hyperbolic plane is illustrated below.

\scalebox{0.2}{\includegraphics[viewport= -490 200 580 770]{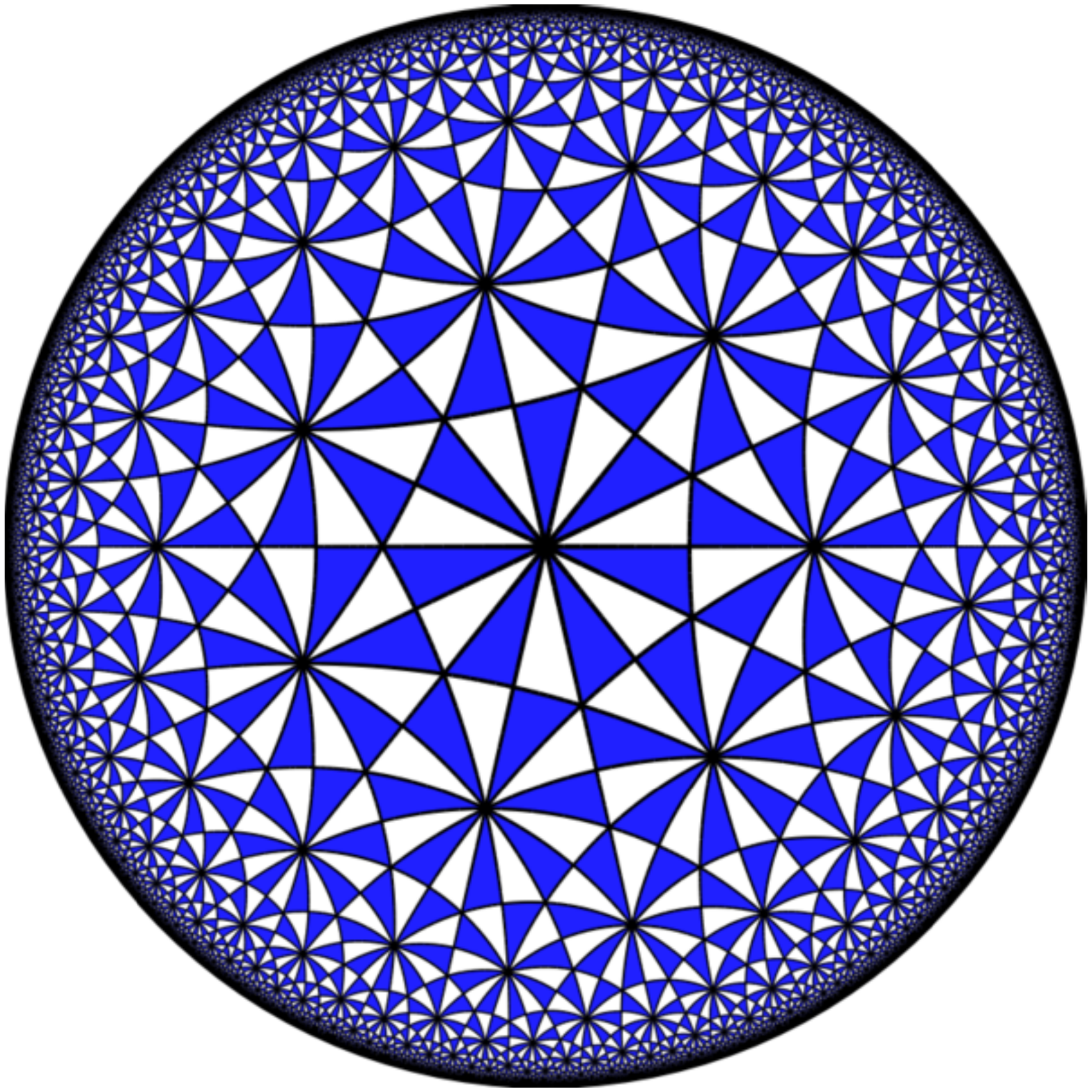}}

\bigskip
\begin{center} The $(2,3,7)$-tessellation of the hyperbolic plane.
\end{center}

Siegel in fact proved what has come to be known as the signature formula from which one may deduce the complete spectrum of coareas of lattices of
the hyperbolic plane.  Ka{\u z}dan and Margulis \cite{KM} showed in 1968 that for each $n$ the infimum of the covolume of lattices is
positive and achieved,  while answering a question of Selberg.

At the time of Siegel's result the theory of covering spaces was not well developed and he could only suggest a connection
between minimal coarea lattices and Hurwitz's $84g-84$ theorem of 1892 \cite{Hur}, bounding the order of the symmetry group of a
Riemann surface in terms of its genus.  This suggested connection was due to the many close analogies and ideas used in the proof and was confirmed by Macbeath
in 1961, \cite{Macbeath}. Selberg's Lemma \cite{Selberg} established the existence of torsion-free subgroups of finite-index in
hyperbolic lattices,  among other things.  As a consequence of the Mostow rigidity theorem \cite{Mostow} the $84g-84$ theorem takes its expression in
terms of bounding the order of the symmetry group of a hyperbolic manifold by its volume,  see Theorem \ref{3dhur} and its preamble below for a discussion. 

\medskip  

At the same time as the excitement over the connections between $3$-manifold theory and hyperbolic geometry was running, the theory of arithmetic Kleinian groups was developing with many similar objectives.  This is mostly admirably accounted in the book of  C. Maclachlan and A. Reid,  \cite{MRbook}.  Due to formulas of Borel,  \cite{Borel},  various explicit calculations can be made in the case a Kleinian group is arithmetic (in particular one can determine the maximal arithmetic Kleinian group in which it embeds and various number-theoretic constructs associated with the group).  This makes the associated manifolds and orbifolds amenable to the use of techniques from algebra and number theory, as well as the well developed topological, analytical and geometric tools.  Earlier work of ours \cite{GMMR}, motivated by the problems discussed here,  give explicit criteria which determine when a two-generator group is arithmetic.   It turns out that nearly all the extremal problems one might formulate are realised by arithmetic groups,  perhaps the number theory forcing additional symmetries in a group and therefore making it ``smaller'' or ``tighter''.  Indeed we will see that many of these externals are arithmetic groups generated by two elements of finite order.  With C. Maclachlan we proved that there are only finitely many such arithmetic groups \cite{McM} -- of course the orbifold Dehn surgeries on two-bridge links give infinitely many examples of finite covolume groups generated by two elements of finite order.  It would seem this finite family of arithmetic groups contains all the extremal groups -- for instance the minimal covolume lattice,  the minimal covolume non-uniform lattice,  the minimal covolume torsion-free lattice  and the minimal covolume torsion-free non-uniform lattice are all two-generator and arithmetic.

\medskip

In order to solve Siegel's problem and other related extremal problems such as those mentioned above,  the procedure we adopt is the following.  We begin by describing the space of two-generator discrete groups using various generalisations of J\o rgensen's inequality \cite{Jorgensen} based around a very intriguing family of polynomial trace identities in $\mathrm{SL}(2,\IC)$ found by us in  \cite{GM3}.  From this we get general {\em a priori} bounds (universal constraints) which hold when various assumptions are made which,  in general,  will not be true for the extremal case.  These assumptions will be along the following lines:  the vertices of the singular set are no closer than a given number,  or the collaring radius of a torsion element is at least such and such.  We then go about trying to prove that if these assumptions fail,  then specific geometric configurations must occur which produce a two-generator subgroup of a specific type.  This is basically possible since the failure of these assumptions confines us to a small part of the space of two-generator   discrete groups which we will have completely described.  This  allows us in turn to identify a finite list of candidate groups to consider.  These candidates will (fortunately) be arithmetic.  This then leads us to another part of the programme which is to enumerate all the two-generator arithmetic Kleinian groups.  We will report on what we know about that problem here as well.  However,  once we have shown that the examples which our {\em a priori} bounds do not cover are arithmetic,  we can use other techniques to identify the extremals from this finite list.  In this way we will see that most examples of small volume hyperbolic orbifolds and manifolds have an arithmetic structure,  and more generally the same is true for other sorts of extremal problems.  

\medskip

Basically our programme of proof is predicated on the belief that the externals for these geometric problems are two-generator arithmetic,  and luckily this is the case.

\section{The volumes of hyperbolic orbifolds}
  
As we have mentioned above,  in 1943 Siegel posed the problem of  identifying the smallest covolume hyperbolic lattices in $n$ dimensions.  He defined the numbers
\[  \mu(n) = \inf_\Gamma  {\rm vol}_{\IH}\Big(\IH^n/\Gamma\Big) \]
where the infimum is taken over all the lattices $\Gamma\subset {\rm Isom}^+(\IH^n)$ of hyperbolic $n$-space. 
For Euclidean lattices (Bieberbach or Crystallographic groups) such an infimum is obviously equal to $0$ as one can tessellate $\IR^n$ by arbitrarily small hypercubes.  It is in fact the negative curvature that is important here to ensure $\mu(n)>0$. After Siegel,  the first key result was that of D. Ka{\u z}hdan and G. Margulis to show $\mu(n)>0$ and is attained.  (Indeed there is a more general discussion to be had about geometric constraints in negative curvature,  see for instance \cite{Martinnegcurv}.)

In two dimensions Siegel showed 
\begin{theorem} 
$\mu(2) = \frac{\pi}{21}$ and this value is uniquely attained for the $(2,3,7)$-triangle group.
\end{theorem}
In fact Siegel refined the Poincar\'e-Klein area formula to the well known {\em signature formula}
\begin{theorem}  Let $\Sigma$ be a Riemann surface of genus $g$ with $n$ punctures (cusp points) and $M$ cone points of angle $2\pi/m_j$, $ j=1,2,\ldots,M$.  The the hyperbolic area of $\Sigma$ is
\[ A = 2\pi \Big(2g-2+N+\sum_{j=1}^{M} (1-1/m_j) \Big)  .\]
\end{theorem} 
There can be no such simple formula for the covolume of lattices in higher dimensions.  In even dimensions $n=2\ell$,  the Chern-Gauss-Bonnet theorem \cite{Chern} gives an explicit formula in the case of manifolds.  For instance in four dimensions this volume is a constant multiple of the Euler characteristic,  which is even.  But it is not known if the smallest multiple is one or two,  see \cite{CM} for  recent results and discussion.

In odd dimensions the Euler characteristic vanishes and very little is known in general while there is good information in the arithmetic case, \cite{Belo}.

As we mentioned,  Siegel suggested a connection with Hurwitz' $84g-84$ theorem on the symmetries of Riemann surfaces based on the nature of his method of proof.  This connection was later confirmed by McBeath in 1961 as the theory of covering spaces developed in the geometric setting and ideas from combinatorial group theory came into play.

 \bigskip
 
In the first instance we know what the spectrum of possible covolumes of hyperbolic lattices in three-dimensions looks like.

\begin{theorem}[J\o rgensen, Thurston]  The spectrum of volumes of lattices of hyperbolic three-space is well ordered and of type $\omega^\omega$.  That is there is a smallestvalue, then there is a smallest limit value,  then a smallest limit of limit values and so forth.  Limits only occur from below and each particular value is obtained with only finite multiplicity.  That is,  up to conjugacy,  there are only finitely many lattices with a given volume.
\end{theorem} 
\noindent See for instance \cite{Gromov} for the manifold case.  This theorem is in striking contrast to Wang's theorem \cite{Wang}:
\begin{theorem} If $n\geq 4$,  for each $V\geq 0$ there are only finitely many   isometry classes of $n$-dimensional hyperbolic manifolds $M^n$ with volume ${\rm vol}_{\IH}(M^n)\leq V$.
\end{theorem}
\noindent This actually implies that for each $n$  the spectrum of possible volumes of  lattices of hyperbolic $n$-space is discrete.
 
 \medskip
 
In a sequence of works,  culminating in papers with Gehring \cite{GMannals} and with Marshall \cite{MMannals} we proved the following theorem identifying the two smallest covolumes of lattices.
\begin{theorem}\label{mainthm}  Let $\Gamma$ be a Kleinian group.  Then either
\begin{eqnarray*}\label{1}
{\rm vol}_{\IH}(\IH^3/\Gamma) &=& {\rm vol}_{\IH}(\IH^3/\G_0) = 275^{3/2} 2^{-7} \pi^{-6}
\zeta_k(2)\sim 0.0390  \;{\rm and}\; \Gamma=\Gamma_0,  \\
{\rm vol}_{\IH}(\IH^3/\Gamma) &=& {\rm vol}_{\IH}(\IH^3/\G_1)= 283^{3/2} 2^{-7} \pi^{-6}
\zeta_{k'}(2) \sim  0.0408 \;{\rm and}\;
\Gamma=\Gamma_1, \\
{\rm or}\;\;\; {\rm vol}_{\IH}(\IH^3/\Gamma) &>& 0.041.
\end{eqnarray*}
\end{theorem}
\noindent The equality of groups here is  up to conjugacy.  Here is a description of 
the  two groups $\Gamma_0$ and $\Gamma_1$ and the associated arithmetic data;
\begin{itemize}
\item 
$\Gamma_0$ is a two-generator  arithmetic Kleinian group obtained as 
a $\IZ_2$--extension of the index-2
orientation-preserving subgroup of the group generated by reflection
in the faces of the 3-5-3--hyperbolic Coxeter
tetrahedron and $\zeta_k$ is the Dedekind zeta function of the 
underlying number field $\IQ(\gamma_0)$,  with $\gamma_0$ a complex 
root of $\gamma^4+6\gamma^3+12\gamma^2+9\gamma+1=0$, of discriminant 
$-275$.  The associated quaternion algebra is unramified.  This group 
has a discrete and faithful representation in $\mathrm{SL}(2,\IC)$, determined 
uniquely up to conjugacy,  generated by  two matrices $A$ and $B$ 
with $\tr^2(A)=0$, $\tr^2(B) = 1$ and $\tr(ABA^{-1}B^{-1})-2 = 
\gamma_0$.  

\item $\Gamma_1$   is a two-generator  arithmetic 
Kleinian group and $\zeta_{k'}$ is the Dedekind zeta function of the 
underlying number field $\IQ(\gamma_1)$,  with $\gamma_1$ a complex 
root of $\gamma^4+5\gamma^3+7\gamma^2+3\gamma+1=0$, of discriminant 
$-283$.  The associated quaternion algebra is unramified.  This group 
has a discrete and faithful representation in $\mathrm{SL}(2,\IC)$, determined 
uniquely up to conjugacy,  generated by  two matrices $A$ and $B$ 
with $\tr^2(A)=0$, $\tr^2(B) = 1$ and $\tr(ABA^{-1}B^{-1})-2 = 
\gamma_1$. 
\end{itemize}

\bigskip

In particular we solve Siegel's problem in dimension $3$,
\begin{corollary}
\[ \mu_3= 275^{3/2} 2^{-7} \pi^{-6} \zeta_k(2) = 0.03905\ldots . \]
\end{corollary}
In fact, as the minimal volume orbifold group $\G_0$ above is the orientable double cover of a non-orientable orbifold group,  Theorem \ref{mainthm} also identifies the unique  minimal volume non-orientable orbifold as having hyperbolic volume equal to $\mu_3/2$.
 
 \medskip
 
 We now develop the following theorem which must be regarded as an analogue of Hurwitz's $84g-84$ theorem.  A modern proof of Hurwitz's results would be the following:
 Suppose a group $G$ acts faithfully by homeomorphism on a finite-area hyperbolic surface $F^2$ of genus $g$.  The  solution of the Nielsen realisation problem by S. Kerchoff, \cite{Ker}, shows that there is a possibly different hyperbolic structure on the same surface $F^2$ in which $G$ acts faithfully by isometries.  The area of this surface is a topological invariant,  $2\pi(g-1)$,  and does not depend on the choice of metric.  Next,  the quotient space $F^2/G$ will be a hyperbolic $2$-orbifold and Siegel's result in two dimensions shows $\mathrm{Area}_{\IH^2}(F^2/G) \geq \pi/42$ (nonorientable case) with equality uniquely attained for the  $(2,3,7)$-reflection group.  But then  \[ \pi/42 \leq Area_{\IH^2}(F^2)/G =  Area_{\IH^2}(F^2)/|G| = 2\pi(g-1)/|G|, \]
 which yields $|G|\leq 84 g-84$.  This also identifies the finite quotients of the $(2,3,7)$-reflection  group as those groups which act maximally on a surface of genus $g$ (however not every $g$ can be realised,  and this is another interesting story).  In three dimensions the argument is the same except that the use of the Nielson realisation problem is obviated as a consequence of Mostow rigidity \cite{Mostow} which says that if $G$ acts by homeomorphism,  then $G$ acts by isometry (it of course says more than this). Here is a version of the Mostow rigidity theorem.
 
 \begin{theorem}[Mostow rigidity] Let $G$ and $H$ be algebraically isomorphic subgroups of the isometry group of hyperbolic three-space and suppose that the volume of $\IH^3/G$ is finite. Then $\IH^3/G$ and $\IH^3/H$ are isometric; there is an isometry $f:\IH^3\to\IH^3$ with $G=f\circ H\circ f^{-1}$.  
\end{theorem}
 As $\IH^3$ is contractible, the condition that $G$ and $H$ are isomorphic is equivalent to an assumption such as $\IH^3/G$ and $\IH^3/H$ are homotopy equivalent. Therefore Mostow's theorem implies,  among many other things,  that volume is a topological invariant of hyperbolic manifolds.  In particular a homeomorphism of a finite-volume hyperbolic manifold $M$ is homotopic to an isometry. A stronger version of this result can be found in  \cite{GMT}.
 
  Consequently we have the following.
 
 \begin{theorem}\label{3dhur}  Let $M$ be a finite-volume hyperbolic $3$-manifold and $G$ a group of orientation-preserving homeomorphisms acting faithfully on $M$.  Then
 \begin{equation}\label{maxorder}
 |G| \leq \frac{1}{\mu_3} {\rm vol}_{\IH}(M).
 \end{equation}
 \end{theorem}
 We see $\mu_{3}^{-1} \approx 26.082$  and in the non-orientable case one simply multiplies by $2$.  This leads to the following elementary consequence (given the rather large hammers above).
 
  \begin{theorem}  Let $\G$ be a finite-covolume hyperbolic lattice, $\G\subset \mathrm{Isom}(\IH^3)$.  Let $G$ be a group of isomorphisms of $\G$.  Then
 \begin{equation}\label{maxorder}
 |G| \leq \frac{2}{\mu_3} {\rm covol}_{\IH}(\G).
 \end{equation}
 This results is sharp and attained infinitely often.  Further $G$ can be realised as a group of automorphisms in $Isom(\IH^3$) of $\G$.
 \end{theorem}
 
 We next recall a very nice and quite relevant result of S. Kojima \cite{Kojima}.
 
\begin{theorem}
 Every finite group can be realized as the full isometry group of some closed hyperbolic 3-manifold.
 \end{theorem}
 Here ``full'' means that there is no larger symmetry group acting.   The proof is based on  modifications of a similar theorem of Greenberg for compact surfaces.  Next,  the natural question is what groups can act on a hyperbolic $3$-manifold so (\ref{maxorder}) holds.  These are groups of maximal order.  The question becomes one of identifying the finite-index normal torsion-free subgroups of $\G_0$. The existence of such groups is a fairly general result known as Selberg's lemma, \cite{Selberg},  and applies to all linear groups.  In our case we write it as follows.   
 
\begin{theorem}[Selberg's lemma] Let $G$ be a finitely generated subgroup of the isometry group of hyperbolic three-space.  Then $G$ contains a normal torsion-free subgroup of finite-index.
\end{theorem}

In fact Selberg's lemma will imply there are infinitely many normal torsion-free subgroups $H$.  If $G$ is discrete,  then $H$ will obviously be discrete and hence $\IH^3/H$ is a hyperbolic manifold on which $G/H$ will act by isometry.  This connects the theory of hyperbolic orbifolds and manifolds.   In the particular case at hand we know of the following examples from \cite{CMT}

\begin{theorem} We have that
\begin{itemize}
\item for every prime $p$, there is some $q=p^k$, with $k\leq8$, such that either $\mathrm{PSL}(2,q)$ or $\mathrm{PGL}(2,q)$ has maximal order, 
\item for all but finitely many $n$, both the alternating group $A_n$ and the symmetric group $S_n$ have maximal order.
\end{itemize}
\end{theorem}
   Actually in that paper we  describe all torsion-free subgroups of index up to $120$ in $\G_0$ and explain how other infinite families of quotients can be constructed.
   
   \medskip  
   
As concerns arithmeticity,  we note that previously T. Chinburg and E. Friedman had identified $\G_0$ as the smallest covolume arithmetic hyperbolic lattice in three dimensions, \cite{CF}.   M. Belolipetsky has found the smallest volume arithmetic hyperbolic lattices in higher even dimensions,  \cite{Belo}.

In the case that the lattice is known not to be compact, using an important result of K. B\"or\"oczky \cite{Boroc},  R. Meyerhoff \cite{Meyer1} has shown that
\begin{theorem}
$PGL(2,{\cal O}(\sqrt{-3}))$ is the smallest covolume non-compact lattice,  its covolume is $\approx 0.0846$.  
\end{theorem}
The volume of the associated  orbifold is $1/12$ the volume of the ideal tetrahedron and contains the  Coexter group $3-3-6$ with index $2$. This work was further extended by C. Adams who identified the next six such smallest cusped orbifolds,  \cite{Adams2}.

We will need to establish the following improvement of Meyerhoff's result,  see \cite{GMJGD}.  The proof is based around new inequalities for discrete groups similar to J\o rgensen's inequality.  The approach to these inequalities was first discussed in \cite{GM0} using new polynomial trace identities and iteration.

\begin{theorem}
If $\G$  is a Kleinian group containing an elliptic element of order
$p\geq 6$,  then
\[ {\rm vol}_{\IH}(\IH^3/\G) \geq {\rm vol}_{\IH}(\IH^3/PGL(2,{\cal O}_{\sqrt{-3}}))=0.0846\ldots. \]
This result is sharp.
\end{theorem}

C. Adams \cite{Adams} has identified the smallest three limit volumes as $0.30532\ldots$,  $0.444451\ldots$ and $0.457982\ldots$. The smallest volume here is a quotient of the Borromean rings group.   Further,  Adams  gave a proof of the following ``folklore'' theorem.

\begin{theorem}  Let ${\cal O}_i = \IH^3/\Gamma_i$ be a sequence of hyperbolic $3$-orbifolds all of whose volumes are bounded by a constant $\alpha$.  Then there exists a subsequence ${\cal O}_{i_j}$ so that all the orbifolds come from Dehn filling a single hyperbolic $3$-orbifold.
\end{theorem}

The minimal volume hyperbolic $3$-manifold has been identified by D. Gabai, R. Meyerhoff and P. Milley, \cite{GabMM}.  This is the culmination of a long and difficult sequence of work  and involves quite different techniques than those we develop here.

\begin{theorem} The Weeks-Matveev-Fomenko  manifold is the smallest volume hyperbolic $3$-manifold.  Its volume is
\[ \frac{ 3 \times 23^{3/2} \zeta_k(2)}{4 \pi^4 } \approx 0.942707\ldots . \]
Here $\zeta_k$ is the Dedekind zeta function of the field of degree $3$ with one complex place and discriminant $-23$, $k=\IQ(\gamma)$,  $\gamma^3-\gamma+1=0$.  This manifold can be obtained as $(5,2)$ and $(5,1)$-Dehn surgery on the  Whitehead link complement.
\end{theorem}

Finally,  in the non-compact case the following  is known  through the work of C. Cao and Meyerhoff \cite{CaoMeyer}.
\begin{theorem}
The figure-eight knot complement admits a hyperbolic structure of finite covolume,  and this with its sister,  $(5,1)$-Dehn surgery on one component of the Whitehead link complement are the smallest volume non-compact manifolds.  Their volume is $2.02988\ldots$ which is twice the volume of the ideal tetrahedron.
\end{theorem}
 
 The figure eight knot complement is actually a double-cover of the nonorientable Gieseking manifold, which therefore has the smallest volume among non-compact hyperbolic 3-manifolds.  Remarkably A. Reid has shown the figure-eight knot is the only arithmetic knot complement,  \cite{Reid}.
 
 \bigskip
 
 Taken together,  these results show we have a fairly good understanding of the small volume hyperbolic manifolds and orbifolds.

\section{The Margulis constant for Kleinian groups}

Let $M=\IH^3/\Gamma$ be an orientable hyperbolic $3$-manifold. So $\G$ is a torsion-free Kleinian group.  A {\em Margulis number} for $M$ (actually for $\Gamma$) is a number ${\bf m}>0$ such that if there exists a point $p\in \IH^3$ that is moved a distance less than ${\bf m}$ by two elements $f,g\in \G$, then $f$ and $g$ commute, \cite{Margulis2}.  There is a universal Margulis number that works for all hyperbolic $3$-manifolds and the largest such is called the Margulis constant ${\bf m}_{3}$. (And of course this is true in all dimension). The first real estimate  for ${\bf m}_3>0.104$ was given by Meyerhoff \cite{Meyer2} and improved in \cite{GM3}. M. Culler has given numerical evidence that for a particular $3$-manifold, $0.616$ fails to be a Margulis number and thus probably  ${\bf m}_3<0.616$.  (Actually  ${\rm arccosh}\left((\sqrt{3}+1)/{2}\right)= 0.831446\ldots$ is shown to be the Margulis number for the  arithmetic  four fold cover of $(4,0)$ \& $(2,0)$ Dehn surgery on the 2 bridge link complement $6^{2}_{2}$ of Rolfsen's tables in \cite{MMannals}).  

\medskip

   One of the primary interests in the Margulis number of a manifold is that it provides a decomposition of the manifold into the thin part, consisting of cusps and solid tubes about short geodesics, and the thick part, consisting of points with injectivity radius at least ${\bf m}/2$. If one seeks to classify hyperbolic $3$-manifolds below any given volume, the larger that ${\bf m}$ is, the fewer the possibilities for the thick part. 
An interesting recent result of P. Shalen says that, up to isometry, there are only finitely many hyperbolic $3$-manifolds that do not have $0.29$ as a Margulis number with all the exceptions being closed.  In an earlier result with M. Culler,  they showed the following, \cite{CS}.

\begin{theorem}
If $M$ is a closed hyperbolic $3$-manifold whose first Betti number is at least $3$, then the Margulis number for $M$ is greater than $\log 3 = 1.0986\ldots$
\end{theorem}

More generally the same geometric analysis will apply to hyperbolic $3$-orbifolds.  However,  the commutativity of $f$ and $g$ must be replaced by another assumption (which will be equivalent in the torsion-free case). 

\medskip

Thus we make the following definition.  Let  ${\cal O}=\IH^3/\Gamma$ be  an orientable hyperbolic $3$-orbifold.  A Margulis number for ${\cal O}$ is a number ${\bf m}_\G>0$ with the property that if there is a point $p\in \IH^3$ that is moved a distance less than ${\bf m}_{\G}$ by two elements $f,g\in \G$, then $\langle f,g\rangle$ is elementary (that is virtually abelian).  We will see in a moment that  the elementary groups are completely classified.  This change,  from commutativity to virtually abelian, is necessary as the finite spherical triangle groups act as isometries of $\IH^3$ fixing a finite point of $\IH^3$ and their generators do not commute.

We let ${\bf m}_{\cal O}$ be the largest constant which works for all Kleinian groups.  Obviously ${\bf m}_{\cal O}\leq  {\bf m}_{3}$.    Later, in \S 10,   we will construct   the following group.

\begin{theorem}  There is a Kleinian group $\G_{{\cal O}}$ generated by three elements of order two for which the Margulis constant   is  $\bar{\bf m}_{\cal O} = 0.1324\ldots$
\end{theorem}
We conjecture this value to be the Margulis constant for Kleinian groups, $\bar{\bf m}_{\cal O} ={\bf m}_{\cal O}$.

The following results are  known about this constant,  see \cite{GM1}.
\begin{theorem}  Let $\G=\langle f_1,f_2, \ldots \rangle$ be a Kleinian group.
\begin{itemize}
\item  If some $f_i$ is parabolic,  then ${\bf m}_\G\geq 0.1829\ldots$.
\item	 If some $f_i$ is elliptic of order $p\geq 3$ and the axis of $f_i$ does not
meet the axis of some $f_j$,  then ${\bf m}_\G\geq 0.189\ldots$.
\item If $\langle f_i,f_j \rangle$ is elementary for all $i,j$ and the $f_i$ are not
all of order 2, then
${\bf m}_\G\geq 0.2309\ldots$.
\item If there are three $f_i$'s each of order two and with coplanar axes, then
${\bf m}_\G\geq 0.2088\ldots$.
\end{itemize}
Each estimate is sharp.
\end{theorem}
It is an important observation of T. Inada \cite{Inada} that this analysis covers all cases and hence we have the following.
\begin{corollary} The Margulis constant is achieved in a Kleinian group generated by three elliptic elements of order 2 and 
 ${\bf m}_{{\cal O}} \leq 0.132409\ldots$.
\end{corollary}
 
The most interesting case occurring in the proof of the above theorem consists of the study of  $(p,q,r)$--Kleinian groups.  These are Kleinian groups generated by three rotations of orders $p,q$ and $r$ whose axes of rotation are coplanar. These groups provide very many examples of groups with ``tight'' singular set and come up in determining the minimal volume orbifold as we look for extremal configurations.  It is therefore necessary to study these groups closely and we will do this a bit later.

\section{The general theory}  In order to motivate what follows we will recount a little of this general theory.  In three-dimensions Thurston has shown that there are eight different model geometries \cite{Thurston}.  These are the natural geometries that might occur on a three-manifold.

\subsection{Geometries in three-dimensions.} The geometries identified by Thurston are
\begin{itemize}
\item {\bf spherical geometry},  giving rise to closed 3-manifolds with finite fundamental group such as Poincar\'e homology spheres and lens spaces.   
\item {\bf Euclidean geometry},  giving rise to  exactly six orientable finite closed $3$-manifolds.   More generally we have the Bieberbach groups  acting as lattices in the group of Euclidean isometries.  
\item   $\IS^2\times \IR$,  there are four finite volume manifolds with this geometry.
\item $\IH^2\times\IR$,  giving examples that include the product of a hyperbolic surface with a circle, or more generally the mapping torus of an isometry of a hyperbolic surface. Orientable finite volume manifolds with this geometry are Seifert fiber spaces.  
\item The geometry  $\widetilde{\mathrm{SL}(2,\IR)}$,  the universal cover of $\mathrm{SL}(2,\IR)$.  Examples  include the manifold of unit vectors of the tangent bundle of a hyperbolic surface  and most Brieskorn homology spheres.
\item The geometry {\bf Nil} which is basically the geometry of the Heisenberg group.  Compact manifolds with this geometry include the mapping torus of a Dehn twist of a 2-torus, or the quotient of the Heisenberg group by the ``integral Heisenberg group''.
\item The geometry {\bf Sol}.   All finite volume manifolds with sol geometry are compact and are either the mapping torus of an Anosov map of the 2-torus  or their finite quotients.
\item  {\bf Hyperbolic geometry}, $\IH^3$.  This is the most prevalent geometric structure. Examples are given by the Seifert--Weber space, or  most Dehn surgeries on (most) knots and links, or most Haken manifolds. The solution of the geometrization conjecture implies that a closed 3-manifold is hyperbolic if and only if it is irreducible, atoroidal (It does not contain an embedded, non-boundary parallel, incompressible torus), and has infinite fundamental group. 
\end{itemize}
\subsection{Classification by fundamental group.}
Given a group $G$ acting by isometry on one of these spaces $X$ we can form the orbit space $X/G$,  that is, the quotient space of $X$ under the relation $x\sim y$ if there is $g\in G$ with $g(x)=y$.  If  $G$ acts {\em properly discontinuously},  that is, for every compact set $ K\subset X$
\[ \#\{g\in G: g(K)\cap K \neq \emptyset\} < \infty, \]
then the space $X/G$ with the metric it inherits is called an orbifold and the group $G$ is called the orbifold fundamental group (usually just fundamental group).  These orbifolds will be manifolds if $G$ acts without fixed points.  In terms of the fundamental group $\pi_1(M)$ of a finite volume quotient with a geometric structure we have the following:
\begin{itemize}
\item   if $\pi_1(M)$  is finite, then the geometric structure on $M$ is spherical and $M$ is compact,
\item   if $\pi_1(M)$  is virtually cyclic but not finite,  then the geometric structure on $M$ is $\IS^2\times \IR$, and $M$ is compact, 
\item   if $\pi_1(M)$  is virtually abelian but not virtually cyclic,  then the geometric structure on $M$ is Euclidean and $M$ is compact,
\item   if $\pi_1(M)$  is virtually nilpotent but not virtually abelian, then the geometric structure on $M$ is nil geometry, and $M$ is compact,
\item   if $\pi_1(M)$ is virtually solvable but not virtually nilpotent, then the geometric structure on $M$ is sol geometry, and $M$ is compact,
\item   if $\pi_1(M)$  has an infinite normal cyclic subgroup but is not virtually solvable, then the geometric structure on $M$ is either $\IH^2\times\IR$ or the universal cover of $SL(2, \IR)$, 
\item   if $\pi_1(M)$  has no infinite normal cyclic subgroup and is not virtually solvable, then the geometric structure on M is hyperbolic, and M may be either compact or non-compact.
\end{itemize}

Every closed $3$-manifold has a prime decomposition as a connected sum of prime $3$-manifolds,  that is, manifolds which do not admit nontrivial connect sum decompositions.  This decomposition is unique for orientable manifolds.   Thurston's conjecture (now the Geometrization Theorem) states that every oriented prime closed $3$-manifold can be cut along tori  so that the interior of each of the resulting $3$-manifold pieces has a geometric structure with finite volume.

\medskip

This discussion leads to the conclusion that the majority of the work in studying the geometry and topology of $3$-manifolds boils down to understanding hyperbolic manifolds and more generally groups of isometries of hyperbolic space and their quotients.

 \section{Basic concepts.}
Hyperbolic $3$-space with the metric $ds^{2}_{hyp}$ is denoted 
\begin{equation}
\mathbb{H}^3 = \{(x,y,z)\in \IR^3:z>0\},  \hskip10pt ds_{hyp}^{2}= \frac{dx^2+dy^2+dz^2}{z^2}.
\end{equation}
It has constant Riemannian sectional curvature equal to $-1$.  We identify its boundary $\partial\IH^3=\oC=\IC\cup\{\infty\}$ as the Riemann sphere.  Here of course,  $\IC$ denotes the complex plane.
The geodesics of hyperbolic three-space are subarcs of the lines $\{(z,t):z\in\IC, t>0\}$ or the subarcs of circles orthogonal to the boundary $\IC$.  These geodesics lie in the totally geodesic hyperbolic planes $\Pi_{z_0,r} = \{(z,t):z\in \IC, t>0, |z-z_0|^2+t^2 =  r^2\}$.  

\medskip

The action of an isometry $g:\IH^3\to\IH^3$ extends continuously to the boundary $\oC$  via the obvious action on geodesics.  Put $g_0:\oC\to\oC$ for the restriction $g|\oC$.  Since totally geodesic hyperplanes are mapped one to the other by an isometry,  the continuous boundary map $g_0$ must map lines and circles in $\oC$ among themselves.  It is an elementary exercise to see that $g_0$ is then a conformal map of $\oC$.  Such (orientable) maps are well known to be the linear fractional transformations.  The group of all such transformations is the M\"obius group   M\"ob$(\oC)$.  Consequently $g_0$ has the form
\[ g_0(z) = \frac{az+b}{cz+d}, \hskip20pt a,b,c,d\in \IC, \;\; ad-bc\neq 0 . \]
The converse is known as the Poincar\'e extension.  Write a linear fractional transformation $g_0$ as the finite composition of translations ($z\mapsto z+z_0$),  dilations ($z\mapsto\lambda z)$ and inversions ($z\mapsto 1/z$). For instance when $c\neq 0$, 
\[ z \mapsto   cz+d \mapsto   \frac{1 }{ cz+d} \mapsto \frac{1 }{cz+d}  + \frac{a}{bc-ad} =\frac{c}{bc-ad}\; \frac{az+b}{ cz+d} \mapsto  \frac{az+b}{cz+d}  = g_0(z) . \]
Since each of these factors obviously extends to $\IH^3$,  so too does $g_0$.  Evidently this association between extension and boundary values is unique and continuous.  The only slight problem is that the representation of the map $g_0=(az+b)/(cz+d)$ is not unique (multiply top and bottom by the same nonzero constant).  This is taken care of once we normalise so $ad-bc=1$ and note that as conformal mappings of the Riemann sphere 
\[ \frac{az+b}{cz+d}=\frac{-az-b}{-cz-d}. \]
In this way it is natural to identify the group of conformal mappings of $\oC$ as the group $\mathrm{PSL}(2,\IC)=\mathrm{SL}(2,\IC)/\pm$ by the mapping
\[g\in \mathrm{Isom}^+(\IH^3) \leftrightarrow g|\partial \IH^3 = \frac{az+b}{cz+d} \in \mbox{ M\"ob}(\oC) \leftrightarrow \pm \Big(\begin{array}{cc}a&b\\c& d\end{array}\Big)\in \mathrm{PSL}(2,\IC) . \]
In particular this association is a topological isomorphism of Lie groups.

\subsection{Classification of elements}
Let $\Gamma, \tilde{\Gamma} \subset \mbox{M\"ob}(\oC) \cong \mathrm{Isom}^+(\IH^3)$ be two subgroups.  We say that $\Gamma$ and $\tilde{\Gamma}$ are {\em conjugate} if there is a M\"obius transformation (or isometry) $f$ such that
\begin{equation}
\tilde{\Gamma} = f^{-1} \circ \Gamma \circ f
. \end{equation}
Most things of interest in the geometry of  discrete groups are conjugacy invariants and in effect the conjugating map merely changes coordinates by an isometry.  Thus we will almost always only be interested in information``up to conjugacy''.

The trace  is a natural conjugacy invariant of matrices and therefore on $\mathrm{SL}(2,\IC)$,  but unfortunately it is not well defined in $\mathrm{PSL}(2,\IC)$.  It is well defined up to sign however.  Therefore we make the following definition.  For  
\begin{equation}\label{lft} f(z)  = \frac{az+b}{cz+d}, \hskip10pt ad-bc=1
\end{equation} we set ${\rm tr}(f)=a+d$ and set
\begin{equation}
\beta(f) = {\rm tr}^2(f) - 4
. \end{equation}
Then $\beta:\mathrm{Isom}^+(\IH^3)\to \IC$ is a well defined conjugacy invariant of elements.  We note the following easy consequences of the same results for matrices.
\[ \beta(id) = 0, \;\;\; \beta(f) = \beta(f^{-1}), \;\;\; \beta(g\circ f\circ g^{-1}) = \beta(f), \;\;\; \beta(g\circ f) = \beta(f\circ g) . \]

Given a linear fractional transformation normalised as at (\ref{lft}) we can solve the fixed point equation $f(z) = z$, to find
\begin{eqnarray}\label{fix}
cz^2+(d-a)z-b & = & 0 \nonumber\\
z  & = & \frac{1}{2c} \Big((a-d) \pm \sqrt{(d-a)^2+4bc}\Big),  \hskip20pt c\neq 0 \\
z  & = & \infty, -b/(a-1/a), \hskip20pt c=0, a\neq 1  \\
z  & = & \infty, \hskip20pt c=0, a = d = 1 . \end{eqnarray}
The discriminant here is 
\[ (d-a)^2-4bc = a^2+d^2-2ad+4bc = (a+d)^2 - 4 = \beta(f) . \]
If a linear fractional transformation $f$ has two fixed points, say $z_0,z_1$,  then we can choose another linear fractional transformation $g$ with $g(z_0)=0$ and $g(z_1)=\infty$, (if there is one fixed point $z_0$,  then choose $g(z_0)=\infty$ and follow these arguments).  Now $h(z) = (g\circ f \circ g^{-1})(z)$ fixes $0$ and $\infty$.  Thus $h(z)=\lambda z$,  $\lambda\in\IC\setminus\{0\}$.
There is an obvious dichotomy in the dynamics of the map $h$ acting on $\oC$ here.  If $|\lambda|=1$,  then the group $\overline{\langle h \rangle}$ is a compact group of rotations.  This leads to the following definitions:
\begin{itemize}
\item {\bf Elliptic:} $f$ is {\em elliptic} if $\beta(f)\in [-4,0)$.  In this case $f$ is conjugate to a rotation $z\mapsto \zeta z$, $|\zeta|=1$.
\item {\bf Loxodromic:} $f$ is {\em loxodromic} if $\beta(f)\not\in [-4,0]$.   In this case $f$ is conjugate to a dilation $z\mapsto \lambda z$, $0< |\lambda| \neq 1$.
\item {\bf Parabolic:} $f$ is {\em parabolic} if $\beta(f)=0$.   In this case $f$ is conjugate to the translation $z\mapsto z+1$.
\end{itemize}
There are a couple of values for $\beta(f)$ it is worth observing.  When $\beta(f)\in[-4,0)$,  $f$ is conjugate to a rotation. This rotation may have a certain period - the order.
\begin{itemize}
\item If $\beta(f)=-4$, then $f$ is elliptic of order two.
\item If $\beta(f)=-3$, then $f$ is elliptic of order three.
\item If $\beta(f)=-2$, then $f$ is elliptic of order four.
\item If $\beta(f)=-1$, then $f$ is elliptic of order six.
\item If $\beta(f)=-4\sin^2(p\pi/n)$ with $(p,n)=1$, then $f$ is elliptic of order $n$.
\end{itemize}
We further say that $f$ is a {\em primitive elliptic} of order $n$ if $\beta(f)=-4\sin^2(\pi/n)$,  so both $f$ and $f^{-1}$ are simultaneously primitive or not.
\subsection{The geometry of elements}
We now discuss the action of each of the elements on hyperbolic space.  If $\beta(f)\neq 0$,  then $f$ has two fixed points, say $z_1,z_2$.  The hyperbolic line whose endpoints are $z_1$ and $z_2$ is called the {\em axis} of $f$,   denoted $\ax(f)$.  

The axis of $f$ is clearly left setwise invariant under the map $f$,  $f(\ax(f))=\ax(f)$ while under conjugacy it is elementry to see that
\[ \ax(g\circ f\circ g^{-1}) = g(\ax(f)) .\]  
In fact,  up to conjugacy we could assume $z_1=0$ and $z_2=\infty$.  Then the axis would be the $t$-axis, $\{(0,0,t):t>0\}$, and we would have $f$ conjugate to $g:z\mapsto \lambda^2 z$.  We put $\lambda^2=re^{i\theta}$,  $r\geq 1$ (note that $g$ and $g^{-1}$ are conjugate,  so we can assume $r\geq 1$) and see that with $\j=(0,0,1)\in\IH^3$, 
\[  g(\j) = r \j,  \hskip10pt \tau(g) = \rho_{\IH^3}(\j,g(\j)) = \int_{1}^{r}\frac{dt}{t} = \log(r) \]
so $g$ translates along its axis by a hyperbolic distance $\tau = \tau(g)$.  This of course must be the same for $f$ acting along its axis.  We call this distance the {\em translation length} and denote it $\tau(f)$.  As $g$ translates along the $t$-axis it also rotates around it by some angle (recall loxdromic means screwlike).  This angle,  $\theta=2\arg(\lambda)$,  is called the holonomy of $g$,  denoted $\theta(g)$.  Of course this is a conjugacy invariant and $\theta(g)=\theta(f)$.  Once we notice
\[ g(z) = \lambda^2  z = r e^{i\theta} z \leftrightarrow \Big( \begin{array}{cc}\sqrt{r} e^{i\theta/2} & 0 \\ 0 &   e^{-i\theta/2}/\sqrt{r} \end{array}\Big). \]
Therefore
\[ \beta(f) =\beta(g) = (\sqrt{r} e^{i\theta/2} + e^{-i\theta/2}/\sqrt{r})^2 -4 = -4 \sinh^2\Big(\frac{\tau}{2}+i\frac{\theta}{2}\Big) . \]
We record this as a lemma.
\begin{lemma}  Let $f\in \mathrm{PSL}(2,\IC)$ be a hyperbolic isometry.  Let $\tau_f=\tau(f)$ be the translation length of $f$ along its axis,  and $\theta_f=\theta(f)$ the holonomy of $f$ about its axis.  Then
\begin{equation}
\beta(f) = 4 \sinh^2\Big(\frac{\tau_f}{2}+i\frac{\theta_f}{2}\Big) = 2\left( \cosh(\tau_f+i\theta_f)-1\right) 
. \end{equation}
\end{lemma}
Of course knowing $\beta$ identifies both $\tau(f)$ and $\theta(f)$ (modulo $2\pi$).
\section{Two-generator groups} 
We
now define, for $f,g \in \mathrm{Isom}^+(\IH)$,
\begin{equation}
 \gamma(f,g)=\tr[f,g]-2
 . \end{equation}
Note that $\tr[f,g] = \tr( fgf^{-1}  g^{-1})$ is defined unambiguously for $f,g\in \mathrm{PSL}(2,\IC)$ because of the even number of appearances of both $f$ and $g$.

\medskip

The parameters defined above conveniently encode various other
geometric quantities.  One of the most important is the  
complex hyperbolic distance $\delta+i\theta$ between the axes of two non-parabolic M\"obius transformations $f$ and $g$.  Here $\delta$ is the hyperbolic distance between the lines $\ax(f)$ and $\ax(g)$ and $\theta$ is the dihedral angle between  two hyperplanes - if $\alpha$ denotes the common perpendicular between $\ax(f)$ and $\ax(g)$,  then  one plane contains $\ax(f)$ and $\alpha$ while the other contains $\ax(g)$ and $\alpha$.  The easiest way to see the angle $\theta$  is to use a conjugacy to arrange things so that the common perpendicular $\alpha$ lies on the $t$-axis.  Then the angle $\theta$ is simply the angle between the vertical projections to $\oC=\partial\IH^3$ of $\ax(f) $ and $\ax(g)$ at the origin. Alternatively,  $\theta$ is the holonomy of the transformation whose axis contains the common perpendicular and moves $\ax(f)$ to $\ax(g)$.

If $f$ and $g$ are each elliptic or loxodromic, with translation
lengths $\tau_f$ and
$\tau_g$ respectively and holonomies $\eta_f$ and $\eta_g$
respectively, then we have the following identities,  \cite{GM1}.
   \begin{eqnarray}
\label{betageom}
\beta(f) & = & 4 \sinh^2\left(\frac{\tau_f+i\eta_f}{2} \right) \\
\beta(g) & = & 4 \sinh^2\left(\frac{\tau_g+i\eta_g}{2} \right) \\
\label{gammageom}
\gamma(f,g) & = & \frac{\beta(f) \beta(g)}{4} \sinh^2(\delta+i\theta)
\end{eqnarray}
where $\delta+i\theta$ is the complex distance between the axes of $f$ and $g$.  From (\ref{betageom})-(\ref{gammageom}) we derive the following useful formulas.
\begin{eqnarray}\label{cosh(tau)}
\cosh(\tau_f) &=& \frac{|\beta(f)+4|+|\beta(f)|}{4}\\
   \cos(\eta_f) &=& \frac{|\beta(f)+4|-|\beta(f)|}{4} \\
\label{cosh2r}
\cosh(2\delta)&=&
\left|\frac{4\gamma(f,g)}{\beta(f)\beta(g)}+1\right|
+\left|\frac{4\gamma(f,g)}{\beta(f)\beta(g)}\right|\\
\label{cos2theta}
\cos(2\phi)&=&
\left|\frac{4\gamma(f,g)}{\beta(f)\beta(g)}+1\right|
-\left|\frac{4\gamma(f,g)}{\beta(f)\beta(g)}\right|.
\end{eqnarray}
   We are often concerned with the case where one of the
isometries,  say $g$,  is of order 2, in which case $\beta(g)=-4$, and
(\ref{cosh2r}) and (\ref{cos2theta}) take the simpler form
\begin{eqnarray}
\label{order2cosh2r}
\cosh(2\delta)&=&
|1-\gamma(f,g)/\beta(f)|+|\gamma(f,g)/\beta(f)|\\ 
\label{order2cos2theta}
\cos(2\phi)&=&
|1-\gamma(f,g)/\beta(f)|-|\gamma(f,g)/\beta(f)|.
\end{eqnarray}
Notice that for fixed $\beta(f)\in\IC$ and fixed $\cosh(2\delta)$ at (\ref{order2cosh2r}),  the set of possible values for $\gamma(f,g)$  form an ellipse,  while for  fixed $\cosh(2\theta)$ at (\ref{order2cos2theta}) we get hyperbola. Thus $\delta$ and $\theta$ give very appealing geometric orthogonal coordinates on $\IC\setminus [\beta,0]$.

 \subsection{Complex parameters for two-generator groups}\label{Riley}

The discussion above allows us to  view the space of all two-generator Kleinian groups (up to  
conjugacy) as a subset of the three complex
dimensional space $\IC^3$ via the map
\begin{equation} \langle f,g \rangle \to
(\gamma(f,g),\beta(f),\beta(g)). 
\end{equation}
The next theorem is quite elementary,  the proof simply amounting to solving equations after normalising by conjugacy.  However the restriction $\gamma\neq 0$ cannot be removed.
\begin{theorem}\label{param}
If $\gamma\neq 0$,  then the complex triple $(\gamma,\beta,\beta^\prime)$ uniquely determines a two-generator M\"obius group $\langle f,g\rangle$ up to conjugacy.
\end{theorem}
It is worth recording the following elementary but useful fact here concerning the exception $\gamma\neq 0$ in Theorem \ref{param},
\begin{lemma}  If $f$ and $g$ are M\"obius transformations and if $\gamma(f,g)=0$,  then $f$ and $g$ share a common fixed point in the Riemann sphere $\oC$.  In particular,  if $\langle f,g\rangle$ is discrete and if $\gamma(f,g)=0$,  then $\langle f,g\rangle$  is virtually abelian (elementary).
\end{lemma}
Now a fundamental problem is to describe the space of all points in
$\IC^3$  corresponding to
discrete two-generator groups. Here we will actually describe
portions of this space, usually slices of codimension one or two.  For instance,  an important case will be the following.  We fix $p,q$ and let $f$ and $g$ be primitive elliptics of order $p$ and $q$ respectively.  Then any $\gamma \in \IC\setminus \{ 0\}$
uniquely determines the conjugacy class of   a two-generator group $\G = \langle  f , g \rangle $.   For fixed $p$ and $q$ it is an elementary consequence of a theorem of J\o rgensen \cite{Jorgensen}, which we outline below, that the  set of all such $\gamma$ values which correspond to discrete non-elementary groups is closed,  and computer generated pictures suggest 
that it is highly fractal in nature - for instance the well known Riley slice, 
corresponding to two parabolic generators, with our parameters correspond to groups generated by an element of order $p=2$ and a parabolic,  $q=\infty$.  The boundary of the space of such groups freely generated by two such elements is illustrated below.

\scalebox{0.5}{\includegraphics[viewport= -5 460 580 790]{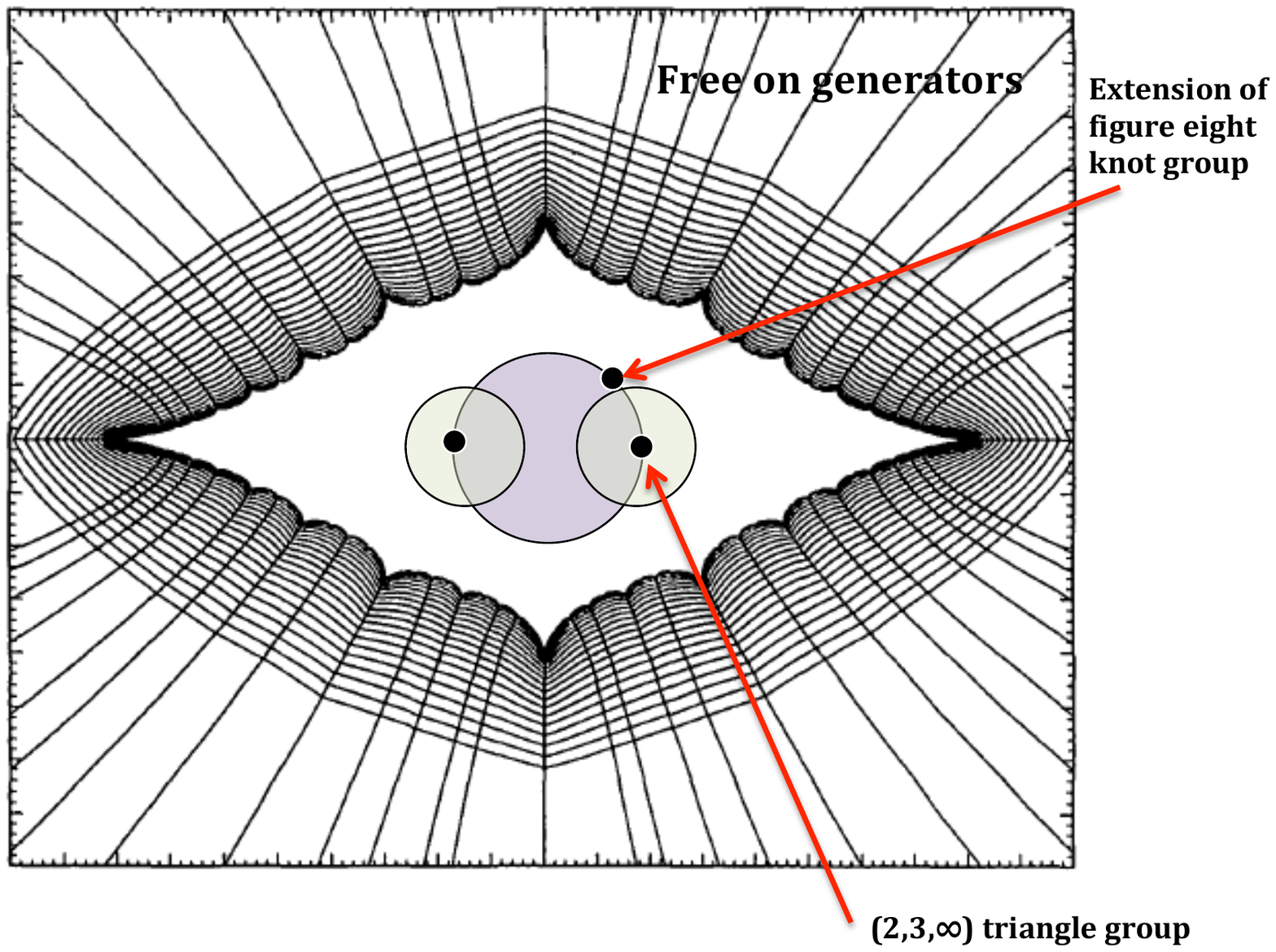}}
\bigskip
\begin{center} ``The Riley Slice''
\end{center} 

\bigskip 
This then is a picture of the slice $\beta(f)=0$,  $\beta(g)=-4$,  $\gamma=\gamma(f,g)\in\IC$ a one complex dimensional slice in the space $\IC^3$ of all two-generator M\"obius groups.  For $\gamma$ outside the bounded region,  the group $\G$ is freely generated by $f$ and $g$ (and $(\oC\setminus L(\G))/\G$ is a marked sphere,  two cone points of order $2$ and two punctures).  

Inside the bounded region we will find isolated $\gamma$ values which correspond to $\gamma$ values for groups which are not freely generated by $f$ and $g$.  For instance the  value $\gamma=(1+i\sqrt{3})/2$ for the $\IZ^2$-extension of the figure eight knot complement  is identified.  We have also  illustrated the unit disk in the centre,  it is a consequence of the Shimitzu-Leutbecher inequality that there are no Kleinian groups with $|\gamma|<1$ at all.  Our polynomial trace identities (discussed in \S 7) extend this to further inequalities giving the two additional disks $|\gamma\pm 1|< r_0$,  where $r_{0}^{2}(r_0+1)=1$, and these additional disks also contain no Kleinian groups except for those $\gamma$ values at their  isolated centres $\gamma=\pm 1$ corresponding to the $(2,3,\infty)$ triangle groups.  Inside this region there are additionally infinitely many $\IZ_2$-extensions, the two-bridge knot and  two-bridge links groups as well as the $(2,p,\infty)$,  $p\geq 3$ triangle groups and other web groups.  Consider for a moment an example; $f$ is parabolic and $g$  elliptic of order $2$ generating the $(2,p,\infty)$ triangle group. Then $h=fg$ is primitive elliptic of order $p$.  The axis of $g$ and $h$ are both perpendicular to a hyperplane on which a $(2,p,\infty)$ hyperbolic triangle is formed,  so the angle between these axes is $\theta(g,h)=\theta =0$.  Next,  a little hyperbolic trigonometry tells us that the distance $\delta = \delta(g,h)$ between the axis of $g$ and $h$ satisfies 
\[ \cosh(\delta) = \frac{1}{\sin(\pi/p)} . \]
Then from (\ref{gammageom}) we see that
\begin{eqnarray*}
\gamma(h,g)  & =   & \frac{\beta(h) \beta(g)}{4} \sinh^2(\delta+i\theta) = - \beta(h)  \sinh^2(\delta)  \\
& = &  4\sin^2(\pi/p) (\cosh^2(\delta)-1) =  4-4 \sin^2(\pi/p) 
. \end{eqnarray*}
Then
\[ \gamma(f,g)=\gamma(fg,g)=\gamma(h,g) . \]
This gives the infinite sequence of points 
\[ \left\{\big(  4-4 \sin^2(\pi/p),  0, -4\big) :p\geq 3 \right\}\]
 lying in the bounded region of the space illustrated above.  Notice that as $p\to \infty$,  $4-4 \sin^2(\pi/p)\to 4$,  a point lying in the boundary of this space,  actually corresponding to the $(2,\infty,\infty)$-triangle group which is actually free on these two-generators.   

 We will see another of these spaces in a bit,  when the generators are order two and three.  Before that we  have to develop a bit more theory.

\subsection{Projection to $(\gamma,\beta,-4)$}

There is an important projection from
the three complex dimensional space of discrete groups $(\gamma,\beta,\beta^\prime)$,
to the two complex dimensional slice
$\beta^\prime=-4$  which preserves discreteness.  This projection arises 
roughly because we can identify a closely related group which has one generator of order two.  This goes as follows.  First it is easy to show that if $f_1$ and $f_2$ have the same trace,  then there is an elliptic $\phi_1$ of order two  such that
\[ f_1=\phi^{-1}_{1} \circ f_2 \circ \phi_{1} . \]  
In fact there are two such elliptics of order two,   $\psi_i$,   $i=1,2$ ,  which are two rotations of order $2$  about the  bisector of the common perpendicular between the axes,  see \cite{GM3} for details.  Then  $g\circ f\circ g^{-1} = \phi_i \circ f\circ \phi_{i}^{-1}$,  $i=1,2$ and it is easy to see that $\langle f,g\circ f\circ g^{-1} \rangle $ is an index-two subgroup of $\langle f,\phi_i\rangle$,  $i=1,2$.  We have not yet discussed discreteness closely,  but that one group is of finite-index in another will imply that  if one group is discrete,  then so is the other.  Next we make a few simple calculations to see the following.\begin{eqnarray}
\mbox{Group} && \mbox{parameters} \nonumber \\
\langle f, g \rangle & \leftrightarrow & (\gamma,\beta,\beta^\prime)\nonumber  \\
\langle f, gfg^{-1} \rangle &  \leftrightarrow  & (\gamma(\gamma-\beta),\beta,\beta) \label{above} \\
\cap \;\;\;\;& & \nonumber  \\
\langle f, \phi_i \rangle &  \leftrightarrow  & (\gamma,\beta,-4), \;\;{\rm and} \;\; (\beta-\gamma,\beta,-4). \nonumber 
\end{eqnarray}
At the end of the day we have the following lemma which turns out to be very useful.  As an example,  if we seek collaring theorems then we need to consider a loxodromic axis of $f$ and its nearest translate,  say under $g$.  This then gives us a group with two-generators of the same trace $f$ and $gfg^{-1}$.  The two groups found by the above constructions are discrete and give exactly the same collaring radius for $f$.  So,  with $\gamma=\gamma(f,g)$ we can go from $(\gamma,\beta(f),\beta(g))$ to $(\gamma,\beta(f),-4)$ and still retain the geometric information we might want.  

\begin{lemma}\label{beta=-4}  Suppose $(\gamma,\beta,\beta^\prime) $
are the parameters of a discrete group.  Then so are
$(\gamma,\beta,-4)$.
\begin{equation} (\gamma,\beta,\beta^\prime) \mbox{ discrete}
\Rightarrow (\gamma,\beta,-4) \mbox{
discrete} . \end{equation}
Further, away from a small finite set of exceptional parameters this
projection preserves the
property of being non-elementary as well.
\end{lemma}
This exceptional set  occurs only in
the presence of finite spherical and Euclidean triangle subgroups and motivates us to identify all the possible parameters for these groups.

\subsection{Symmetries of the parameter space}

There are a few natural symmetries of the parameter space which will further cut down the region we have to investigate.  The first obvious one is complex conjugation.  There are others and  they can be deduced from the following two  trace identities that we use quite a lot.  The first is well known,
\begin{equation}
\tr(fh)+\tr(fh^{-1}) = \tr(f)\tr(h)
\end{equation}
and the second is the the  Fricke identity
\begin{equation}\label{ti2}
\tr[f,h]=\tr^2(f)+\tr^2(g)+\tr^2(fg)-\tr(f)\tr(g)\tr(fg)-2
. \end{equation}
This last one (\ref{ti2}) reads as
\begin{equation}
\gamma(f,g)=\beta(f)+\beta(g)+\beta(fg)-\tr(f)\tr(g)\tr(fg) + 8 
\end{equation}
and so if $g=hf^{-1}h^{-1}$,  then
\begin{eqnarray*}
\gamma(f,hfh^{-1})& = & 2\beta(f)+((\gamma(f,h)+2)^2-4)-(\beta(f)+4)(\gamma(f,h)+2)+ 8\\
& = & \gamma(f,h)(\gamma(f,h)-\beta(f))
. \end{eqnarray*}
We have already used this map $\gamma\mapsto\gamma(\gamma-\beta)$   above at (\ref{above}).
Further,  if $g$ has order two,  then $\tr(g)=0$ and (\ref{ti2}) gives $\tr[f,g]-2=\tr^2(f) + \tr^2(fg) -4$ and so
\begin{equation} 
\gamma(f,g)=\beta(f) + \beta(fg) + 4
. \end{equation}
This tells us that when  $g$ has order two
\[\beta(fg)=\gamma(f,g)-\beta(f)-4 . \]
We then deduce the following result.
\begin{theorem}  If $(\gamma,\beta,-4)$ are the parameters of a discrete group,  then so too are
\begin{itemize}
\item $(\gamma,\beta,-4)$,  corresponding to $\langle f,g\rangle$
\item $(\beta - \gamma,\beta,-4)$, corresponding to $\langle f,\varphi g\rangle$
\item $(\gamma,\gamma-\beta-4,-4)$, corresponding to $\langle fg, g\rangle$
\item $(-\beta-4,\gamma-\beta-4,-4)$, corresponding to $\langle fg, \varphi g\rangle$
\item $(-\beta-4,\gamma-\beta-4,-4)$, corresponding to $\langle fg, \varphi g\rangle$ 
\end{itemize}
\end{theorem}
Notice at all these symmetries have preserved $\beta(g)=-4$.

There is also another action on the space of parameters of discrete groups. This is the multiplicative action by shifted Chebychev polynomials on these slices. These results are more carefully discussed in \cite{GMCheby}.
\begin{theorem}  Consider $\langle f^n,g \rangle\to \langle f^n,g \rangle$.  Then 
\[ \beta(f)=-4\sinh^2\Big(\frac{\tau+i\theta}{2}\Big), \hskip10pt \beta(f^n) =-4\sinh^2\Big(n\frac{\tau+i\theta}{2}\Big) \]
and
\begin{equation}
(\gamma(f^n,g),\beta(f^n) ) = \frac{\beta(f^n)}{\beta(f)} \; (\gamma,\beta) 
. \end{equation}
That is 
\[(\gamma,\beta)\mapsto T_n(\beta) (\gamma,\beta) \]
where $T_n$ is a shifted Chebychev polynomial,
\begin{equation}
z T_n(z) = 2U_n\big(\frac{z}{2}+1\big)-2
\end{equation}
and $U_n$ is the usual Chebychev polynomial,  $U_n(\cosh(z))=\cosh(nz)$.
\end{theorem}

\subsection{Kleinian groups}

A {\em  Kleinian group} is a discrete subgroup $\Gamma$ of the group of orientation-preserving isometries of hyperbolic $3$-space  which is not virtually abelian.  Let us work through this definition for a moment.  According to our discussion above there are at least three different ways of thinking about these groups.  As subgroups of the group of isometries of hyperbolic space $\mathrm{Isom}^+(\IH^3)$,  as subgroups of the M\"obius group M\"ob($\IC$) and also as subgroups of $\mathrm{PSL}(2,\IC)$.  Each of these groups has its own topology,  however the topological isomorphism we outlined tells us that the concept of discreteness is the same.  To fix ideas we say that a subgroup $\G\subset$M\"{o}b($\oC$) is discrete if the identity element is isolated in $\G$,  that is there is no sequence $\{g_j\}_{j=1}^{\infty}\subset\Gamma$ with 
\begin{equation}
g_j\to{\rm identity}\;\;\;\; \mbox{as $j\to\infty$}
. \end{equation}
Notice that since $\G$ is a group,  there can be no convergent sequences whatsoever.  If $g_j\to g$,  with $g$ not necessarily in $\G$ we would still have  $g_{j+1}\circ g_{j}^{-1}\to {\rm identity}$ in $\G$.  Thus,  in a discrete group the identity is isolated.  What is a little surprising initially is that this isolation can be quantified in the case of negative curvature - indeed the Margulis constant is an example of this and there are many others.
\subsubsection{The convergence properties of M\"obius transformations}
In terms of the relationship between the topology induced from the group $\Gamma$ being viewed as a subgroup of the matrix group $\mathrm{PSL}(2,\IC)$ and as a group of linear fractional transformations of the Riemann sphere we note the following theorem.
\begin{theorem}\label{convergence}  Let $\{f_j\}_{j=1}^{\infty}$ be a sequence of M\"obius transformations of the Riemann sphere.  Then there is a subsequence $\{f_{j_k}\}_{k=1}^{\infty}$ such that either
\begin{enumerate}
\item there is a linear fractional transformation $f:\oC\to\oC$ such that  $f_{j_k}\to f$ and $f_{j_k}^{-1}\to f^{-1}$ uniformly in $\oC$ as $k\to\infty$,  or
\item there are two points $z^*$ and $w^*$ such that 
\begin{itemize} 
\item $f_{j_k}\to z^*$ locally uniformly in $\oC\setminus\{w^*\}$, and
\item  $f_{j_k}^{-1}\to w^*$ locally uniformly in $\oC\setminus\{z^*\}$ as $k\to\infty$.
\end{itemize}
\end{enumerate}
\end{theorem}
This topological theorem actually is all that is needed to develop much of the theory of discrete groups of conformal transformation of the Riemann sphere $\oC$ (and indeed higher dimensional spheres) and leads to the notion of convergence groups \cite{GMDQCG}.  Of course condition (1) will not occur in a discrete group.  One can find a proof for Theorem \ref{convergence} simply by studying the way a sequence of matrices in $\mathrm{SL}(2,\IC)$ can degenerate.
\subsubsection {Limit sets,  elementary and non-elementary groups.}
The dynamical study of a discrete group of M\"obius transformations acting on the Riemann sphere is a deep and interesting subject.  Unfortunately it is largely irrelevant to our considerations here but we will make a very brief digression so as to introduce a few concepts.  Basically we want to get a condition to guarantee that a group is discrete and free on its generators.  So we will need the combination theorems and a little knowledge of how a Kleinian group acts on the Riemann sphere $\oC$.

Thus let $\G$ be a discrete group of M\"obius transformations acting on $\oC$.  Basically the limit set of $\G$ is the set of accumulation points of a generic orbit. If we consider the group acting on $\IH^3$,  and if we put $\j=(0,0,1)$, then
\[ L(\G) = \mbox{the set of accumulation points of the orbit $\G(\j)=\{g(\j):g\in \G\}$} . \]
Notice that discreteness will imply that the orbit can only accumulate on the boundary $\oC$.  Another definition of the limit set would be the set of all points $z^*$ and $w^*$ identified in Theorem \ref{convergence}.  It is easy to see that the limit set is closed and $\G$ invariant;  if $g\in\G$,  then $g(L(\G))=L(\G)$. The following dichotomy is fairly well known and also not too hard to prove.
\begin{theorem}  Let $\G$ be a discrete group of M\"obius transformations.  Then either
\begin{enumerate}
\item $ \# L(\G) \leq 2 $,  or
\item $L(\G)$ is perfect (a closed and uncountable set,  with every point a point of accumulation).
\end{enumerate}
\end{theorem}
This dichotomy reflects the elementary/nonelementary dichotomy that we have already discussed.  It also provides an alternative definition of a non-elementary group.
\begin{lemma} A discrete group of M\"obius transformations of $\oC$ is elementary (virtually abelian) if and only if $\#L(\G)\leq 2$.
\end{lemma}

Next,  the ordinary set of a Kleinian group is the set $O(\G)=\oC\setminus L(\G)$.  The group acts ``chaotically'' on $L(\G)$ and nicely on $O(\G)$.  The reader can imagine that there is a significant literature in making these statements precise.  In particular Theorem \ref{convergence}  implies that $\G$ acts properly discontinuously on $O(\G)$ and so $O(\G)/\G$ is a Riemann surface.  A famous theorem of Ahlfors shows this surface to be of finite topological type if $\G$ is finitely generated - Ahlfors finiteness theorem \cite{Ahlfors}.  While more recent results imply that $L(\G)$ has Lebesgue measure $0$ for finitely generated Kleinian groups.  In any case,  that $\G$ acts properly discontinuously on $O(\G)$ tells us that there is a {\em fundamental domain} for this action.  That is a connected open set $U$ such that if $g,h\in \G$,  $g\neq h$,  then $g(U)\cap h(U) = \emptyset$ and
\[ O(\G) = \bigcup_{g\in \G} g(\overline{U}) . \]
So the images of $U$ tessellate the ordinary set.

\subsubsection {Schottky groups and a combination theorem.}
The first simple examples one meets of discrete nonelementary groups - perhaps after meeting Fuchsian groups - are the Schottky groups.  Here the limit set will be a Cantor set in $\oC$ and every element in the group will be loxodromic. Let $ f(z)= \frac{az+b}{cz+d}$.  There is a natural way to construct these groups through pairings of isometric circles.  The {\em isometric circles} of $f$ are
\[ \{z\in \IC: |cz+d|=1\},  \{z\in \IC: |cz-a|=1\} \]
and $f$ pairs these circles.  Of course the isometric circles of $f$ and $f^{-1}$ are the same.  Let $U_1=\{z\in \IC: |cz+d|<1\}$ and $U_2=\{z\in \IC: |cz-a|<1\}$.  The terminology arises as $|f'(z)|\equiv 1$ on these circles.  Next
\begin{equation} \label{fpair}
f(U_1)=\oC\setminus \overline{U_2},  \hskip10pt f^{-1}(U_2)=\oC\setminus \overline{U_1} 
. \end{equation}
Given another M\"obius transformation $g$ we can find a similar pair of disks $V_1$ and $V_2$ bounded by the isometric circles of $g$ with 
\begin{equation} \label{gpair}
g(V_1)=\oC\setminus \overline{V_2},  \hskip10pt g^{-1}(V_2)=\oC\setminus \overline{V_1} 
. \end{equation}
The following theorem is often referred to as the ``ping-pong'' lemma,  a name suggested by its proof.  The picture is illustrated below.

\scalebox{0.5}{\includegraphics[viewport= -40 400 480 750]{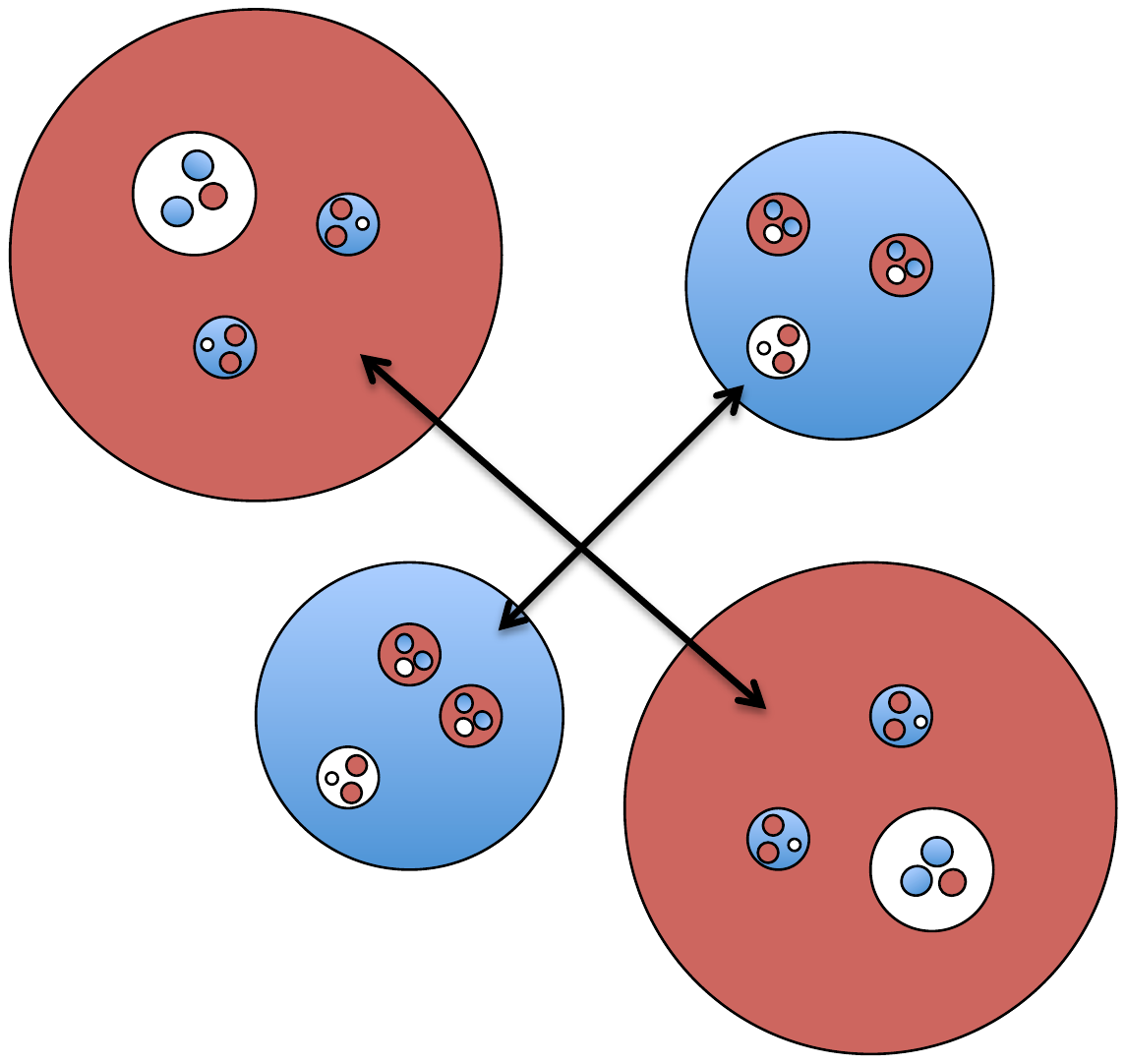}}

\medskip
Paired isometric circles giving nested sequences of disks whose intersection is the Cantor set limit set.

\bigskip

\begin{theorem}Suppose the four disks $U_1,U_2,V_1, V_2$ are all disjoint.  Then the group $\langle f,g\rangle$ is Kleinian and has limit set a Cantor set.  The group is algebraically isomorphic to the free group on two-generators.
\end{theorem}
\noindent{\bf Sketch of proof.} Let $w=f^{p_1}\circ g^{q_1} \circ \cdots f^{p_m}\circ g^{q_m}$ with $p_i,q_i\neq 0$ (possibly allowing $p_1=0$ or $q_m=0$)  be a word in the group $\langle f,g \rangle$.  Let $z\in \Omega = \oC\setminus \{U_1,U_2,V_1, V_2\}$.  Then when we evaluate $w$ on $z$ we see that  every occurrence of $f^{\pm1} $ or $g^{\pm 1}$ buries the image deeper and deeper in the nested sequence of circles.   In particular $w(\Omega)\cap \Omega=0$ where upon the group is discrete and free,  since no $w=\mathrm{identity}$. \hfill $\Box $

\medskip

This result is part of a much more general theory about the combinations of discrete groups.  This includes $HNN$ extensions,  free product and other amalgamations. This is very well explained in B. Maskit's book,  \cite{Maskit}.  We record the obvious generalisation here.  

\begin{theorem}[Klein Combination Theorem] \label{KCT} Let $\G_1$ and $\G_2$ be Kleinian groups acting on $\oC$ with fundamental domains $U_1$ and $U_2$ respectively.  Suppose that
\[ \oC\setminus U_2 \subset U_1,  \hskip10pt \oC\setminus U_1 \subset U_2 . \]
Then the group $\langle \G_1,\G_2\rangle$ generated by $\G_1$ and $\G_2$ is discrete and isomorphic to the free product $\G_1* \G_2$.
\end{theorem}

We want to use this theorem to describe various parameter spaces of discrete groups with elliptic generators.  If $f$ is elliptic with $f(\infty)\neq\infty$,  then the two isometric circles of $f$ intersect at the fixed points of $f$ and the exterior of the two bounded disks is a fundamental domain for $f$,  as is the intersection of these two disks.

Suppose $\G$ is a   group generated by two primitive elliptic elements $f,g$ of
orders $p,q $.   We  normalise a choice of matrix representatives for the elements
$f,g$ as
\begin{equation}
f \sim \left(\begin{array}{cc} \cos \pi/p & i \sin \pi/p \\ i \sin \pi/p &
\cos \pi / p
\end{array}\right), \quad g \sim \left(\begin{array}{cc} \cos \pi/q & i \omega \sin
\pi/q \\
i  \omega ^{-1} \sin \pi/q & \cos \pi/q  \end{array}\right)
. \end{equation}
Here $ \omega $ is  a complex parameter which by further conjugation, if necessary,
we can assume lies in the unit disk,  $|\omega| < 1$.
Next,  according to our discussion about fundamental domains for $\langle f\rangle$ and $\langle g \rangle$ and the Klein combination theorem \ref{KCT},  if the isometric circles of $g$ lie inside the
region bounded by the intersection of the isometric circles for $f$, then the group $\langle f,g \rangle $ is discrete and isomorphic to a free product,  $\langle f\rangle*\langle g \rangle\cong \IZ_p*\IZ_q$.  Thus, with the
normalisation
above, if the inequality
\begin{equation}
\mid \sin(\pi/q) \cos(\pi/p) \pm  \omega  \cos(\pi/q) \sin(\pi/p) \mid + \mid  \omega t
\mid\sin(\pi/p) \leq \sin(\pi/q)
\end{equation}
holds for both choices of sign, then we have $G \cong  \IZ_p*\IZ_q$.
These in particular hold if
\begin{equation}\label{omegabound}
\mid  \omega  \mid \leq \frac{\sin(\pi/p) \sin(\pi/q)}{(1 + \cos(\pi/p)) ( 1 +
\cos(\pi/q))}.
\end{equation}
We then calculate that 
\begin{equation}
\gamma(f,g) = (\omega-1/\omega)^2 \sin^2(\pi/p)\sin^2(\pi/q)
. \end{equation}
To get an idea,  a picture of the configurations of isometric circles for discrete groups freely generated by two elliptics is illustrated below.

\scalebox{0.33}{\includegraphics[viewport= -20 350 660 770]{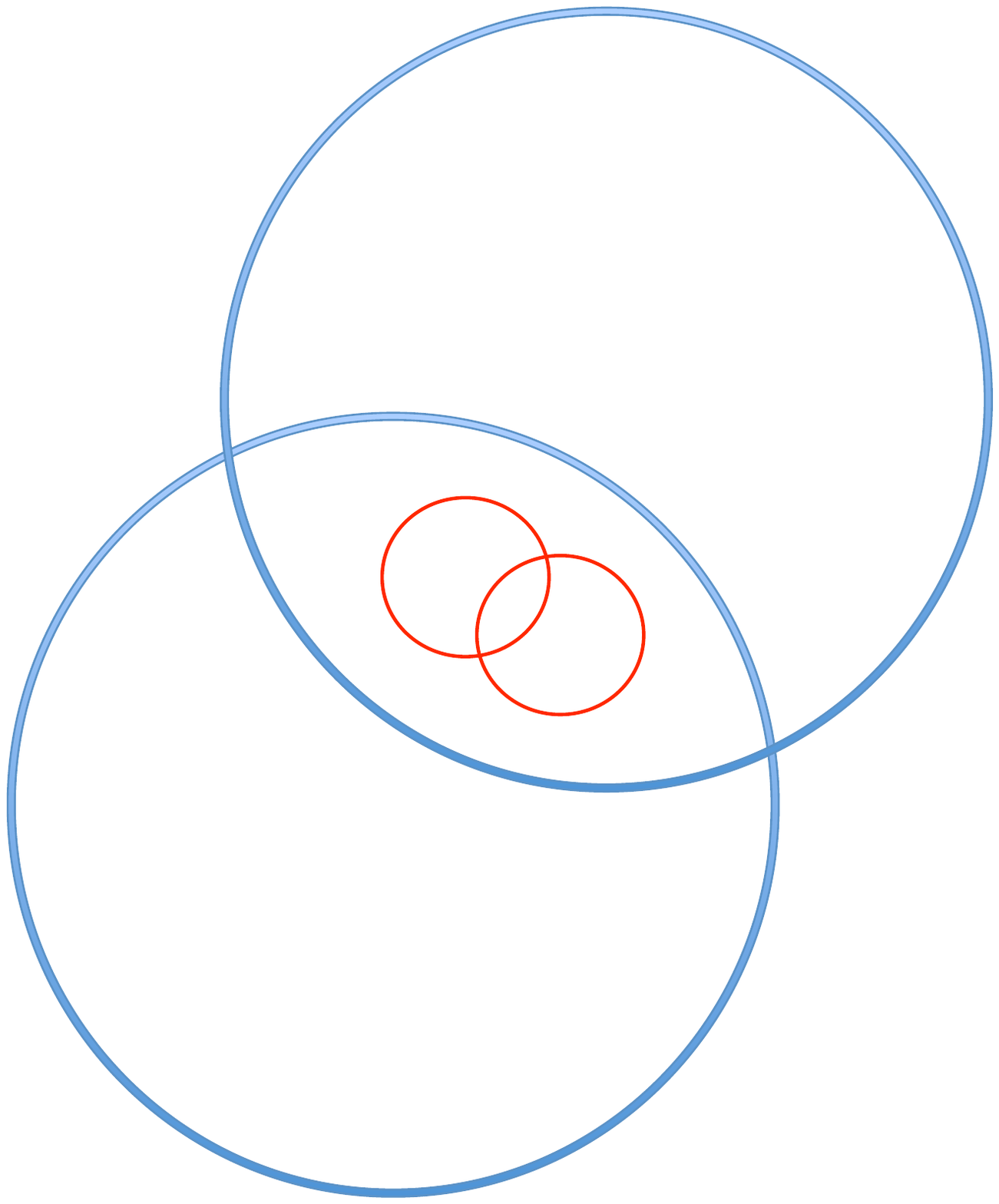}}

\subsection{The free part.}

Now,  putting these together we obtain the following two results,  \cite{GMM2elliptics}.
 
 \begin{theorem}\label{free} Put
\begin{equation}
\lambda_{p,q} = 4(\cos\pi/p  + \cos\pi/q)^2 + 4(\cos \pi/p \;  \cos\pi/q + 1)^2
 . \end{equation}
Let $\G$ be a M\"obius group generate by primitive elliptics $f$ and $g$ of orders $p$ and $q$ respectively,  $2\leq p,q \leq \infty$ and not both equal to $2$.  Then with $\gamma=\gamma(f,g)$ the parameters for $\G$ are 
\[  (\gamma, -4 \sin^2(\pi/p), -4 \sin^2(\pi/q)). \]
Further, if $\gamma$ lies outside of the open ellipse defined by the equation
 \begin{equation}
\{z: |z| + |z + 4 \sin^2(\pi/p) \sin^2(\pi/q)| < \lambda_{p,q}\} ,
\end{equation} then $\G$ is discrete and isomorphic to the free product of cyclics.
This result is sharp.  
\end{theorem}
In particular this theorem shows that $\G$ cannot be isomorphic to a finite covolume lattice of $\mathrm{Isom}^+(\IH^3)$ when $\gamma$ lies outside this ellipse. This is because in three-dimensions  a group of finite covolume has Euler characteristic  zero and  so $\G$ cannot split into a free product of cyclic  groups.  We note a case we discussed earlier and will use later. If $n=\infty$ and $m=2$ we have the ellipse actually being the disk $\ID(0,4)$.  This is our Riley slice picture.  We identified the points $\gamma=4-4\sin^2(\pi/p)$ and by symmetry  $\gamma=-4+4\sin^2(\pi/p)$ as being sequences corresponding to discrete groups which are not free - they are in fact triangle groups or their $\IZ_2$-extensions.  This shows sharpness and is   actually part of a more general phenomenon;  the $(p,q,r)$ triangle groups lying on the real part of the axis accumulating on the part freely generated as $r\to\infty$,  the limit being the $(p,q,\infty)$-triangle group,  which is freely generated.

 In terms of the distance between the two axes of the generators we also have the following theorem which is simply an application of our earlier formulas.  Nevertheless it places some constraints on our more general efforts to identify these distances in a moment.

\begin{theorem}
For each $p$ and $q$ not both equal to $2$, we set
\[ b_\infty(p, q) = \arccosh\Big(\frac{\cos(\pi/p) \cos(\pi/q) + 1}{  \sin(\pi/p) \sin(\pi/q)} \Big) . \]
Suppose that $f$ and $g$ are elliptics of order $p$ and $q$ respectively with the distance between their axes bounded below by
\[ \delta(f, g) \geq \delta_\infty (p, q) . \]
Then $\G = \langle f,g\rangle$  is discrete and isomorphic to the free product $\IZ_p*\IZ_q$. The lower bound is sharp in the sense that it is attained in the $(p,q,\infty)$ - triangle group and for every $\epsilon > 0$ there are infinitely many Kleinian groups $\langle f,g\rangle$ generated by elliptics of order $p$ and $q$ with
$\delta_\infty(p,q)-\epsilon ² \delta(f,g) < \delta_\infty(p,q)$ which are not isomorphic to the free product of cyclic groups.
\end{theorem}

We will identify the smallest possible value of the distance between axes for discrete groups in a moment.

In the applications we want we will have to be quite a bit more clever than just  relying on isometric circles.  We ultimately want  to get stronger versions of Theorem \ref{free} and the ideas are basically the same,  but more carefully identifying the fundamental domains and their combinatorial relations yields definitely better results.  There are a couple of remarks worth making.  First, isometric circles are
not conjugacy invariants.  Choosing different conjugacies
of the group $\langle f,g \rangle$ (which amounts to choosing different
matrix normalisations) can significantly simplify the
combinatorial patterns of the intersecting isometric circles.  Even
interchanging $p$ and $q$ (when they are different) sometimes
leads to great simplification. 

\bigskip

When it comes to the parameter spaces for two-generator discrete groups generated by a pair of elliptic elements, Theorem \ref{free} shows that as soon as $|\gamma|$,  ($\gamma=\gamma(f,g)$),  is big enough,  then the group with complex parameters $(\gamma,-4\sin^2(\pi/p),-4\sin^2(\pi/q))$ is discrete, nonelementary and cannot be a lattice. Note here that for all $p$ and $q$,  the trivial bound from Theorem \ref{free} has 
\begin{eqnarray*} |\gamma| &\geq & 4(\cos\pi/p + \cos\pi/q)^2 + 4(\cos\pi/p \cos\pi/q + 1)^2 + 4\sin^2\pi/p \sin^2\pi/q \\
& = &  4(1+\cos(\pi/p)\cos(\pi/q))^2
. \end{eqnarray*}
Thus in our search for lattices we wind up looking for  $\gamma$-values confined to the region bounded by an ellipse in  the complex plane.  We want to get a more thorough description of these spaces and we do this through excluding regions of possible $\gamma$ values by finding inequalities on the parameters for discrete groups - such as the omitted  disks in the Riley slice.

\subsubsection{The space of Kleinian groups.} 

These days the  first important universal constraint one meets when studying the geometry of discrete groups is J\o rgensen's inequality, \cite{Jorgensen}.  This inequality is a generalisation of the earlier Shimitzu-Leutbecher inequality where one assumes that one of the generators in parabolic.  We give the usual proof of this result here,  but another interesting and much more revealing proof of this inequality and other generalisations a bit later after we have discussed polynomial trace identities.  

\begin{theorem}[J\o rgensen's inequality]  Let $\langle f, g\rangle$ be a Kleinian group.  Then
\begin{equation}\label{ji}
|\gamma(f,g)|+|\beta(f)| \geq 1.
\end{equation}
This inequality is sharp
\end{theorem}
\noindent{\bf Proof}.  Suppose that  $r=|\gamma(f,g)|+|\beta(f)| <1$.  Clearly $r>0$ for otherwise $\gamma(f,g)=0$ and  $\langle f,g\rangle$ is elementary.  Put
\[ f \sim \Big(\begin{array}{cc} \lambda & 0 \\ 0 & 1/\lambda \end{array}\Big), \hskip15pt g \sim \Big(\begin{array}{cc} a & b \\c & d \end{array}\Big),  \hskip10pt ad-bc=1 . \]
We recursively define $g_0=g$ and $g_{n+1} = g_{n} \circ f \circ g_{n}^{-1}$.  We calculate that
\[ \beta(f) = (\lambda-1/\lambda)^2, \hskip20pt\gamma(f,g) = - (\lambda-1/\lambda)^2 bc, \]
\[ g_{n+1} = \Big(\begin{array}{cc} a_{n+1} & b_{n+1} \\c_{n+1} & d_{n+1}    \end{array}\Big)  = \Big(\begin{array}{cc} a_nd_n \lambda-b_nc_n/\lambda  & -a_nb_n(\lambda-1/\lambda) \\ c_nd_n(\lambda-1/\lambda) & a_nd_n/ \lambda-b_nc_n\lambda    \end{array}\Big)  . \]
Then 
\[ | b_{n+1}c_{n+1}| = | a_nb_nc_nd_n (\lambda-1/\lambda)^2| = |b_n c_n(1+b_nc_n) (\lambda-1/\lambda)^2|   . \]
Since $r = (1+|bc|) |\lambda-1/\lambda|^2 $ we see by induction that 
\[ | b_{n}c_{n}|   \leq r^n |bc| \to 0, \hskip20pt \mbox{as $n\to\infty$}  . \]
Thus $a_nd_n=1+b_nc_n \to 1$,  $a_n\to \lambda$ and $d_n\to 1/\lambda$.  Next
\[ |b_{n+1}/b_n |= | a_n( \lambda-1/\lambda)| \to |\lambda(\lambda-1/\lambda)| \leq |\lambda| \sqrt{r} \]
which gives
\[ |b_{n+1}/ \lambda^{n+1} |   \leq   |b_n/\lambda^n| (1+\sqrt{r})/2 \]
for all sufficiently large $n$.  Then by induction $b_n/\lambda^n\to 0$ and similarly $c_n\lambda^n \to 0$.  Thus
\[ f^{-n}\circ g_{2n} \circ f^n = \Big(\begin{array}{cc} a_{2n} & b_{2n}/\lambda^{2n} \\c_{2n}\lambda^{2n} & d_{2n} \end{array}\Big) \to \Big(\begin{array}{cc} \lambda & 0 \\ 0 & 1/\lambda \end{array}\Big) = f . \]
This convergence cannot happen in a discrete group,  and so for sufficiently large $n$,  $ f^{-n}\circ g_{2n} \circ f^n = f$ whereupon $g_{2n}= g_{2n-1}\circ f \circ g_{2n-1}^{-1}=f$.  Then $g_{2n-1}$ and $f$ have the same fixed point set,  as $g_n$ is never elliptic of order $2$ as it is a conjugate of $f$ and $|\beta(f)|\leq r <1$.  Inductively we conclude that $f$ and $g$ share fixed points and therefore generate an elementary group contrary to our hypothesis. \hfill $\Box$

\bigskip

We give two important applications of J\o rgensen's inequality.   

\begin{theorem}  A nonelementary subgroup $\G$ of $\mathrm{PSL}(2,\IC)$ is discrete if and only if every two-generator subgroup is discrete.
\end{theorem}
\noindent{\bf Proof.}  The first step in the proof is to show that the supposition that  $\G$ is nonelementary (not virtually abelian) implies  there are loxodromic $g_1$ and $g_2$ in $\G$ with disjoint fixed point sets. We leave this elementary observation to the reader. Next suppose that $\G$ is not discrete.  Then there is a sequence $\{f_j\}_{j=1}^{\infty}$ with $g_j\to \mathrm{identity}$ uniformly in $\oC$.  We may assume that
\[ \{f_j\}_{j=1}^{\infty}\cap \{g_1,g_2\}=\emptyset \]
and further that no $f_j$ has order two.  Now
if every two-generator group is discrete,  for $j$ sufficiently large we have that both groups $\langle f_j,g_1\rangle$ and $\langle f_j,g_2\rangle$ are discrete.  We see that $\beta(f_j)\to 0$ and also,  since $f_j\to id$,  $\gamma(f_j,g_j)\to 0$.  Thus
\[ |\beta(f_j)|+|\gamma(f_j,g_i)|<1,  \hskip20pt i=1,2, \hskip10pt j\geq j_0 \]
and thus $\langle f_{j_0},g_i\rangle$ is elementary.  As $g_i$ is loxodromic either  $f_{j_0}$ fixes or interchanges the fixed points of $g_i$, $i=1,2$.  Since $f_{j_0}$ is not of order two we conclude that  
\[ {\rm fix}(g_1)={\rm fix}(f_{j_0}) = {\rm fix}(g_2)\]
which contradicts our original choice of $g_i$,  $i=1,2$. \hfill $\Box$

\medskip

Actually it is this basic argument that is used quite a lot as we will see. A sequence of $n$-generator Kleinian groups  $\G_i=\langle g_{1}^{i},g_{2}^{i},\ldots,g_{n}^{i}\rangle$ is said to converge algebraically to  $\G=\langle g_{1},g_{2},\ldots,g_{n} \rangle$ if for each $k=1,2,\ldots, n$ we have $g_{k}^{i}\to g_k$.  This convergence will be uniform in the spherical metric of the Riemann sphere $\oC$.

\begin{theorem}\label{closed1} The space of $n$-generator discrete non-elementary groups is closed in the topology of algebraic convergence.
\end{theorem}
\noindent{\bf Proof.}  Let 
\[ \G_i=\langle g_{1}^{i},g_{2}^{i},\ldots,g_{n}^{i}: R_{1}^{i}=R_{2}^{i}= \cdots = R_{m}^{i}=1\rangle \]
be a sequence of discrete non-elementary groups with
\[ g_{j}^{i}\to g_{j}\hskip10pt\mbox{uniformly on $\oC$} . \]
Let $\G=\langle g_{1},g_{2},\ldots,g_{n}\rangle$ and suppose that $\G$ is not discrete.  Then there is a sequence $\{h_k\}_{k=1}^{\infty}\subset \G$ with $h_k\to id$ uniformly on $\oC$ as $k\to\infty$.  We write $h_k$ in terms of the generators $g_{j}$ of $\G$ and then we get approximants $h_{k}^{i}\in \G_i$ by replacing $g_{j}$ with $g_{j}^{i}$.  Then clearly,  as each $h_{k}$ is a finite length word in the generators, 
\[ h_{k}^{i} \to h_k \hskip10pt\mbox{uniformly on $\oC$ as $i\to\infty$} . \]
Upon passing to a subsequence we obtain a sequence $\{f_\ell\}_{\ell=1}^{\infty}$,  
\[ f_{\ell} = h_{k(\ell)}^{i(\ell)} \in \Gamma_{i(\ell)} \]
converging uniformly to the identity. Thus,  for all $\ell$ sufficiently large $|\beta(f_\ell)|<1/4$ and $|\gamma(f_\ell,g_j)|<1/4$ for $j=1,2,\ldots, n$.  Everything here is continuous in the variables and we are in a compact region of $\mathrm{PSL}(2,\IC)$,  so this certainly implies that for all $i$ sufficiently large (and independent of $\ell$), that $|\gamma(f_\ell,g_{j}^{i})|<1/2$ for $j=1,2,\ldots, n$. In particular
\[  |\beta(f_{\ell})|+|\gamma(f_\ell,g_{j}^{i(\ell)})|<3/4 < 1 \]
so J\o rgensen's theorem tells us that 
$\langle g_{k}^{i(\ell)} ,f_{\ell} \rangle\subset \G_{i(\ell)}$ is elementary and discrete.   Since $f_\ell$ is not elliptic of order less than $6$ (as $\beta(f_\ell)\to 0$) it must be that
all the generators of $\G_{i(\ell)}$ fix,  or possibly interchange,  the fixed point set of $f_\ell$ and the group $\G_{i(\ell)}$ cannot be nonelementary. This contradiction shows $\G$ to be discrete. It is just a little more work to show it is nonelementary as well and we leave this as an exercise. \hfill $\Box$

\medskip

More or less the same arguments establish the following theorem.

\begin{theorem}  Suppose that $\G_i=\langle g_{1}^{i},g_{2}^{i},\ldots,g_{n}^{i} \rangle$ converge algebraically to $\G= \langle g_{1},g_{2},\ldots,g_{n}\rangle$ as $i\to\infty$.  
Then the map back is an eventual homomorphism.  That is for all sufficiently large $i$ the map $\G\to\G_i$ given by $g^{j} \mapsto g^{j}_{i}$ extends to a homomorphism of the groups.
\end{theorem}

We leave it to the reader to develop these results into theorems about the space of complex parameters.  For instance

\begin{theorem}\label{closed2}  Let $(\gamma_i,\beta_i,\tilde{\beta}_i)$ be a sequence of parameters for discrete nonelementary (Kleinian) groups.  Suppose that
\[ \gamma_i\to\gamma,  \hskip10pt\beta_i\to\beta, \hskip10pt {\rm and} \;\;\; \tilde{\beta}_i\to\tilde{\beta} . \]
Then $(\gamma,\beta, \tilde{\beta})$ are the parameters for a Kleinian group.
\end{theorem}
With a bit of work one can deduce,  using the eventual homomorphism back,  the further implication that if $\beta=-4\sin^2(m\pi/n)$ and $\tilde{\beta}=-4\sin^2(\tilde{m}\pi/\tilde{n})$,  (generators of finite order),  then for all sufficiently large $i$,  $\beta_i=\beta$ and  $\tilde{\beta}_i=\tilde{\beta}$

 \subsection{The elementary groups} 
 In order to describe the various spaces of two-generator groups with elliptic generators it is necessary to give a complete description of the elementary groups.
 In fact these elementary discrete subgroups of $\mathrm{Isom}^+(\IH^3)$  can be completely classified without too much difficulty.   A. Beardon's excellent book \cite{Beardon} goes through this classification quite carefully and that book is a useful reference to have at hand for many concepts we will meet in this article.  Since this classification will be very important in what follows we record it as follows.
\begin{theorem}[Classification of elementary discrete groups]  Let $G$ be a discrete subgroup of $\mathrm{Isom}^+(\IH^3)$ which is virtually abelian.  Then either
\begin{itemize}
\item $G$ is isomorphic to a cyclic group,  $G = \IZ_p$,  $p=1,2,3,\ldots,\infty$,
\[ \langle a: a^p=1 \rangle; \]  
\item $G$ is isomorphic to a dyhedral group,  $G= D_p \cong \IZ_p \ltimes \IZ_2$,  $p=1,2,3,\ldots,\infty$.
\[ \langle a,b: a^p=b^2 = 1, bab^{-1}=a^{-1} \rangle ;\]  
\item $G$ is isomorphic to the group,  $G= (\IZ_p \times \IZ)\ltimes \IZ_2$ or $G= \IZ_p \times \IZ$.
\[ \langle a,b,c: [a,b]=b^p =c^2= 1, cac^{-1}=a^{-1}, cbc^{-1}=b^{-1} \rangle; \] 
\item $G$ is isomorphic to a finite spherical $(p,q,r)$ triangle group,  one of $A_4=\Delta(2,3,3)$, $S_4=\Delta(2,3,4)$ or $A_5=\Delta(2,3,5)$
\[ \langle a,b : a^p=b^q=(ab)^r= 1 \rangle; \] 
\item $G$ is isomorphic to a Euclidean translation group $\IZ\times \IZ$;
\item $G$ is isomorphic to a Euclidean triangle group,  $\Delta(2,3,6)$, $\Delta(3,3,3)$, $\Delta(2,4,4)$
\[ \langle a,b : a^p=b^q=(ab)^r= 1 \rangle . \] 
\end{itemize}
\end{theorem}
These discrete elementary groups are easy to identify geometrically as their action on $\overline{\IH^3}$ has a point with finite orbit - giving an alternate definition of elementary.

 \subsection{The structure of the spherical triangle groups}
 
 It is going to be very important to us to have a thorough understanding of the finite spherical triangle subgroups as well as the Euclidean triangle subgroups of  a Kleinian group.
Therefore we next recall a few facts about them that we will need. We will simply
tabulate results.  The verification of these tables is a sometimes lengthy calculation in spherical trigonometry.

\begin{lemma}
Let $\theta_{2,3}$  and $\theta_{3,3}$ denote,  respectively,  the angle subtended at the origin between the axes of
order 2 and 3 and the axes of order 3 of a spherical (2,3,3)--triangle group (fixing the origin).  Then
\begin{equation}
\cos(\theta_{2,3}) = 1/\sqrt{3},  \hskip40pt
\cos(\theta_{3,3}) = 1/3.
\end{equation}
\end{lemma}

\begin{lemma}
Let $\phi_{2,3}$, $\phi_{2,4}$ and $\phi_{3,4}$ denote, respectively, the angle subtended at
the origin between the axes of order 2 and 3, the axes of order 2 and 4  and the axes
of order 3 and 4 of a spherical (2,3,4)--triangle (fixing the origin).  Then
\begin{equation}
\cos(\phi_{2,3}) = \sqrt{2/3}, \hskip10pt
\cos(\phi_{2,4}) = 1/\sqrt{2}, \hskip10pt
\cos(\phi_{3,4}) = 1/\sqrt{3}
. \end{equation}
\end{lemma}

\begin{lemma}
Let $\psi_{2,3}$, $\psi_{2,5}$ and $\psi_{3,5}$ denote,respectively, the angle subtended at
the origin between the axes of order 2 and 3, the axes of order 2 and 5 and the axes
of order 3 and 5 of a spherical (2,3,5)--triangle (fixing the origin).  Then
\begin{equation} \left.
\begin{array}{ccc}
\cos(\psi(2,3)) & = & 2\cos(\pi/5)/\sqrt{3}, \\
\cos(\psi(2,5)) & = & 1/2\sin(\pi/5), \\
\cos(\psi(3,5)) & = & \cot(\pi/5)/\sqrt{3}.
\end{array} \hskip30pt \right\}
. \end{equation}
\end{lemma}

From these elementary observations and some spherical trigonometry we obtain the following
tables which list the possible angles, up to four decimal places, at which elliptic axes in a Kleinian group can
meet.  We also include the dihedral angle $\psi$ opposite this angle in the spherical triangle formed when the adjacent
vertices are $\pi/p$ and $\pi/q$.

\bigskip
{\scriptsize
\begin{center}
\begin{tabular}{|c|c|c|c|c|}
\multicolumn{5}{c}{\bf  Angles: $2,\;m$ axes}\\
\hline
$m$&$ \sin(\theta)$ & $\theta $ & $\psi$ & Group \\
\hline
$ 3 $&$ \sqrt{\frac{2}{3}} $ & $.9553 $ & $ \pi/3$  &$ A_4 $\\
\hline
$ 3 $&$ -\sqrt{\frac{2}{3}} $ & $2.1462 $ & $ 2\pi/3$  &$ A_4 $\\
\hline
$ 3 $&$ \sqrt{\frac{1}{3}} $ & $.6154$ & $\pi/4 $ &$ S_4 $\\
\hline
$ 3 $&$ -\sqrt{\frac{1}{3}} $ & $2.5261$ & $3\pi/4 $ &$ S_4 $\\
\hline
$ 3 $&$ \sqrt{\frac{3-\sqrt{5}}{6}} $& $.3648 $ & $\pi/5$ &$A_5 $\\
\hline
$ 3 $&$ -\sqrt{\frac{3-\sqrt{5}}{6}} $& $2.7767 $ & $4\pi/5$ &$A_5 $\\
\hline
$ 3 $&$ \sqrt{\frac{3+\sqrt{5}}{6}} $& $ 1.2059$ & $2\pi/5 $ &$ A_5 $\\
\hline
$ 3 $&$ -\sqrt{\frac{3+\sqrt{5}}{6}} $& $1.9356$ & $3\pi/5 $ &$ A_5 $\\
\hline
$ 3 $&$ 1 $& $ \pi/2 $ & $ \pi/2$ &$ D_3 $\\
\hline
$ 4 $&$ \frac{1}{\sqrt{2}} $& $\pi/4 $ & $ \pi/3$ & $ S_4 $\\
\hline
$ 4 $&$ \frac{1}{\sqrt{2}} $& $3\pi/4 $ & $ 2\pi/3$ & $ S_4 $\\
\hline
$ 4 $&$ 1 $& $\pi/2 $ & $\pi/2 $ & $ D_4 $\\
\hline
$ 5 $&$ \sqrt{\frac{5-\sqrt{5}}{10}} $& $ .5535 $ & $\pi/3 $ &$ A_5 $\\
\hline
$ 5 $&$ -\sqrt{\frac{5-\sqrt{5}}{10}} $& $ 2.5880 $ & $2\pi/3 $ &$ A_5 $\\
\hline
$ 5 $&$ \sqrt{\frac{5+\sqrt{5}}{10}} $& $ 1.0172 $ & $ 2\pi/5 $ &$ A_5 $\\
\hline
$ 5 $&$ -\sqrt{\frac{5+\sqrt{5}}{10}} $& $ 2.1243 $ & $ 3 \pi/5 $ &$ A_5 $\\
\hline
$ 5 $&$ 1 $& $\pi/2 $ & $\pi/2 $ &$ D_5 $\\
\hline
\end{tabular}
\begin{tabular}{|c|c|c|c|c|}
\multicolumn{5}{c}{\bf  Angles: $3,\; m$ axes}\\
\hline
$m$&$ \sin(\theta) $ & $\theta $& $\psi $ & Group \\
\hline
$ 3 $&$ \frac{2}{3} $& $ .7297 $ & $2\pi/5$ &$ A_5 $\\
\hline
$ 3 $&$ -\frac{2}{3} $& $ 2.4118 $ & $3\pi/5$ &$ A_5 $\\
\hline
$ 3 $&$ \frac{2\sqrt{2}}{3} $& $1.2309 $ & $\pi/2 $ &$A_4 $\\
\hline
$ 3 $&$ -\frac{2\sqrt{2}}{3} $& $1.9106 $ & $2\pi/3 $ &$A_4 $\\
\hline 
$ 4 $&$ \sqrt{\frac{2}{3}} $& $.9553$ & $\pi/2 $ &$S_4 $\\
\hline
$ 4 $&$ -\sqrt{\frac{2}{3}} $& $2.1462$ & $3\pi/4 $ &$S_4 $\\
\hline
$ 5 $&$ \sqrt{\frac{10-2\sqrt{5}}{15}} $& $.6523$ & $\pi/2 $ &$ A_5 $\\
\hline
$ 5 $&$ -\sqrt{\frac{10-2\sqrt{5}}{15}} $& $2.4892$ & $4\pi/5 $ &$ A_5 $\\
\hline
$ 5 $&$ \sqrt{\frac{10+2\sqrt{5}}{15}} $& $1.3820$ & $3\pi/5 $ &$ A_5 $\\
\hline
$ 5 $&$ -\sqrt{\frac{10+2\sqrt{5}}{15}} $& $1.7595$ & $2\pi/3 $ &$ A_5 $\\
\hline
\end{tabular}
\end{center}
}

We include the following additional remarks.

\begin{enumerate}
\item The angle between intersecting axes of elliptics of order $4$ in a discrete group
is always either $0$ when they meet on the Riemann sphere or $\pi/2$.  In the first case the dihedral angle
$\psi$ as above is $\pi/2$ while in the second case it is $2\pi/3$.
\item The angle between intersecting axes of  elliptics of order $5$ in a discrete
group is either $\arcsin(2/\sqrt{5})=1.107149\ldots$ or its complement $\arcsin(-2/\sqrt{5})=2.03444\ldots$.  In the former
case the dihedral angle $\psi$ is $2\pi/3$. In the latter it is $4\pi/5$.
\item The axes of elliptics both of order $2$ can intersect at an angle
$k\pi/m$ for any $k$ and $n\geq 2$.  In this case the dihedral angle $\psi$ is $k\pi/m$.
\item The axes of elliptics of order $p$ and $q$,  $p\leq q$, in a discrete group meet on the
sphere at infinity, i.e., meeting with angle $0$, if and only if  
\[ (p,q)\;\; \in \{(2,2),(2,3), (2,4), (2,6), (3,3), (3,6), (4,4), (6,6)\}. \]
In each case the dihedral angle $\psi$ is $\pi-\pi/p -\pi/q$.
\end{enumerate}  
 
\section{Polynomial trace identities and inequalities}

We shall now discuss a very interesting family of polynomial
trace identities which will be used to obtain geometric information about
Kleinian groups, \cite{GM3}.

\medskip

Let $\langle a,b \rangle$ be the free group on the two letters $a$
and $b$.  We say that a
word $w\in \langle a, b \rangle$ is a {\em good word} if $w$ can be written as
\begin{equation}
\label{goodword}
w = a^{s_1} b^{r_1} a^{s_2} b^{r_2} \ldots a^{s_{m-1}} b^{r_{m-1}}a^{s_{m}}
\end{equation}
where $s_1\in \{\pm 1\}$, $s_j=(-1)^{j+1} s_1$ and $r_j\neq 0$ but
are otherwise
unconstrained.  The good words start and end in $a$ and the exponents
of $a$ alternate in
sign.  The following theorem is a key tool used in the study of the
parameter spaces of
discrete groups.

\begin{theorem} \label{words} Let $w=w(a,b)\in \langle a,b \rangle$
be a good word.   Then
there is a
monic polynomial
$p_w$  of two complex variables having integer coefficients with the following property.

Suppose that $f,g\in \mathrm{PSL}(2,\IC)$ with  $\beta=\beta(f)$ and $\gamma=\gamma(f,g)$.  Then
\begin{equation}
\gamma(f,w(g,f)) = p_w(\gamma,\beta).
\end{equation}
where $w(g,f)$ is the word in $\langle f,g\rangle$ found from the assignment $f\leftrightarrow b$, $g\leftrightarrow a$
\end{theorem}

There are three things to note. The first is that if we assume that
$a^2=1$,  then the alternating 
sign condition is redundant since $a=a^{-1}$ and every word $w(a,b)$ is good.
   The second thing is that there is a
natural semigroup operation on the good words. If $w_1=w_1(a,b)$ and
$w_2=w_2(a,b)$ are good words,  then so is
\begin{equation} w_1*w_2 = w_1(a,w_2(a,b)).
\end{equation}
That is,  we replace every instance of $b$ in $w_1$ with $w_2(a,b)$.  So for example
\[  (bab^{-1}ab)*(bab^{-1})= bab^{-1}a(bab^{-1})^{-1}abab^{-1} = bab^{-1}aba^{-1} b^{-1}abab^{-1}. \]
It is not too difficult to see the remarkable fact that
\begin{equation} p_{w_1*w_2}(\gamma,\beta) =
p_{w_1}(p_{w_2}(\gamma,\beta),\beta) \end{equation}
which corresponds to polynomial composition in the first slot.   

Notice the obvious fact that $\langle f,g\rangle$ Kleinian,  implies $\langle f,w(g,f)\rangle$ discrete.  Finally note that for any word $w=w(g,f)$ and $m,n\in \IZ$,  $\gamma(f,f^mwf^n)=\gamma(f,w)$ so that the requirement that the word start and end in a nontrivial power of $b$ is simply to avoid some obvious redundancy.

\bigskip

Let us give two simple examples of word polynomials and how they
generate inequalities.  This might seem an aside to our task of
studying small covolume
lattices, but actually our  searches for the $\gamma$ values of discrete groups with given generators amount in large part to
systematising the following arguments.  The direct calculation of the polynomial $p_w$ from $w$  by hand can be a little tricky,  but can be done for some short words $w$,  see  \cite{GM3}.

 \medskip

First
is the classical example:
\begin{equation}
\label{aba}
   w = aba^{-1}, \hskip30pt  p_w(z) = z(z-\beta). \end{equation}
Here $z=\gamma(f,g)$  and $\beta=\beta(f)$.  In our applications we will suppress $\beta$ (that is treat it as a coefficient) and
treat $z$ as the variable.  This is typically because we will be restricting to slices in the space of discrete groups where $\beta$ is constant.
The next example comes from the good word
\begin{equation} w = aba^{-1}b^{-1}a, \hskip30pt  p_w(z) = z(1+\beta-z)^2.
\end{equation}
Another is
\[
w  = bab^{-1}a^{-1}b, \hskip10pt
p_w(z) = z(1-2\beta+2z-\beta z+z^2) . \]

Let us indicate how these words are used to describe parts of the parameter
space for  two-generator Kleinian groups.  We take for granted the
well known fact that the
space of finitely generated discrete non-elementary groups is closed.  We have used J\o rgensen's inequality earlier to prove this,  but it is a very
general fact concerning groups of isometries of negative curvature, see for instance \cite{Martinnegcurv}.

\subsection{J\o rgensen's inequality revisited}
Let
us recover J\o rgensen's inequality from the first trace polynomial.  Consider
\[ \min\{|\gamma|+|\beta|: (\gamma,\beta,\tilde{\beta}) \; \mbox{ are the
parameters of a Kleinian group }\} . \]
This minimum is attained by some Kleinian group $\G=\langle f,g
\rangle$.  If $\langle gfg^{-1} , f \rangle$ is Kleinian (it is
certainly discrete),  then,  by minimality,
\begin{eqnarray} |\gamma|+|\beta| & \leq & |\gamma(\gamma-\beta)|+|\beta| \nonumber \\
1 & \leq & |\gamma-\beta| , \label{jorgsharp}
\end{eqnarray}
as $\gamma\neq 0$ (we reiterate $\gamma=0$ implies that $f$ and $g$ share a fixed
point on $\oC$  so $\G$ would not be Kleinian).
If $f$ has order $2,3,4$ or $6$,
$|\beta|\geq 1$.  In all other cases $\gamma\neq \beta$ as $\G$ is
non-elementary.  Thus at the mimimum we have $|\gamma-\beta|\geq 1$
and so
$|\gamma|+|\beta|\geq 1$ in general.  This is J\o rgensen's
inequality.  This inequality is attained with equality (indeed (\ref{jorgsharp}) holds with equality) for representations of the (2,3,p)-triangle groups.  A point to observe here is that  we must examine and eliminate, for some geometric reason, the zero
locus of $p_w$ (here the variety
$\{\gamma=\beta\}$) to remove the possibility that $\langle
f,w(f,g) \rangle $ is a discrete elementary group.  Further,  it seems we are tacitly using Lemma \ref{beta=-4} in these
arguments to remove the relevance of the parameter $\tilde{\beta}$.

In a similar fashion,  if we minimize $|\gamma|+|1+\beta|$ and use
the second polynomial $z(1+\beta-z)^2$ we see that at the minimum
\begin{eqnarray*} |\gamma|+|1+\beta| & \leq &
|\gamma(1+\beta-\gamma)^2|+|1+\beta| \\
1 & \leq & |1+\beta-\gamma|,
\end{eqnarray*}
and so $|\gamma|+|1+\beta|\geq 1$ at the minimum.  The zero locus
this time is the set $\{\gamma=1+\beta\}$,  these groups are Nielsen
equivalent to groups
generated by elliptics of order $2$ and $3$ \cite{GM3}.  In
particular this gives us the following inequality which we will use
later.

\begin{lemma} \label{modjorg} Let $\langle f,g\rangle$ be a Kleinian
group.  Then
\[ |\gamma|+|1+\beta|\geq 1 \]
unless $\gamma=1+\beta$ and $fg$ or $fg^{-1}$ is elliptic of order $3$.
\end{lemma}
As a consequence,  if $f$ has order $6$,  $\beta=-1$ and we have
   \begin{corollary}  If $\langle f,g \rangle$ is a Kleinian group and
$f$ is elliptic
of order $6$,  then
\begin{equation}\label{SL6} |\gamma(f,g)| \geq 1.\end{equation}
This result is sharp.
\end{corollary}
This result is entirely analogous to the Shimitzu--Leutbecher inequality above.
\medskip

There are a couple of further points we wish to make here.  Suppose
we have eliminated a certain region from the possible values for
$(\gamma,\beta,-4)$ among
Kleinian groups.  Then $(p_w(\gamma,\beta),\beta,-4)$ also cannot lie in
this region.  Thus for instance,  from J\o rgensen's inequality we
have
\begin{lemma} \label{geqlem} Let $\langle f,g\rangle$ be a Kleinian group.  Then
\begin{equation}\label{genineq}
|p_w(\gamma,\beta)|+|\beta|\geq 1 
\end{equation}
unless $p_w(\gamma,\beta)=0$.
\end{lemma}
We have already seen this idea in an example,  discussed earlier when we were looking at the Riley slice.  For a parabolic generator we have $\beta=0$ and the Shimitzu-Leutbecher inequality $|\gamma(f,g)|\geq 1$.  If $p_w(\gamma)=\gamma(1+\beta-\gamma)^2 = \gamma(1-\gamma)^2$, for $\beta=0$,  then $|p_w(\gamma)|\geq 1$ unless $p_w(\gamma)=0$ and the group $\langle f,w\rangle$ is elementary.  This gives the following lemma.

\begin{lemma} Let $\langle f,g\rangle$ be a Kleinian group with $f$ parabolic and $\gamma=\gamma(f,g)$.  Then either $\gamma=1$ or
\[ |\gamma(1-\gamma)^2|\geq 1 \]
and so in particular $|1-\gamma|\geq r_0$,  $(1+r_0)r_{0}^{2}=1$.
\end{lemma}

A useful special case (using (\ref{aba})) is

\begin{lemma}
\label{specialcase}
Let $\langle f,g\rangle$ be a Kleinian group.  Then
\begin{equation}
|\gamma(\gamma-\beta)|+|\beta|\geq 1
\end{equation}
unless $\gamma= \beta$.
\end{lemma}

We remark that the condition $p_w(\gamma,\beta)=0$ implies that  $f$  and $w=w(g,f)$
share a fixed point on $\oC$.  We may use discreteness and geometry
to find implications of this fact which can be used to eliminate it.  For instance if
$f$  is loxodromic, this implies that
$w$ shares both fixed points with $f$  since we know that $\langle
f, w\rangle$ is discrete by hypothesis.  Thus $[f,w]=id$.  If
$f$ has some other special
property such as being primitive,  then we either have $w$ elliptic
or a power of $f$,  further relations.  These additional properties
of $w$ place greater
constraints on $\gamma$ and $\beta$.

More information is garnered in special cases.  Here are a number of identities that we use in studying the distance between octahedral vertices in a Kleinian group.  It is fairly typical of the sorts of polynomial trace identities that we need in other cases too.
\begin{lemma} Let  $f$ and $g$ be elliptic M\"obius transformations of order $4$ and $2$ respectively.  Set $\gamma=\gamma(f,g)$.  Then
\begin{itemize}
\item $\gamma((gf)^4g,f)= \gamma(-1+\gamma+\gamma^2)^2$
\item $\gamma((gf)^3g,f)= \gamma^3(2+\gamma)$
\item $\gamma((gf)^3(gf^{-1})^3g,f)=  -2+(2+\gamma)(1+\gamma^2+\gamma^3)^2$
\item $\gamma((gf)^3(gf^{-1})^3 (gf)^3g,f)= \gamma(2+\gamma)(1+2\gamma+\gamma^2+2\gamma^3+\gamma^4)^2$
\end{itemize}
\end{lemma}
 Of course,  the composition of these polynomials gives further polynomials corresponding to trace identities.   Further, at the roots of these polynomials we can identify some potential $\gamma$ values which {\em may} correspond to discrete groups - as we will see when we consider arithmeticity.  As examples here,  the roots of the first polynomial, that is of $ -1+\gamma+\gamma^2 $ are $(-1\pm\sqrt{5})/2=4\cos^2(\pi/5)-2$ corresponding to the $(2,4,5)$-triangle group.  A root of the last polynomial $1+2\gamma+\gamma^2+2\gamma^3+\gamma^4=0$ is
 \[-\frac{1}{2}+\frac{1}{\sqrt{2}}+\frac{1}{2} i \sqrt{1+2 \sqrt{2}}\approx 0.207107 + 0.978318 i \] 
corresponding to a two-generator arithmetic lattice,  generated by elements of order $2$ and $4$ and of covolume approximately $1.0287608$.

\subsection{Simple axes}

An elliptic or loxodromic element $f$ in a Kleinian group $\G$ is called {\em simple} if for all $h\in \Gamma$
\begin{equation}\label{simple}
h(\ax(f)) \cap \ax(f) = \emptyset \hskip15pt {\rm or} \hskip15pt h(\ax(f))
= \ax(f).
\end{equation}
Thus,  for simple $f$,  the translates of the axis of $f$ will form a disjoint collection of
hyperbolic lines.  It is a well know fact that the shortest geodesic (should it exist) in a hyperbolic orbifold is always simple. To see this we simply cut at a crossing point and smooth the angles producing two shorter curves,  one of which must be homotopically nontrivial and therefore represented by a shorter geodesic. Finite volume orbifolds always have a shortest geodesic,  as do all geometrically finite Kleinian groups.  In fact in these cases the spectrum of traces $\{\beta(g):g\in\G\}\subset\IC$ is discrete,  and so there is a shortest translation length.

More generally,  given a set $X\subset \IH^3$ we will say that $X$ is   {\em precisely invariant} if for
all $h\in \Gamma$
\begin{equation}\label{simple2}
h(X) \cap X = \emptyset \hskip15pt {\rm or} \hskip15pt h(X)
= X,
\end{equation}
This term was introduced by  Maskit.  A {\em collar} of
radius $r$ about a subset $E\subset \IH^3$ is
\[ {\cal C}(E,r)=\{x\in \IH^3:\rho(x,E)<  r \}, \]
where $\rho$
denotes hyperbolic distance.  A collar about a simple axis is a solid hyperbolic cylinder.

The {\em collaring radius} radius of a non-parabolic $f \in \Gamma$
is the largest $r$ for which $\ax(f)$ has a
precisely invariant collar of radius $r$.

If  $g$ is an elliptic element which is not simple,  then $g$ lies in a
triangle subgroup which is either spherical (finite) or euclidean, 
should $\IH^3/\Gamma$ not be compact.  Accordingly 

\begin{theorem}\label{psimple}
An elliptic element of order $p\geq 7$ in a Kleinian group $\G$ is simple.
\end{theorem}
In a moment we will be able to quantify this result.

The following theorem is obvious,  but it underpins our volume estimates.

\begin{theorem}  Let $X$ be a precisely invariant set for $\G$ acting on $\IH^3$ and let 
\[ \G_X = \{g\in \G:g(X)=X\}\]
denote the stabiliser of $X$.  Then
\begin{equation}\label{simpleineq}
{\rm vol}_\IH(\IH^3/\G) \geq {\rm vol}_\IH(X/\G_X)
. \end{equation}
\end{theorem}
The point here is that for collars about axes, or $X$ a ball in hyperbolic space, the group $\G_X$ is elementary and is completely described.  Thus we can compute ${\rm vol}_\IH(X/\G_X)$.  Of course these arguments are never going to be sharp since the inequality at (\ref{simpleineq}) will be strict unless $X$ is a fundamental domain.

In \cite{GM2} we show that given elliptic elements $f$  and  $g$  of
order $p$ and $q$ generating
a Kleinian group,  the allowable (hyperbolic) distances $\delta_i(p,q)$
between their axes has an initially discrete spectrum and, crucially,
the first several initial values of the spectrum (at least for $p,q \leq 6$)
are uniquely attained for arithmetic lattices \cite{GMMR}.  These estimates came from a description of the space of two-generator discrete groups,  with one generator of order $2$ and the other $p\geq 3$.  We'll discuss this in much more detail for $p$ small,  but for $p\geq 7$ we have the following refinement of Theorem \ref{psimple}.  For $p,q\geq 2$ define
\[ \delta_0(p,q) = \inf_\G \delta(f,g) \]
where the infimum is taken over all Kleinian groups $\G=\langle f,g\rangle$ generated by elliptics of order $p$ and $q$ respectively.  Notice that as soon as $p$ or $q$ is at least $7$ and neither $p$ nor $q$ is $2$,  then discreteness alone will imply the group is Kleinian. Further, if one of $p$ or $q$ is $2$,  then the group is Kleinian unless it is dihedral.

\begin{theorem}[Collaring radii for elliptics of order $p\geq 7$]  For $p,q\geq 7$,
\begin{equation}
\cosh(\delta_0(p,q)) = \frac{1}{2\sin(\pi/p)\sin(\pi/q)}
. \end{equation}
Each value is achieved.
\end{theorem}
As we will see in a moment,  this theorem can be used to provide good volume bounds when there is large order torsion ($p\geq 7$) in a group. A large part of the work in solving the general problem is establishing a similar (but much stronger) result for low-order torsion.

\medskip

Given a simple axis $\alpha$,  its stabiliser is an elementary group.  This group consists of possibly a rotation about $\alpha$,  possibly a loxodromic which shares its axis with $\alpha$ and has translation length $\tau$  and possibly an involution interchanging the endpoints of $\alpha$.  Therefore given  a Kleinian group $\G$ which has  a simple axis $\alpha$ with collaring radius $r$,  stabilised by an elliptic of order $p$ we have,  with $X={\cal C}(\alpha,r)$, 
\begin{equation} 
{\rm vol}_\IH(\IH^3/\G) \geq {\rm vol}_\IH(X/\G_X)  = \frac{\pi \tau \sinh^2(r)}{2p}
\end{equation}
and with $p\geq 7$ we may put in the value $r = \delta_0(p,p)/2$ (note that the collar radius is half the distance between the axis and its nearest translate) to get
\begin{equation} \label{vest1}
{\rm vol}_\IH(\IH^3/\G) \geq   \frac{\pi \tau (\cosh(\delta_0(p,p))-1)}{4p} =   \frac{\pi \tau \cos(2\pi/p)}{8p \sin^2(\pi/p) }  
. \end{equation}
So now we need to get a handle on $\tau$.

 In
\cite{GM8} a
``collar-volume'' formula is established which gives covolume
estimates simply in terms of the collaring radius of a simple axis and is found by considering the ways maximal collars pack.  But here we will just proceed with our current argument. 

\subsection{Killing holonomy}

We should expect that if a loxodromic $f$ has short translation length $\tau$,  then something like J\o rgensen's inequality will give us the control we need.  The problem is that $\tau$ small does not imply that $\beta(f)$ is small (and hence $|\gamma(f,hfh^{-1})|$ is correspondingly large,  where $hfh^{-1}$ is the closest translate of $f$ to itself). The idea is to use powers of $f$ to kill the holonomy since that is reduced mod $2\pi$.  But this grows the translation length,  so a balancing argument is needed.  When Meyerhoff first addressed this problem in \cite{Meyer1} he used a lemma of Zagier.  In \cite{CaoGM} we considered a different number theoretic approach using old ideas about lattice constants in $\IC$. It is rather more complicated,  but gives cleaner and sharper formulas. Set
\[c_1 =2.97, c_2 =1.91, \hskip5pt{\rm and} \hskip10pt c_p = \frac{\sqrt{3}\pi}{p},  \hskip10pt p\geq 3 . \]
\begin{theorem}\label{betachoice}
Let $p\geq 1$ be an integer and $0<\tau \leq c_p$.  Then there are integers $m\geq 1$ and $n$ so that
\[ |4\sinh^2(m(\tau+i\theta)/2+in\pi/p)|\leq \frac{4\pi}{\sqrt{3}p} \tau . \]
\end{theorem}
\noindent This theorem is sharp.  What we now want is the following corollary. It is simply obtained by choosing $h=f^mg^n$ where $g$ is primitive elliptic of order $p$ sharing its axis with $f$ and $m,n$ realise the minimum on the left-hand side of (\ref{beta choice}).
\begin{corollary} Let $\alpha$ be a simple axis in a Kleinian group stabilised by a loxodromic with translation length $\tau$ and an elliptic of order $p\geq 1$.  If $\tau\leq c_p$, then there is $h\in \G_\alpha$ with
\begin{equation}\label{beta choice} |\beta(h)| \leq \frac{4\pi}{\sqrt{3}p} \tau. \end{equation}
\end{corollary}
\subsection{Lattices with high-order torsion, $p\geq 6$}

There are now two cases to consider.

\medskip

\noindent {\bf (i)}. 
If $\tau\geq c_p$ we have the volume estimate from (\ref{vest1})
\begin{equation} 
{\rm vol}_\IH(\IH^3/\G) \geq    \frac{\pi^2\sqrt{3}\cos(2\pi/p)}{8p^2 \sin^2(\pi/p) }  \geq 0.1444\ldots
. \end{equation}
since the minimum occurs when $p=7$ (actually this gives good estimates down to $p=4$,  but the formula for $\delta_0(p,p)$,  $p=4,5,6$ is quite different.)

\medskip

\noindent {\bf (ii)}.  If $\tau\leq c_p$,  then we set $g$ as the nearest translate of $f$ and choose $h\in \G_\alpha$ with $|\beta(h)|\leq \frac{4\pi}{\sqrt{3}p} \tau$ using (\ref{beta choice}).  Then from  (\ref{gammageom}) we set $\delta+i\theta$ as the complex distance between $h$ and its nearest translate (the same distance as $f$ and its nearest translate of course) to see by J\o rgensen's inequality
\begin{eqnarray*}
|\sinh(\delta+i\theta) |^2& = & \Big|\frac{4\gamma(h,h')}{\beta(h)^2} \Big| \geq \frac{4(1-|\beta(h)|)}{|\beta(h)|^2} \geq \frac{ 3p^2-4\pi \tau p\sqrt{3} }{4\pi^2 \tau^2}. 
\end{eqnarray*}
Then
\begin{eqnarray} 
\cosh(\delta) & \geq &\sqrt{ \frac{ 3p^2-4\pi \tau p\sqrt{3} }{4\pi^2 \tau^2} }\nonumber  \\  
  \frac{\pi \tau \sinh^2(\delta/2) }{2p} & \geq &\frac{\sqrt{3}}{8 }\sqrt{ 1 -\frac{4\pi \tau}{p \sqrt{3}}  } -\frac{\pi\tau}{4p} 
. \end{eqnarray}
This estimate is actually decreasing in $\tau$,  so we put $\tau=c_p$ to get 
\begin{equation} \label{vest2}
{\rm vol}_\IH(\IH^3/\G)   \geq  \frac{\sqrt{3}}{8 }\sqrt{ 1 -\frac{4\pi^2 }{p^2 }  } -\frac{\sqrt{3}\pi^2 }{4p^2} 
. \end{equation}
For $p\geq 9$ this estimate gives volume at least $0.1$.  For $p=7$ and $8$ the best way forward is to balance the two estimates at (\ref{vest1}) and (\ref{vest2}).  In effect,  this amounts to choosing a smaller $\tau$ for the first bound. 

\scalebox{0.6}{\includegraphics[viewport= 30 500 280 800]{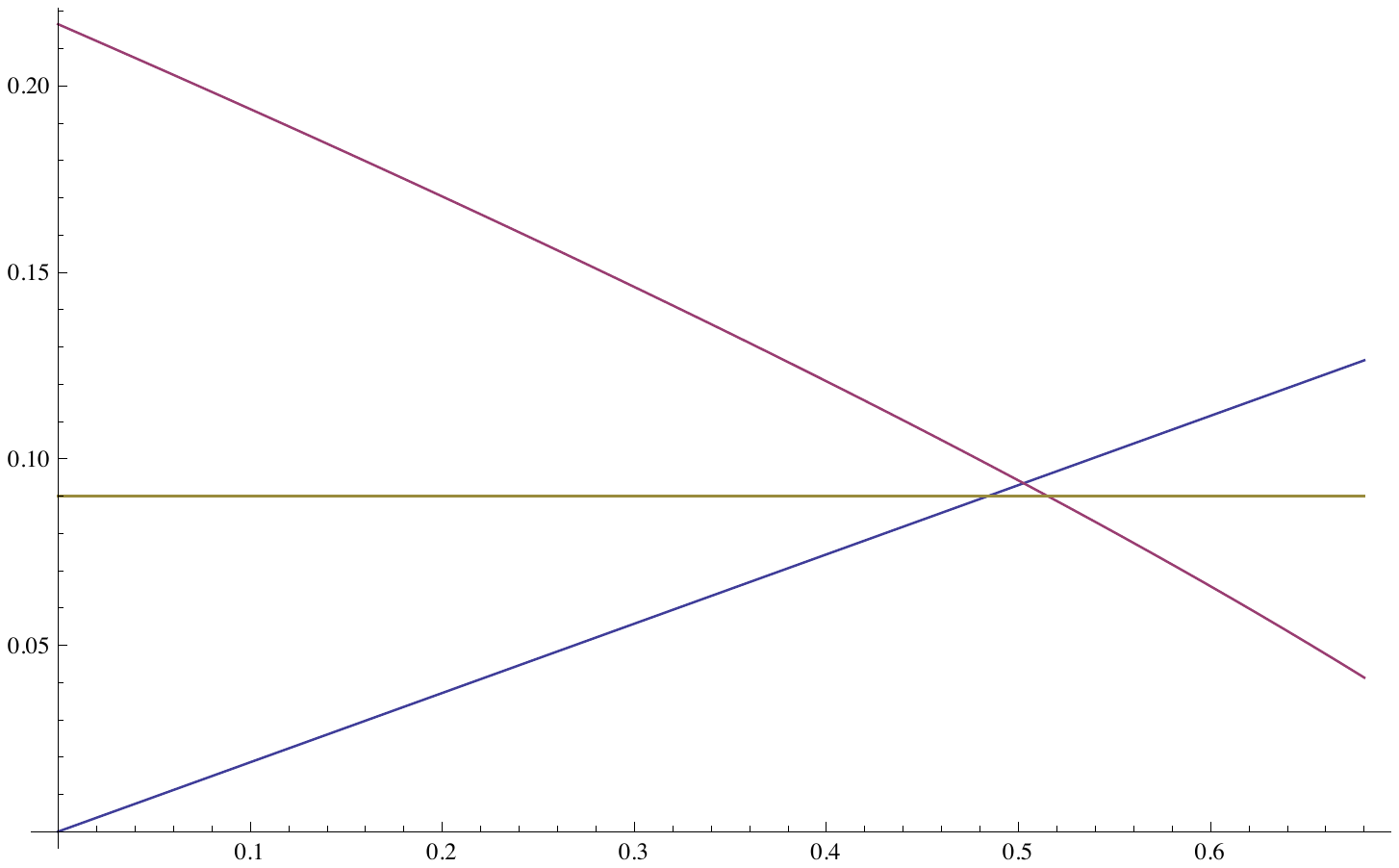}}

The graphs for $p=7$ from the estimates at (\ref{vest1}) and (\ref{vest2}) and the level $v=0.09$  showing that $0.09$ is a lower bound for the volume of a Kleinian group with an elliptic of order $7$.

\medskip

Finally we note that when $p=6$ we are able to use the   Shimitzu-Leutbecher type inequality at (\ref{SL6}) rather than the less strong but more general J\o rgehsen inequality.  This will give substantially better estimates, the only exceptional case to analyse is  when $f$ and its nearest translate share a single  fixed point.  However,  this will imply the existence of a parabolic element and in turn identify  the $(2,3,6)$-Euclidean triangle group as a subgroup.  This discussion has therefore sketched a proof of   the main result of this section.

\begin{theorem} \label{ST1}Let $\Gamma$ be a Kleinian group containing an elliptic of order $p\geq 6$.  Then
\begin{equation}
{\rm vol}_\IH(\IH^3/\G) \geq {\rm vol}_\IH(\IH^3/ PGL(2,{\cal O}_{\sqrt{-3}})) = 0.0846\ldots
. \end{equation}
This result is sharp.
\end{theorem}

The sorts of arguments described above can be used to get volume estimates for simple torsion of other orders ($p=3,4$ and $5$) as well.  Here one only tries to achieve a bound greater than $0.04$ since the minimal covolume lattice has volume smaller than this (and does not have a simple elliptic),  but it is necessary to know much more about the spectrum of possible axial distances in these cases.  This means identifying all the possible $\gamma$ values. This requires much more information and we develop this over the next two sections.  But ultimately the ideas are the same once one has knowledge of the collaring radii.  We record what comes out of that analysis.

\begin{theorem}\label{ST2} Let $\Gamma$ be a Kleinian group containing a simple elliptic of order $p\in\{3,4,5\}$.  Then
\begin{equation}
{\rm vol}_\IH(\IH^3/\G) \geq 0.041
. \end{equation} 
\end{theorem}

This result is rather far from being sharp and we do not have a good idea of what the extremal might be.  A candidate would be $\G_{3,6}$ of covolume $0.0785\ldots$ described in \S 9.  That the elliptic of order three is simple is shown in \cite{GMMR}.

\medskip

In order to advance our attempts to find the $\gamma$ values for discrete groups which are not free on their generators we need to develop some arithmetic machinery.   First note that inequalities such as that at (\ref{genineq})  in Lemma \ref{geqlem} have the restriction that $p_w(\gamma)\neq 0$.  It does not assert that if  $p_w(\gamma)= 0$,  then we have a discrete group. In fact this is most often  to the case and some other arguments are need to eliminate these points,  though we will have enough information to explicitly construct the group to test for discreteness in other ways.  It would be nice to give a condition that tell us when the equation $p_w(\gamma)= 0$ implies we do have a discrete group,  at that $\gamma$ value  it would be also be good to have at hand a  description of the precise group in question.  This problem seems set up for arithmetic considerations which we will now turn to.

 \section{Arithmetic hyperbolic geometry}
 
The main aim of this section is to give an arithmetic criterion which is sufficient to imply the discreteness
of various two-generator subgroups of $\mathrm{PSL}(2,\IC)$.  We also want to establish the finiteness of the family of all two-generator arithmetic groups generated by elements of finite order (or parabolic) and then say as much as we can about this finite family of groups and their associated orbifolds.  

\medskip

Basically an {\em arithmetic Kleinian group} $\Gamma$ is one algebraically isomorphic to something commensurable with $SL(m,\IZ)$ for some $m$.  So there's a representation $\rho:\Gamma\to GL(m,\IC)$ so that 
\[ |\rho(\Gamma)\cap SL(m,\IZ): SL(m,\IZ)|  +  |\rho(\Gamma)\cap SL(m,\IZ): \rho(\Gamma)|  <\infty . \]
Borel gave a nice description of these groups in three dimensions, \cite{Borel}, and  a thorough account is given by Maclachlan and Reid in their book \cite{MR}.   We point out that 
in higher rank,  all lattices in semi-simple Lie groups are arithmetic.  This is Margulis super-rigidity \cite{Margulis1,Margulis2} and even in rank one only hyperbolic and complex hyperbolic lattices need not be arithmetic.

\medskip

For completeness we recall the precise definition of an arithmetic group obtained from the  work of Borel and Vinberg \cite{Vin} when restricted to the three dimensional case.

A {\em quaternion algebra} $A$ over a field $k$ is a four-dimensional $k$-space with basis vectors $1, \i,\j,\k$ where multiplication is defined on $A$ by requiring that $1$ is a multiplicative identity and that for some $a,b\in k^*=k\setminus\{0\}$
\begin{equation}
\i^2=a\,1,\hskip10pt   \j^2= b\,1,\hskip10pt   \i\j=-\j\i=\k
. \end{equation}
We then extend the multiplication linearly so that $A$ becomes an associative algebra.  This algebra is actually simple (no proper two-sided ideals) and $A$ will be denoted by the Hilbert symbol
\begin{equation} \HS{a}{b}{k}
. \end{equation}
We note that $\k^2 = (\i\j)^2 = -ab$ so any pair of the basis vectors anti-commute. The Hilbert symbol determines a quaternion algebra,  but different Hilbert symbols may determine the same quaternion algebra,  for instance if $z,w\in k^*$,  then
\[\HS{a}{b}{k}  = \HS{az^2}{bw^2}{k} . \]
If $F$ is a field extending $k$,  then
\[ \HS{a}{b}{F} \otimes_F k \approx \HS{a}{b}{k} . \]
The two most common quaternion algebras are Hamilton's quaternions
\[ {\bf H}=\HS{-1}{-1}{\IR} \] 
and the matrix algebra
\[M_{2\times 2}(k) = \HS{1}{1}{k} \]
with generators
$ \i = \Big( \begin{array}{cc} 1 & 0 \\ 0 & -1\end{array}\Big)$ and $ \j=\Big( \begin{array}{cc} 0 & 1 \\ 1 & 0\end{array}\Big)$.

\medskip

The following theorem is the first step to determine the arithmetic criteria we want to develop.

\begin{theorem}  Let $k$ be a
number field and let $A$ be a quaternion algebra over $k$. If $\nu$ is a
place of $k$, let $k_{\nu}$ denote the completion of $k$ at $\nu$. Then
$A$ is said to be ramified at $\nu$ if the quaternion algebra $A \otimes_k
k_{\nu}$
is a division algebra over $k_{\nu}$. Now suppose that $k$ has exactly one
complex place and $A$ is ramified at all real places. Let $\rho$ be an
embedding
of $A$ into $M_2(\IC)$, ${\cal O}$ an order of $A$ and ${\cal O}^1$
the elements of norm 1 in ${\cal O}$. Then $P\rho ({\cal O}^1)$ is a finite
covolume Kleinian group and the totality of arithmetic Kleinian groups
consists of all groups
commensurable with some such $P\rho ({\cal O}^1)$
\end{theorem} 

The field of definition  of an arithmetic Kleinian group is recovered as its
invariant trace field 
\[ k\G = \G^{(2)} = \IQ(\{\beta(g):g\in\G\}) \]
and the quaternion algebra is the invariant
quaternion algebra. One then deduces that two arithmetic Kleinian groups are commensurable up
to conjugacy if their invariant quaternion algebras are isomorphic \cite{MR}.

\medskip

Now we come to the following theorem of \cite{GMMR}.  It makes use of the arithmetic description above by making the explicit calculations and examining their implications.
\begin{theorem}
Let $\G$ be a finitely generated non-elementary subgroup of the group
$\mathrm{PSL}(2,\IC)$ such
that
\begin{itemize}
\item $k\G$ is a number field with exactly one complex place,
\item $\tr(\G)$ consists of algebraic integers,
\item $A\G$ is ramified at all real places of $k\G$.
\end{itemize}
Then $\G$ is a subgroup of an arithmetic Kleinian group.
\end{theorem}
If $\G$ is a non-elementary Kleinian group which contains a parabolic element,
then $A\G$ cannot be a division algebra and so cannot be ramified at any places.

\medskip

As a basic reference to the deep relationships between arithmetic and 
hyperbolic geometry we again refer to \cite{MR}.  We note here a few connections.  
In $Isom^+\IH^2$,
 Takeuchi
\cite{Tak} identified all 82 arithmetic lattices generated by two elements
of finite order (equivalently arithmetic Fuchsian triangle groups), and, 
all arithmetic Fuchsian groups with two-generators have been identified.  The
connections between arithmetic surfaces,  number theory and 
theoretical physics can be found in work of Sarnak and coauthors, e.g. \cite{RS}.  In
\cite{V2} Vinberg gave criteria for Coxeter groups in $Isom(\IH^n)$ to be 
arithmetic. Such groups do not exist in the cocompact case for dimension $n \geq 30$. Further there are only finitely many conjugacy classes of 
maximal arithmetic Coxeter groups for $n \geq 10$ \cite{Nik} (recently established for $n=2,3$  as well).

   Returning to dimension 3, the orientation-preserving subgroups of Coxeter
   groups for tetrahedra, some of which are generated by two elements of finite
order, and which are arithmetic  are identified in \cite{V,Vin} (see  \cite{CM} for related results).
  A. Reid
\cite{Reid} identified the figure eight knot complement as the only 
arithmetic knot complement. 

\medskip

As far as volumes are concerned,  the 
signature formula and Takeuchi's identification of arithmetic triangle groups shows that the smallest  9 lattices of $\mathrm{Isom}^+(\IH^2)$ are arithmetic,  the smallest non-arithmetic lattice is the $(2,3,13)$-triangle group. 
Also the two smallest non-compact lattices are arithmetic  - $(2,3,\infty)$ and $(2,4,\infty)$ with the smallest non-arithmetic lattice being $(2,5,\infty)$. 

\medskip
 
For 3-manifolds we have the  ``expected'' data \cite{MR}  
\begin{itemize}
\item five  smallest closed manifolds arithmetic, but only $20$ of first $50$
\item smallest non-compact manifold arithmetic, but  only $2$ of first $50$  
 \item expect  the  dozen smallest lattices to be arithmetic with perhaps the smallest non-arithmetic  commensurable with a tetrahedral group
 \end{itemize}
 
\noindent  What is surprising in all this is the prevalence of two-generator groups.

\subsection{Two-generator arithmetic Kleinian groups}
To start with, suppose that $\G = \langle f, g\rangle$ where $f$ is a primitive elliptic
 of order $p\geq 3$ and $g$ is elliptic of order $2$. Then,  $\beta=\beta(f)=-4\sin^2(\pi/p)$ is an algebraic integer in the field $ \IQ(\cos 2\pi/p)$, which is totally real with $\phi(p)/2$ places,  where $\phi$ is the Euler function. Let $ R_p$ denote the ring of integers in $ \IQ(\cos (2\pi/p))$,  $R_p=\IZ[\beta]$.  The Galois conjugates of $\beta$ are
\[ \sigma(\beta) = -4 \sin^2(m\pi/p)\]
with $(m,p)=1$. They all lie in the interval $(-4, 0)$ with $\sigma(\beta)<\beta$. Next,  for a number field $K$ of degree
$p$ over $\IQ$, the $p$ Galois embeddings $\sigma : K \to \IC$ give rise to valuations on $K$
which fall into equivalence classes - the places - modulo the action of complex
conjugation.  In \cite{GMMR} we prove the following theorem:
\begin{theorem}\label{arithcharact}
Let $\G = \langle f,g\rangle $ be a subgroup of M\"ob($\IC$)  with $f$ a primitive
elliptic element of order $p\geq 3$, and $g$ an elliptic of order $2$.  Suppose that  $\gamma=\gamma(f,g)\neq \beta(f)=\beta$. Then $\G$ is a subgroup of an arithmetic Kleinian group if
\begin{enumerate}
\item  $\IQ(\gamma,\beta)$ has at most one complex place;
\item $\gamma$ is a root of a monic polynomial $P(z) \in \IZ[\beta][z]$;
\item if $\gamma$ and $\bar\gamma$ are not real, then all other roots of $P(z)$ are real
and lie in the interval $(\beta, 0)$;
\item if $\gamma$ is real, then all roots of $P(z)$  are real and lie in the
interval $(\beta, 0)$.
\end{enumerate}
\end{theorem}
The requirement that the field $\IQ(\gamma,\beta)$ has at most one complex
place can be described in terms of the
polynomial $P(z)$. While if in addition $p = 3, 4$ or $6$, then $\beta$ is an integer and the criteria
admit a greatly simplified description.
 
\medskip

The following theorem is the main result of \cite{McM}.

\begin{theorem}[Finiteness Theorem] \label{finiteness}
There are only finitely many arithmetic Kleinian groups generated by two elements of orders $p$ and $q$,  $2\leq p,q\leq \infty$.
\end{theorem}
In particular there are only finitely many arithmetic generalised triangle groups,  contrary to a conjecture of Hilden and Montesinos.  A generalised triangle group has a presentation of the form
\[ \langle x,y| x^p=y^q=W^r\rangle\]
where $W$ is a cyclically reduced word in $x$ and $y$ (and  involving both $x$ and $y$). When $q=r=2$ and
$W(x, y) = (xy)^3(x^{-1}y)^2$,  these groups are commensurable with the Fibonacci groups $F_{2p}$
and are arithmetic Kleinian groups  if and only if $p = 4, 5, 6, 8,12$. When
$W(x, y) = [x, y]$ we denote the group by $(p, q; r)$. These are arithmetic for $(p,q;r) = (3, 3; 3),
(4, 4; 2), (3,4; 2) $. Other examples, e.g. $(3, 6; 2)$, $(6, 6; 2)$, $(6, 6; 3)$ can be shown to be
nearly arithmetic, but neither arithmetic nor finite extensions of Fuchsian groups.

Determining all the  groups of Theorem \ref{finiteness} is a very interesting,  and not yet complete,  problem.  It involves serious computational number theory,  using work of Stark, Odlyzko,  Diaz Y Diaz and Olivier on discriminant bounds for fields as well as computer searches,  then obtaining covolume bounds and a topological description.

We sketch the proof which will actually apply to a larger class of Kleinian groups. A finitely-generated non-elementary Kleinian group $\G$ is said to be {\em nearly
arithmetic} if 
\begin{itemize}
\item $\G$ is a subgroup of an arithmetic Kleinian group, and
\item $\G$ is not isomorphic to a free product of cyclic groups.
\end{itemize}
Certainly an arithmetic Kleinian group is nearly arithmetic since it is of finite
covolume, thus $\chi(\G)= 0$ and hence $\G$ does not split as the product of cyclic groups.
Within the class of nearly arithmetic groups, those that are arithmetic Kleinian groups
are precisely those with Euler characteristic zero as that distinguishes those with finite
covolume.  We can  say a bit more about nearly arithmetic groups which are not arithmetic.  If the
ordinary set $O(\G)\neq \emptyset$, then Ahlfors
finiteness theorem asserts $O(\G)/\G$ is a finite union of hyperbolic Riemann surfaces of
finite type. It can be proved that $\G$ is rigid and so these Riemann surfaces can support no moduli and therefore is a finite union of triply marked spheres with each component $U$ of $O(\G)$ a round disk.  Thus,  the limit set  $L(\G)$ is a circle or  a Sierpinski curve.
We call this latter class of groups Web groups. Furthermore BersÕ area theorem 
implies that the area of $O(\G)/\G$ is at most $4\pi$. Also there are at most $18$ components $O(\G)/\G$ and hence at most $18$ different
conjugacy classes of triangle subgroups. However the stronger ÒBers Area ConjectureÓ
would imply that there are at most two components of $O(\G)/\G$ and hence two conjugacy
classes of triangle subgroup.

\medskip

To establish the finiteness,  we need a lemma of Schur, \cite{Schur}.

\begin{theorem} If $-1\leq x_1<x_2< \cdots < x_r \leq 1$ with $r\geq 3$,  then
\begin{equation}
\prod_{1\leq i<j\leq r} (x_i-x_j)^2 \leq M_r = \frac{2^23^3 \cdots r^r \times 2^23^3 \cdots (r-2)^{r-2}}{3^3 r^5 \cdots (2r-3)^{2r-3}}
. \end{equation}
\end{theorem}
We set $M_1=1$, $M_2=M_3=4$. We first want to estimate the growth
of the function $M_r$. Of course the numbers $M_{r}^{1/r(r-1)}$ are the squares of the $r^{th}$ order transfinite diameters of the interval $[ -1,1]$ and so $M_{r}^{1/r(r-1)}\to 1/4$. This estimate is enough to yield finiteness but does not provide explicit bounds.  We calculate the discriminant of the field $\IQ(\gamma)$ of degree $n$.  According to our characterisation Theorem \ref{arithcharact} there is a polynomial $p(\gamma)=0$ with $p\in \IZ[\beta][z]$ with exactly one complex conjugate pair of roots ($\gamma$ and $\bar\gamma$) and all other real roots in $(\beta,0)$.  Since the field is ramified at all real places we can write find a polynomial $P(z)\in\IZ[z]$ with roots $\gamma,\bar\gamma$ complex and real roots 
\[-4\sin^2(m\pi/p) \leq x_1,x_2,\ldots x_{n-2} < 0 \]
where $m=p/2$ or $(p-1)/2$ depending on which is even.  Then the discriminant is an integer,  so 
\begin{eqnarray}\label{discrbound}
1 & \leq & |\gamma-\bar\gamma|^2\times \prod_{i=1}^{n-2} |\gamma-x_i|^2 \times \prod_{1\leq i<j\leq n-2} (x_i-x_j)^2
. \end{eqnarray}
As the group does not split as a free product,  $\gamma$ lies inside the ellipse described at Theorem \ref{free}
\begin{itemize}
\item $|\gamma-\gamma|^2 \leq 16(1+\cos(\pi/p)\cos(\pi/q))^4 \leq 256$
\item $|\gamma-x_i|^2 \leq 4+4(1+\cos(\pi/p)\cos(\pi/q))^2 \leq 20$
\end{itemize}
Now we need to use a bit more number theory to get an estimate on the term $\prod_{1\leq i<j\leq n-2} (x_i-x_j)^2
$,  though trivial estimates show finiteness for any fixed $p$ or any fixed $q$. This is because the diameter of the interval containing the points $x_i$ is at most 
\[ 4\sin^2((p-1)\pi/2p)\sin^2((q-1)\pi/2q)\sim 4(1-\pi^2/8p^2)^2(1-\pi^2/8q^2)^2=4s,\]
 $s<1$,  and so when appropriately scaled so the $x_i$ are spread out over an interval of length $2$ so we can use Schur's lemma,  the term 
\[ \prod_{1\leq i<j\leq n-2} (x_i-x_j)^2 \sim s^{2(n-2)(n-3)} \] 
showing the right hand side of (\ref{discrbound}) converging rapidly to zero and giving a contradiction.  However,  this is not good enough to bound $p$ and $q$ yet.  The point is the distribution of the roots $x_i$ are controlled by the real embeddings and so clustered into intervals of the form $[-4\sigma_i\big(\sin^2(\pi/p)\sin^2(\pi/q)\big),0]$,  all but the first containing the same number of points (the degree of the number field $\IQ(\gamma)$ over $\IQ(\cos(2\pi/p),\cos(2\pi/q))$ so that we can iterate the Schur estimate to more carefully obtain a better bound.  This will then establish finiteness,  the details are in \cite{McM}.  However,  even this refinement gives quite poor bounds on $p$, $q$ and the degree of the number field,  and in order to find all the candidates thing get a little more serious.  First there are super-exponentially increasing estimates of the discriminant of number fields of degree $n$ with one complex place,  this is important
results of Stark \cite{Stark} and  Odlyzko \cite{Od},  
as well as results 
concerning the discriminants of number fields
of small degree such as those of Diaz Y Diaz and Olivier 
\cite{Di,DO}. Next,  a little computation improves the estimates on the discriminant by shrinking the regions that $\gamma$ lies in and so forth.  Presently we have estimates to show the total degree is not more than $17$,  and as soon as $p$ and $q$ are large enough we know all the examples.    To identify all the groups we use  a
computer search to examine all algebraic integers in the field 
satisfying the given bounds and additional arithmetic restrictions on the real embeddings.  This procedure gives us a
relatively short list of candidate discrete groups which are now known to 
be subgroups of arithmetic Kleinian groups as per our results above.     

\subsection{The arithmetic groups generated by two parabolics} 

Let us work through the easiest possible case, \cite{GMcM}.  It contains all the ideas for more general cases,  but none of the complexity because the existence of a parabolic element  implies {\em a priori} a good bound on the associated number field.  Every group generated by two parabolics admits two (not necessarily different) $\IZ_2$-extensions,  so it make sense to find these (it also makes estimates better).  We are therefore looking for Kleinian groups whose complex parameters are of the form $(\gamma,0,-4)$.  The picture of this space is above in \S \ref{Riley}.

\subsubsection{Degree bounds}

If an arithmetic group $\G\subset \mathrm{PSL}(2,\IC)$ contains a parabolic element $f$,  then the invariant quaternion algebra cannot be a division algebra as $f^2-I$ is nilpotent,  and Wedderburn's structure theorem implies this quaternion algebra must be the matrix algebra $M_{2\times 2}(k\G)$.  But this algebra is not ramified at any places,  so $k\G$ is either $\IQ$ or $\IQ(\sqrt{-d})$.  $k\G=\IQ$ gives groups commensurable with the modular group so they cannot have finite covolume and hence cannot be arithmetic (in $3$-dimensions).  This gives us the following {\em a priori} bound on the degree of the field $\IQ(\gamma)$.

\begin{theorem}\label{quad}  If $\G$ is an arithmetic Kleinian group which contains a parabolic element,  then $\G$ is commensurable with a Bianchi group,  $PSL(2,{\cal O}_d)$.  In particular the invariant trace field $\IQ(\gamma)=\IQ(\sqrt{-d})$ and ${\cal O}_d$ is the ring of integers.
\end{theorem}

Ordinarily we would have to find bounds on the degree of $\IQ(\gamma)$ using arguments such as we discussed in the proof of the finiteness theorem.

\subsubsection{List of possible $\gamma$ values}

We know that $\gamma\in \IC\setminus \IR $ is an integer in a quadratic extension of $\IQ$. Thus
\begin{equation}\label{gamminpoly}
\gamma^2 + b \gamma + c = 0, \hskip20pt b, \; c \in \IZ
. \end{equation}
We know the group with parameters $(\gamma, 0, -4)$ cannot be freely generated and so we have $b=\gamma+\bar\gamma = 2\Re e\{\gamma\} \in (-8,8)$ and $c=|\gamma|^2 < 16$ from Theorem \ref{free}.
Thus $c\in\{\pm1,\pm 2,\pm 3\}$ and $b\in \{0,\pm 1,\pm 2, \pm 3, \pm 4,\pm 5, \pm 6,\pm 7\}$.  The symmetries of this space (in this case $\pm\gamma$, $\pm \bar\gamma$) all give isomorphic groups and so  we need only consider $-7\leq b\leq 0$ and we also need $b^2-4c< 0$.  This eliminates $b=-7,-6,-5,-4$,  so we are down to less than 20 cases.  Refining the condition for freeness using fairly elementary methods,  such as that of Lyndon and Ullman, \cite{LU},  quickly get the list down to $7$ points:
\[ i,i\sqrt{2},i\sqrt{3}, (1+i\sqrt{3})/2, (1+i\sqrt{7})/2,  1+i,   (3+i\sqrt{3})/2. \]
The first three cases give known examples of generalised triangle groups,  or subgroups of infinite-index in reflection groups and so cannot be arithmetic.  The remainder are arithmetic,  but we need a description.

\subsubsection{The orbifolds}

The way we find these four arithmetic groups is simply by compiling a (long) list of $\gamma$ values that do come from discrete groups by direct computation using J. Weeks' snappea program.  Thus we set up tables of two bridge knot and link complements (two-bridge implies two-generator) and ask for matrix representations in $\mathrm{PSL}(2,\IC)$ and from these we calculate various $\gamma$ values,  such as those for the group or its $\IZ_2$-extensions and so forth.  In the case of torsion,  we first perform orbifold Dehn surgery on the knot or link complement.   Once we numerically compare $\gamma$ values between these lists it is easy to check directly that it is correct and we have the right group.  We then obtain theorems of the following sort after passing to an index-two subgroup to get our group generated by two parabolics.
\begin{itemize}
\item $\frac{1+i\sqrt{3}}{2}$: This gives the figure of $8$ knot group, a two-bridge link of slope $\Big(\frac{5}{3}\Big)$.  It has index $12$ in $PSL(2,{\cal O}_{3})$.  It is the only arithmetic knot,  \cite{Reid}.
\item $1+i$: This gives the Whitehead link group. It is a two-bridge link of slope $\Big(\frac{8}{3}\Big)$and has index $12$ in the Picard group $PSL(2,{\cal O}_{1})$.
\item $\frac{1+i\sqrt{3}}{2}$: This gives the two-bridge link of slope $\Big(\frac{10}{3}\Big)$ and has index $24$ in $PSL(2,{\cal O}_{3})$. 
\item $\frac{1+i\sqrt{3}}{2}$: This gives the two-bridge link of slope $\Big(\frac{12}{5}\Big)$. 
\end{itemize}
  
 \scalebox{0.6}{\includegraphics[viewport= 30 550 480 780]{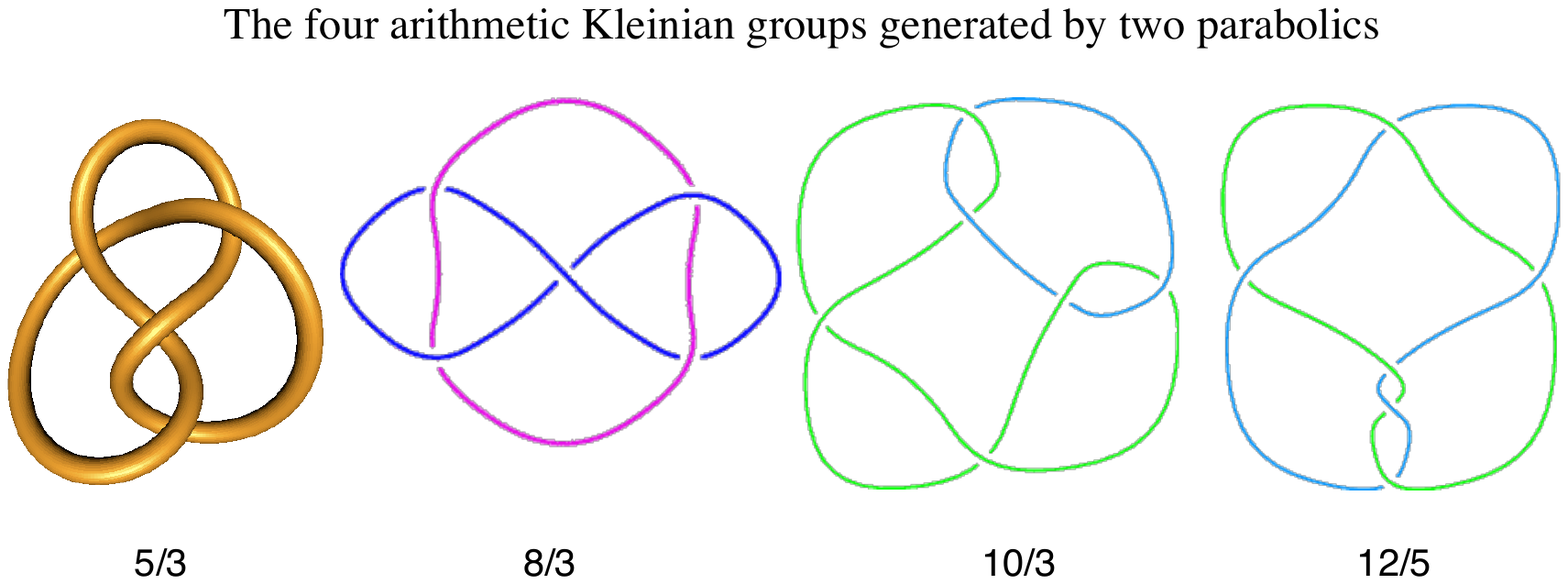}}
 
 Surprisingly there are no orbifolds here.  In summary: 
 \begin{theorem} Let $\G$ be an arithmetic Kleinian group generated by two parabolic transformations.  Then $\G$ is torsion-free and $\IH^3/\G$ is homeomorphic to one of the four two-bridge knots or links described above.
 \end{theorem}
 
\subsection{The arithmetic generalised triangle groups $p,q\geq 6$}

Building on the above ideas we have the following theorem when the two-generators have order six or more.

\begin{theorem}  Let  $\Gamma=\langle f,g \rangle$ be an arithmetic 
Kleinian group generated by elements of
order $p$ and $q$ with $p,q \geq 6$.  Then $p,q$ fall into one of the following 
4 cases:
\begin{enumerate}
\item $p=q=6$, \\ there are precisely  $12$ groups enumerated below
(See comments following the table).
\item $p=q=8$, \\ there is precisely one group obtained by (8,0) 
surgery on the knot 5/3.
\item $p=q=10$, \\there is precisely one group obtained by (10,0) 
surgery on each component of the link 13/5.
\item $p=q=12$, \\ there are precisely two groups obtained by (12,0) 
surgery on the knot 5/3  and on each component
of 8/3. 
\end{enumerate}
\end{theorem}
Here $r/s$ denotes the slope (or Schubert normal form) of a  two 
bridge knot or link.
\bigskip

  As a corollary we are able to verify a condition noticed by  Hilden, 
Lozano and Montesinos
\cite{HLM1} concerning the $(n,0)$ surgeries on two-bridge link complements.
\begin{corollary}  Let $(r/s,n)$ denote the arithmetic hyperbolic orbifold 
whose underlying space is the 3-sphere and whose 
singular set is the two-bridge knot or link with
slope
$r/s$ and has degree $n$.  Then
\begin{equation}
n \in \{2,3,4,5,6,8,10,12,\infty \}
. \end{equation}
\end{corollary}

We also have the following,  slightly surprising, corollary.

\begin{corollary}\label{notrigrp}  There are no arithmetic generalised triangle 
groups with generators of orders at
least  $6$.
\end{corollary}

\begin{table}
\centering
  \caption[]{Arithmetic groups with $p,q=6$ }
\begin{tabular}{|c|c|c|}
\hline
  $\gamma$ value & $k\Gamma$ & description of orbifold \\
\hline
   $ i\sqrt{3} $ &$z^2-z+1$ & $\Gamma_{21}$
\\
\hline
  $-1+i$ & $z^2+1$ & (6,0) surgery on 5/3 \\
\hline
  $-1$ & $z^2+z+1$ & $\Gamma_{20}$\\
\hline
   $1+3i$  &$z^2-2z+2$ & (6,0)-(6,0) surgery on link 24/7
\\
\hline
  $-1+i\sqrt{7}$ & $z^2-z+2$ & (6,0)-(6,0) surgery on link 30/11
\\
\hline
   $-2+i\sqrt{2}$  & $z^2+2$ & (6,0)-(6,0) surgery on link 12/5
\\
\hline
  $4.1096 - i\ 2.4317$ & $1+2z-3z^2+z^3$ &
(6,0) surgery on knot 65/51 \\
\hline
    $3.0674 -i\ 2.3277$ & $2-2z^2+z^3$ & (6,0) surgery on knot 13/3$^*$ \\
\hline
   $2.1244 -i\ 2.7466$ & $1+z-2z^2+z^3$ &  (6,0) surgery on knot 15/11 \\
\hline
    $1.0925 - i\ 2.052 $ & $1-z^2+z^3$ & (6,0) surgery on knot 7/3$^*$ \\
\hline
  $0.1240 -i\ 2.8365$ & $1+z-z^2+z^3$ & (6,0) surgery on knot 13/3$^*$\\
\hline
  $-0.8916 -i\ 1.9540 $ & $1+z+z^3$ & (6,0)-(6,0) surgery on link 8/5 \\
\hline
   $-1.8774 - i\ 0.7448$ & $1+2z+z^2+z^3$ & (6,0) surgery on knot 7/3$^*$ \\
\hline
    $-2.8846 -i\ 0.5897$ & $1+3z+z^2+z^3$ &
(6,0)-(6,0) surgery on link 20/9 \\
\hline
\end{tabular}
\end{table}
\bigskip

As indicated for the case of two parabolic generators the results above were produced as follows.  The 
methods we outlined above on bounding degree and discriminant for $\IQ(\gamma)$  produce for us all possible values of the trace 
of the commutator of a pair of primitive elliptic generators of an arithmetic Kleinian group.  Refining the description of the ``free on generators'' part of this space brings this list down to manageable size. A computer programme gives us a  Dirichlet region and an approximate volume (recall that at this point discreteness of the group is assured by the number theory).   We then use Jeff Weeks' hyperbolic geometry 
package snappea  to try and
identify the orbifold in question by surgering various two-bridge 
knots and links and comparing
volumes.  Once we have a likely candidate,  we use the matrix presentation 
given by snappea and verify that
the commutator traces are the same.  As these traces come as the 
roots of a monic polynomial with integer coefficients of modest
degree,  this comparison is exact.  Since this trace determines the 
group up to conjugacy,  we   identify the orbit space.
Conversely, once the two-bridge knot or link and the relevant surgery 
is determined, a value of $\gamma$ can be recovered from the algorithm
in \cite{HLM1}.

\medskip

We now turn to another family of examples.

\subsection{The non-compact arithmetic generalised triangle groups}

Here we will describe all $21$ of the  non-compact arithmetic Kleinian groups with two elliptic generators. The non-compact arithmetic generalised  triangle groups  are necessarily included and 15 of the 21 groups obtained turn out to be 
generalised triangle groups, the remainder being quotients of generalised
triangle groups. The generalised triangle groups all yield orbifolds 
with a particularly simple singular set in $\IS^3$.
 
 \begin{table}
\begin{center}
\begin{tabular}{||c|c|c|c|c|l|c||}  \hline\hline
$\Gamma_i$ & p & q & $\gamma$ & d & $\Gamma_i \cong \langle x^p=y^q= \cdots
=1 \rangle $ & covolume \\ \hline\hline
1 & 2 & 3 & $(-3 + \sqrt{3}i)/2$ & 3 & $((xy^{-1})^2(x^{-1}y)^2)^3$ &
0.33831 \\ \hline
2 & 2 & 3 & $(-1 + \sqrt{3}i)/2$ & 3 & $((xy^{-1})^3(x^{-1}y)^3)^3$ &
0.67664 \\ \hline
3 & 2 & 3 & $(-5 + \sqrt{3}i)/2$ & 3 & $(xyx^{-1}yxy^{-1})^6$& 0.67664 \\
\hline
4 & 2 & 4 & $-1+i$ & 1 & $((xy^{-1})^2(x^{-1}y)^2)^2$ & 0.45798 \\ \hline
5 & 2 & 4 & $-2+i$ & 1 & $(xyx^{-1}yxy^{-1})^4$ & 0.45798 \\ \hline
6 & 2 & 4 & $i$ & 1 & $((xy^{-1})^3)(x^{-1}y)^3)^2)$ & 0.45798 \\ \hline
7 & 2 & 6 & $(-1+\sqrt{3}i)/2$ & 3 &
$\begin{array}{l}((yx)^2(y^{-1}x)^2)^3\\  =((yx)^2(y^{-1}x)^2y)^2 \\
 = (xy^{-1}x(yx)^2(y^{-1}x)^2)^2 \\
 =(xy^{-1}x(yx)^2(y^{-1}x)^2y)^2\end{array}$ & 0.25374 \\ \hline
8 & 2 & 6 & $(1+\sqrt{3}i)/2$ & 3 &
$\begin{array}{l}([y,x](yx)^2[y^{-1},x]y)^2 \\
 = ([y,x](yx)^2 [y^{-1},x])^2 x y^{-2} \; \cdot \\
\hskip40pt [x,y][x,y^{-1}](xy^{-1})^2x \end{array}$
& 0.50748 \\ \hline
9 & 2 & 6 & $(-3+\sqrt{3}i)/2$ & 3 & $(xyx^{-1}yxy^{-1})^3$ & 1.01496 \\ \hline
10 & 3 & 3 &$-3$ & 3 & $[x,y]^3$ & 0.67664  \\ \hline
11& 3 & 3 & $-2+\sqrt{3}i$ & 3 & $(xyx^{-1}yxy^{-1})^3$ & 1.35328 \\ \hline
12 & 3 & 4 & $-2$ & 1 & $[x,y]^2$ & 0.30532 \\ \hline
13 & 3 & 4 & $2i$ & 1 & $((xy^{-1})^2(x^{-1}y)^2)^2$ & 1.5266 \\ \hline
14 & 3 & 6 & $-1$ & 3 & $\begin{array}{l} [x,y]^3=([x,y]x)^2\\
=(y^{-1}[x,y])^2=(y^{-1}[x,y]x)^2\end{array}$ &
0.08458
\\ \hline
15 & 3 & 6 &$-1+\sqrt{3}i$ & 3 & $(xyx^{-1}yxy^{-1})^2$ & 1.6916 \\ \hline
16 & 3 & 6 & $2+\sqrt{3}i$ & 3 &$((xy)^2xy^{-1}(x^{-1}y^{-1})^2)^2$ &
1.35328 \\ \hline
17 & 4 & 4 & $-2$ & 1 & $[x,y]^4$ & 0.91596 \\ \hline
18 & 4 & 4 & $-1+2i$ & 1 & $(xyx^{-1}yxy^{-1})^2$ & 1.83192 \\ \hline
19 & 4 & 6 & $-1$ & 3 & $\begin{array}{l} [x,y]^3=([x,y]x)^2 \\
=(y^{-1}[x,y])^2=(y^{-1}[x,y]x)^2 \end{array}$ &
0.21145
\\ \hline
20 & 6 & 6 & $-1$ & 3 & $\begin{array}{l}[x,y]^3=([x,y]x)^2\\
=(y^{-1}[x,y])^2=(y^{-1}[x,y]x)^2 \end{array} $
& 0.50748 \\  \hline
21 & 6 & 6 & $\sqrt{3}i$ & 3 &
$\begin{array}{l}
(y^{-1}x)^2y[x^{-1},y][x,y][x,y^{-1}]x^{-1} \\ = ([y^{-1},x]yx^2)^2
\end{array}$ &
1.01496
\\
\hline\hline
\end{tabular}
\caption{All non-compact arithmetic groups with 2 elliptic generators}
\end{center}
\end{table}

Recall that 
if  $\G$ is a non-compact arithmetic Kleinian group,  then there is a parabolic element in $\G$ and,  as we outlined above at Theorem \ref{quad},  $kG =
\IQ(\sqrt{-d})$
is a quadratic imaginary number field, $AG \cong M_2(\IQ(\sqrt{-d}))$ and
$G$ is commensurable, up to conjugacy, with a Bianchi group $PSL(2,O_d)$. 
Apart from the  $15$  generalised triangle groups there are three Coxeter groups
and the others are $\IZ_2$-extensions of the two $p=q=6$ groups. None of
these appear to be
generalised triangle groups.

We have also given an approximation to the
covolume of these
groups. These have already been computed for the generalised triangle groups
$\Gamma(p,q;r)$ described here in \cite{HMR} and is also well known for the
Coxeter  groups. The remaining covolumes can all be computed from the
information above
and the following additional calculations, made using MAGMA. (Recall that
the conjugacy class
of a finite covolume group is determined by its isomorphism class)
\begin{itemize}
\item $\tilde{\Gamma}(2,3;6)$ is a subgroup of index 8 in $PGL(2,O_3)$.
\item $\tilde{\Gamma}(2,6;3)$ is a subgroup of index 12 in $PGL(2,O_3)$.
\item $\tilde{\Gamma}(3,6;2)$ is a subgroup of index 20 in $PGL(2,O_3)$.
\item $\tilde{\Gamma}(2,4;4)$ is a subgroup of index 6 in $PGL(2,O_1)$.
\end{itemize}

Although we have listed 21 groups and we know that 2 of these groups are
the same, we have not yet proved that all
the others are pairwise non-conjugate. For many,
this follows immediately by consideration of the covolumes. For pairs that
have the same covolume,
particularly pairs that are $\IZ_2$-extensions of the same underlying
group, we use the fact that
the trace field of the group itself, which can be determined from the
normalisation
  is a conjugacy invariant. By computing these trace fields, as required,
we can  show that all other pairs are indeed non-conjugate.

\subsubsection{The Orbifolds}

Certain orbifolds whose singular sets are graphs consisting of knots or links
with tunnels were studied in \cite{H,HMV,JR}, and the fundamental groups
of these orbifolds are generalised triangle groups. The quotient spaces of the
15 generalised triangle groups admit such a description. Indeed they all admit a uniform type
of description. Consider the orbifold whose singular set in $\IS^3$ is given by the
graph below, 
where the label $i$ indicates a cone angle of $2 \pi/i$ on that
edge of the singular set.

\medskip
 
\begin{center}
 \includegraphics[scale=0.6]{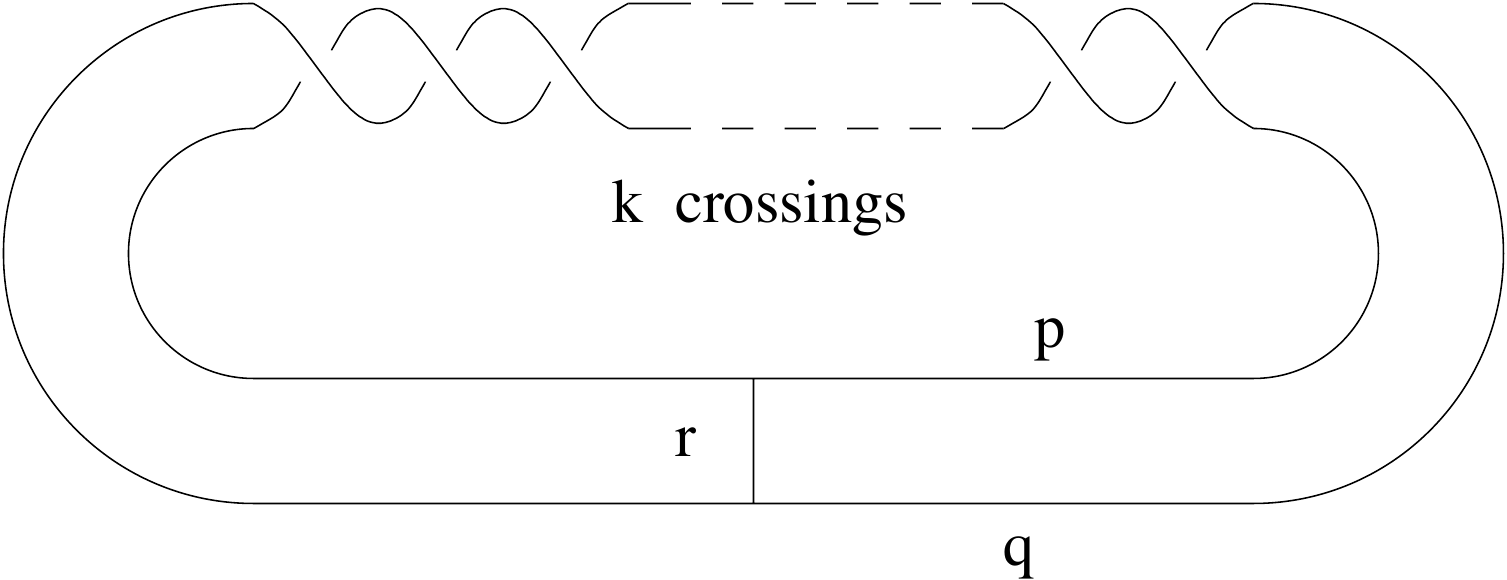}\\
 $(p,q;r)$--orbifolds 
\end{center}

The fundamental group of this orbifold has presentation
$$ \langle x,y \mid x^p=y^q=W(x,y)^r = 1 \rangle$$
where $W(x,y)$ has $2k$ letters, $k$ with exponent $+1$ and $k$ with exponent
$-1$ and 
$$  W(x,y) = \left\{ \begin{array}{ll}
                   xyxy \ldots xyx^{-1}y^{-1} \ldots x^{-1}y^{-1} & {\rm if}~k~{\rm
is~even} \\
                  xyxy \ldots yxy^{-1}x^{-1} \ldots x^{-1}y^{-1} & {\rm if}~k~{\rm
is~odd}
                   \end{array}
            \right.
$$
(c.f \cite{JR}). It is not difficult to see that all 15 of our generalised triangle groups have such presentations. When $k$ is even, these orbifolds were discussed
in \cite{JR,H}; when $k = 3$ in \cite{HMV} and when $k=5$ in \cite{H}.

\subsection{The arithmetic points of the $(2,3)$ commutator plane}

There remains the cases of identifying all the two-generator arithmetic groups generated by two elliptic elements of low order.  We are currently working on these problems but they become far more computational. As an instance here is some data from H. Cooper's PhD Thesis, \cite{Cooper} which identifies all the arithmetic Kleinian groups generated by elliptics of order $2$ and $3$.  It utilises earlier work of Q. Zhang \cite{Zhang} and V. Flammang- G. Rhin \cite{FR}. We presently also have tables for the $(2,4)$ and $(2,5)$ cases as well.

To see how these points fit in the scheme of things,  one should look at the illustration of the $(2,3)$-commutator plane below in the next section \S \ref{23plane}.

\begin{table}
\begin{center}
\begin{tabular}{|c|c|c|c|}  \hline 
$\G_i$ & Minimal polynomial for $\gamma$ & $\gamma$ & covolume \\ \hline \hline
1. & $z^4 +4z^3 +4z^2 +z+1 $&$0.00755 + 0.51312i $&$ 0.3271 $\\  \hline
2. & $z^5 +6z^4 +13z^3 +14z^2 +8z+1$&$ -0.66222 + 0.89978i$&$ 0.5158  $\\ \hline
3.& $z^5 +8z^4 +23z^3 +27z^2 +11z+1$&$ -0.18484 + 0.71242i $&$ 0.5174  $\\  \hline
4.&$z^4 +6z^3 +13z^2 +10z+1$&$ -0.79289 + 0.97832i$&$ 0.6860  $\\ \hline
5.&$z^4 +8z^3 +21z^2 +19z+4$&$0.21021 + 0.41375i$&$ 0.4638 $ \\ \hline
6.&$z^4 +6z^3 +13z^2 +11z+1$&$ -1.29342 + 1.00144i $&$  0.4639 $\\ \hline
7.&$z^5 +9z^4 +29z^3 +38z^2 +16z+1$&$0.02268 + 0.62320i $&$ 0.6462  $\\ \hline
8.&$z^5 +6z^4 +12z^3 +11z^2 +7z+1$&$-0.37053 + 0.84016i $&$ 0.6470 $\\  \hline
9.&$z^3 +2z^2 +z+1$&$-0.12256 + 0.74486i $&$ 0.7859 $ \\ \hline
10.&$z^4 +4z^3 +6z^2 +5z+1$&$-0.75187 + 1.03398i $&$ 1.1451 $\\ \hline
11.&$z^5 +5z^4 +6z^3 -z^2 -z+1$&$0.34815 + 0.31570i $&$ 0.5671 $\\ \hline
12.&$z^5 +7z^4 +19z^3 +25z^2 +14z+1$&$-1.35087 + 1.05848i $&$ 0.5672  $\\ \hline
13.&$z^5 +7z^4 +19z^3 +23z^2 +10z+1$&$-1.02127 + 1.12212i $&$ 1.8946 $\\ \hline
14.&$z^4 +5z^3 +8z^2 +6z+3$&$-0.34861 + 0.75874i $&$  0.4325  $\\ \hline
15.&$2z^2 + 2z + 2$&$ -1.00000 + 1.00000i $&$ 0.6105  $\\ \hline
16.&$z^5 +7z^4 +18z^3 +23z^2 +17z+5$&$-0.60186 + 0.93867i $&$ 0.7111 $\\  \hline
17.&$z^5 +5z^4 +7z^3 +3z^2 +2z+1$&$0.11005 + 0.57190i $&$ 0.7125 $\\ \hline
18.&$z^5 +8z^4 +25z^3 +37z^2 +23z+3$&$-1.86240 + 1.07589i $&$  0.8659  $\\ \hline
19.&$z^3 +3z^2 +2z+2$&$-0.23931 + 0.85787i $&$ 1.2552 $\\  \hline
20.&$z^3 +4z^2 +6z+2$&$-1.22816 + 1.11514i $&$ 1.0590  $\\  \hline
21.&$z^3 +4z^2 +5z+4$&$-0.65219 + 1.02885i$&$ 1.4988 $\\  \hline
22.&$z^5 +6z^4 +11z^3 +6z^2 +2z+3$&$ 0.17229 + 0.58559i $&$ 1.1813 $ \\  \hline
23.&$z^5 +6z^4 +12z^3 +11z^2 +8z+3$&$-0.13972 + 0.82586i $&$ 1.7505  $\\  \hline
24.&$z^4 +3z^3 -2z+1$&$0.46746 + 0.27759i $&$ 1.0197  $\\ \hline
25.&$z^4 +6z^3 +12z^2 +7z+1$&$-0.41847 + 0.93916i $&$ 1.0188  $ \\  \hline
26.&$z^3 +3z^2 +4z+1$&$-1.34116 + 1.16154i $&$ 1.3193 $\\  \hline
27.&$z^3 +2z^2 +1$&$0.10278 + 0.66546i $&$ 1.709 $\\  \hline
28.&$z^4 +3z^3 +z^2 -z+1$&$ 0.33909 + 0.44663i $&$ 1.4729 $ \\ \hline
29.&$z^4 +6z^3 +12z^2 +9z+1$&$-1.50000 + 0.60666i $&$ 0.0390 $\\ \hline
30.&$z^4 +5z^3 +7z^2 +3z+1$&$-0.21190 + 0.40136i $&$ 0.0408  $ \\  \hline
31.&$z^3 +4z^2 +4z+2$&$-0.58036 + 0.60629i $&$  0.1323  $ \\  \hline
32.&$z^3 +5z^2 +8z+5$&$ -1.12256 + 0.74486i $&$  0.1571  $\\  \hline
33.&$z^3 +3z^2 +2z+1$&$-0.33764 + 0.56228i $&$ 0.1569 $ \\  \hline
34.&$z^4 +5z^3 +8z^2 +6z+1$&$-0.92362 + 0.81470i $&$ 0.1274 $ \\ \hline 
35.&$z^2 + 3z^ + 3$&$-1.50000 + 0.86603i  $&$ 0.3383 $\\ \hline
36.&$z^3 +4z^2 +5z+3$&$-0.76721 + 0.79255i $&$ 0.2630 7 $\\ \hline
37.&$z^4 +5z^3 +6z^2 +1$&$0.06115 + 0.38830i $&$0.1026  $\\ \hline
38.&$z^6 +8z^5 +24z^4 +35 z^3 +28z^2+12z+1$&$-0.52842 + 0.78122i  $&$ 0.1026  $\\ \hline
39.&$z^3 +3z^2 +z+2$&$ -0.11535 + 0.58974i$&$  0.3307 $ \\ \hline
 \end{tabular}
\caption{The arithmetic groups of the $(2,3)$-commutator plane}
\end{center}
\end{table}
  
 \newpage
   
\section{Spaces of discrete groups,  $p,q\in\{3,4,5\}$}

Our arguments concerning torsion of order $6$ and identifying the Kleinian group of minimal covolume that contains such an elliptic,  Theorem \ref{ST1}, and also the case of simple torsion Theorem \ref{ST2} (which actually relies on the description we find here),  allows us to focus on Kleinian groups in which every elliptic element has order $2$, $3$, $4$ or $5$.  Here we deal with the case that there is torsion of order $3,4$ or $5$.  The case that there is only $2$-torsion is of a completely different nature.
 
 \medskip
 
 In order to get the information we need, for each pair $p,q\in \{2,3,4,5\}$, $p,q$ not both equal to $2$, we want to describe the space
 \[ {\cal D}_{p,q} = \{\gamma\in \IC : (\gamma, -4\sin^2\frac{\pi}{p},\sin^2 \frac{\pi}{q})\hskip5pt \mbox{is Kleinian}\}. \]
In particular the parameters for a Kleinian group which is not free on its generators will be of most interest to us.  This space is illustrated below for perhaps the most interesting case $p=2$ and $q=3$ and corresponds to the slice $(\gamma, -3,-4)$.

\subsection{The $(2,3)$-commutator plane. }\label{23plane}
\begin{center}  $\gamma = \frac{1}{3} \sinh^2(\delta+i\theta)$ \end{center}
 
\scalebox{0.5}{\includegraphics[viewport= -10 300 480 780]{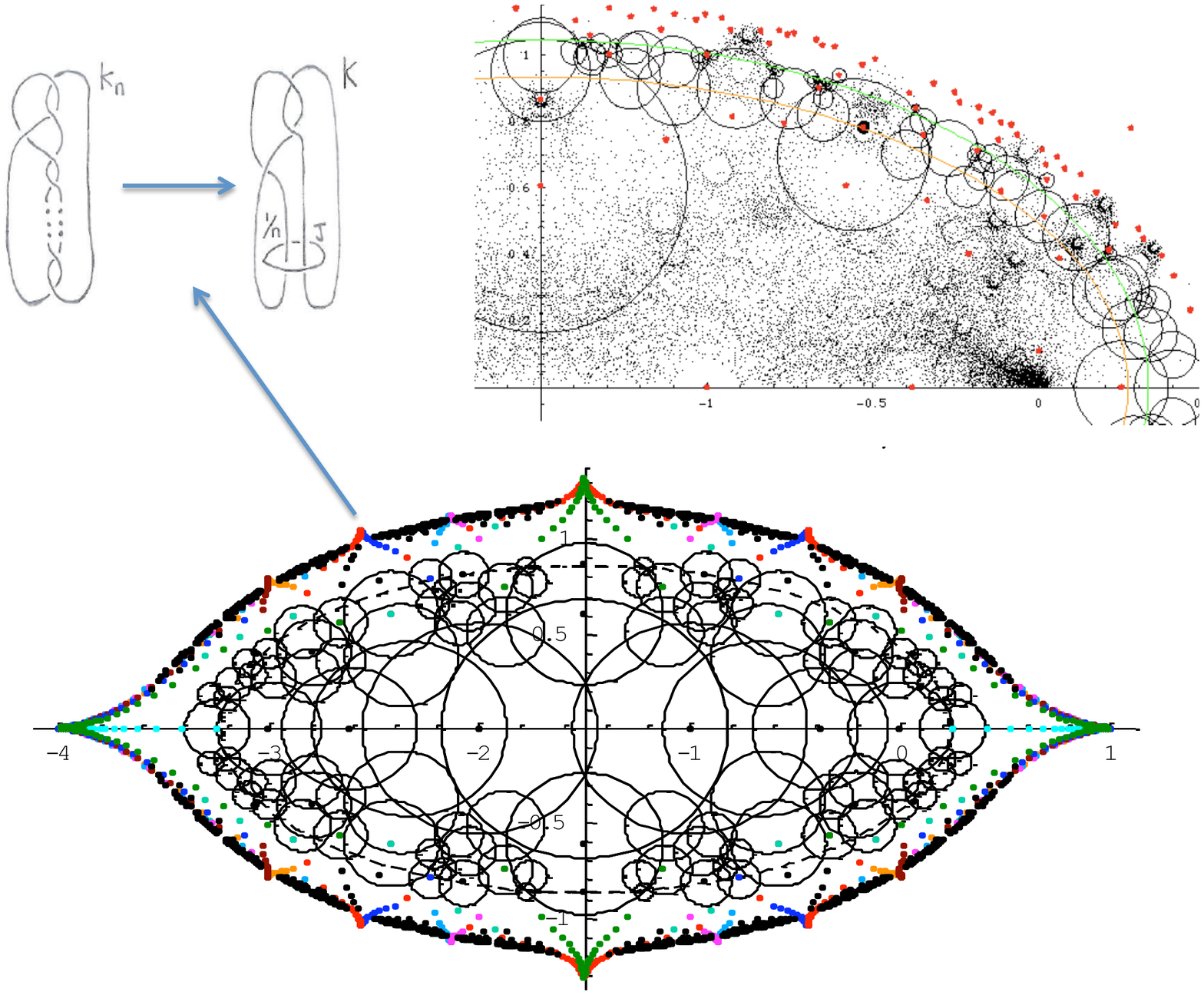}}

 We want to take a moment to describe this picture.  For our application here the most important thing is in identifying the interior points corresponding to Kleinian groups (illustrated in the lower portion) and we will deal with that more carefully in a moment.  In the lower picture we see sequences of points clustering to the boundary.  These points are obtained from $(2,0)$,  $(3,0)$ orbifold Dehn surgeries on various two-bridge links or on $\IZ_2$-extensions of $(3,0)$ orbifold Dehn surgeries on two-bridge knots or links.  The way these sequences cluster and their density in the boundary is interesting,  but will be investigated elsewhere.  We have illustrated a sequence of two-bridge links,  $(3,0)$ orbifold Dehn surgered and then $\IZ_2$-extended, which converge via increasing twists to the limit link as per Thurston's limit theorem - producing a parabolic.  On the upper-right we illustrate the ``dust'' obtained from calculating the roots of about 500 of the good word polynomials with $\beta=-3$. The larger points are actually those of discrete groups,  there are obviously rather fewer of these.  
 
 Each disk in the interior region in the lower picture corresponds to an inequality of the form ``there are no $\gamma$ values for a Kleinian group within the disk $|\gamma-\gamma_0|<r_0$ apart from a finite number exceptions and these are ... ''. We will give a more precise version of this in a moment.

 \medskip
 
 Our applications require us to find the spectrum of distances between the axes of elliptics.  Notice that (\ref{cosh2r}) implies that $\gamma$ values on the the ellipse with foci $-3$ and $0$ have constant axial distance.  One such ellipse is illustrated.  In general the foci of theses ellipses will be $-4\sin^2(\pi/p)\sin^2(\pi/q)$ and $0$.  In fact the paramerisation given by $\gamma = \frac{1}{3} \sinh^2(\delta+i\theta)$ is often a very good way to study these spaces.
 
 We have earlier proven that ${\cal D}_{p,q}$ is closed, Theorems \ref{closed1} \& \ref{closed2} while obviously $0\not\in {\cal D}_{p,q}$. Next,   we recall that the family of good word polynomials is closed under composition.  Thus we may study the iterates of a specific polynomial
and, in this way, generate a sequence of commutator traces which
cannot accumulate at $0$.  In the particular case of order $3$,  if $w = (fgfg^{-1})^3$, then 
\[ p_w(\gamma,  \beta ) = \gamma(\gamma-\beta)(1 +  \beta-\gamma)^2(3 +  \beta-\gamma)^2\]
with $\beta=-3$ we write 
 \[ p_w(z) = \gamma^3(\gamma+3)(2+\gamma)^2. \]
  Thus $p_w$ has a super-attracting fixed point at the origin.  The preimages of $0$ are $0,-2,-3$.  The disk $\ID(0,r_0)$ with $r_{0}^2(r_0+3)(r_0+2)^2=1$,  $r_0 = 2\cos(2\pi/7)-1 = 0.24698\ldots$,  though this is only a small part of the super-attracting basin.  Indeed the filled in Julia set is illustrated below. This connection between the parameters for discrete groups and the Julia sets for polynomials arising from the family of good words is very interesting and yet to be fully explored.  Notice that $z_0=2\cos(2\pi/7)-1$ is a repelling fixed point for the polynomial iterations $z\mapsto z^3(z+2)^2(z+3)$ and so in the boundary of the attracting basin.
  
  \scalebox{0.6}{\includegraphics[viewport= 25 460 500 790]{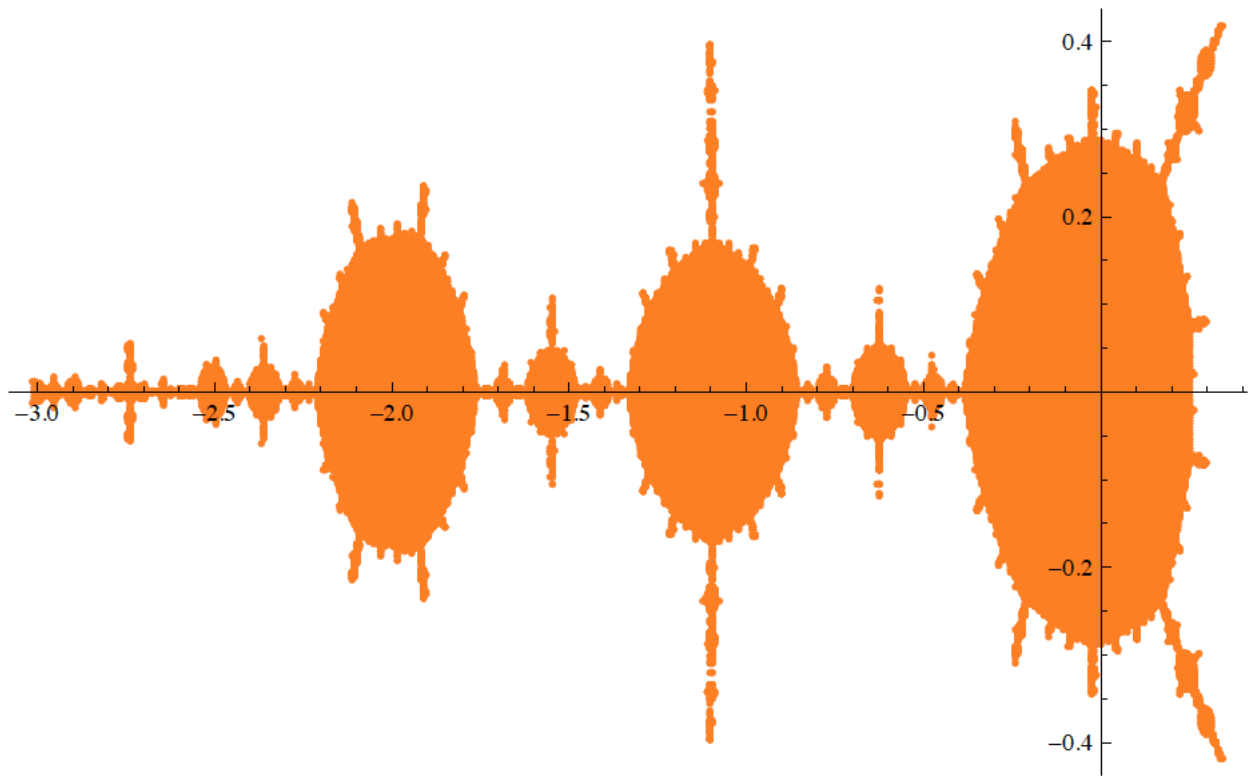}}
  
  \begin{center}
  Filled in Julia set for $z\mapsto z^3(z+2)^2(z+3)$.
  \end{center}
  
  However this gives us the following sharp result.
  \begin{lemma}  Let $\langle f,g\rangle$ be a Kleinian group with  $f$ an elliptic of order $3$.  Then 
  \begin{equation}
  |\gamma|\geq 2\cos(2\pi/7)-1
  . \end{equation}
  This estimate is sharp and uniquely achieved in the $(2,3,7)$--triangle group with presentation
  \[ \langle f,g: f^3=g^7=(fg)^2=1\rangle. \]
  \end{lemma}
   Next,  by  considering the polynomial $p_w(\gamma,\beta)=\gamma(\gamma^2-(\beta-1)\gamma-(\beta-1))^2$ associated with the word $\tilde{w}=gf^{-1}g^{-1}fgfg^{-1}f^{-1}g$,  set $\beta=-3$ to get $p_{\tilde{w}}(z)= z(z+2)^4$,  and follow the above arguments we obtain
   \begin{equation}
   |\gamma(2+\gamma)^4|\geq 2\cos(2\pi/7)-1
   . \end{equation}
   This offers us the chance to find an omitted disk about $-2$.  Namely $\ID(-2,r_1)$ where
   $(2+r_1)r_{1}^{4} = 2\cos(2\pi/7)-1 $, $r_1=0.557459\ldots$.  Notice that we have not eliminated the value $-2$.  The symmetries of the space ${\cal D}_{2,3}$ also tells us that the disk  centred on $-1$ and of the same radius $\ID(-1,r_1)$ is also omitted.  These two disks are illustrated in the picture of ${\cal D}_{2,3}$ above.  More polynomials give the additional disks and following this procedure we begin to cover the region we wish to describe from the inside.  Two points to note are that we never eliminate a root of a polynomial directly,  they are typically eliminated by appealing to their values under other polynomials.  Also,  this way we build up a list of discrete groups.  Some of these polynomials are not very pretty,  for instance the word
   \[w=gfg^{-1}fgf^{-1}g^{-1}f^{-1}gfg^{-1} \]
   gives the polynomial
   \[
   \gamma(\gamma-\beta)(-1 + 2\beta + 2\gamma\beta - \gamma\beta^2 - 2\gamma^2 + \beta\gamma^2 - \beta^2\gamma^2 + 
   2\beta\gamma^3- \gamma^4)\]
   yielding,  when $\beta=-3$, 
   \[ p_w(z) =z (3 + z) (7 + 15 z + 14 z^2 + 6 z^3 + z^4) . \]
    We can solve $p_w(z)=-2$ to find $z=-1$,  $z=-2$ or
\[z=      \frac{1}{2} \left(-3\pm i \sqrt{-3+2 \sqrt{5}}\right),\hskip5pt z= \frac{1}{2} \left(-3\pm \sqrt{3+2 \sqrt{5}}\right) . \]
    We can also obtain disks around these points.  Note here that 
    \[ \gamma_0 = \frac{1}{2} \left(-3\pm i \sqrt{-3+2 \sqrt{5}}\right) \approx -1.5 + 0.606658 i \]
    actually gives the miminal covolume arithmetic hyperbolic lattice.  So,  while not pretty,  some of these polynomials are very important.  Also,  with a little work one can show that if $D_0 = |z-\gamma_0|<0.183$,  then $|p_w(z)+2|<0.557$ and so the disk $D_0$ contains no parameters for Kleinian groups other than $z_0$.  Of course once we have covered a region we are able to update previous calculations.  It is evident that a rather larger disk than the first one we found about $-2$ is excluded.  Using this,  we can improve the estimate $0.183$ that we just found.

    \bigskip
    
For each value of $p,q$ above we can argue as above using iteration with an
appropriate polynomial to obtain an open punctured disk about $0$
which does not contain the parameters for any Kleinian group and then build on the arguments as above. This gives an idea of the structure of the spaces ${\cal D}_{p,q}$ and from this we can extract the following information:  We begin with tables for the groups $\G_{3,i}$.  
 
 \begin{theorem}  Let $\G$ be a Kleinian group with $f\in \G$ elliptic of order $3$.  Let $\delta$ be half the distance from the axis of $f$ to a translate under $\G$. Then either $\delta>0.297$ or $\delta$ is one of the values in the table of groups $\G_{3,i}$ and $\G_{3,i}\subset\G$.
 \end{theorem}
 Thus the collaring radius for $f$ is at least $0.297$ or we can find a specific arithmetic subgroup in $\G$ generated by $f$ and a conjugate.
 \begin{table}
\centering
  \caption{Groups $\G_{3,i}$}
\begin{tabular}{|c|c|c|c|}
\hline
  $i$ & $\gamma$ value & minimal polynomial & $\delta(f,g)$ \\
\hline
\hline
1 & -1& $ z + 1 $& 0 \\
\hline
2 & -.3819 &$ z^2 + 3z + 1$ & 0\\
\hline
3 &-1.5 + .6066i &$z^4 + 6z^3 + 12z^2 + 9z + 1$ &.1970\\
\hline
4 &-.2118 + .4013i &$z^4 + 5z^3 + 7z^2 + 3z + 1 $&.2108\\
\hline
5 &-.5803 + .6062i &$z^3 + 4z^2 + 4z + 2$ &.2337\\
\hline
6 &-1.1225 + .7448i &$z^3 + 5z^2 + 8z + 5 $&.2448\\
\hline
7 &-.3376 + .5622i &$z^3 + 3z^2 + 2z + 1 $&.2480\\
\hline
8 &-.9236 + .8147i &$z^4 + 5z^3 + 8z^2 + 6z + 1$& .2740\\
\hline
9 &-1.5 + .8660i &$z^2 + 3z + 3$ &.2746\\
\hline
10 &-.7672 + .7925i&$ z^3 + 4z^2 + 5z + 3$ &.2770\\
\hline
11 &.0611 + .3882i &$z^4 + 5z^3 + 6z^2 + 1 $&.2788\\
\hline
12 &.2469 &$z^3 + 4z^2 + 3z - 1 $&.2831\\
\hline
13 &-.5284 + .7812i &$z^6 + 8z^5 + 24z^4 + 35z^3 + 28z^2 + 12z + 1 $&.2944\\
\hline
14 &-.1153 + .5897i &$z^3 + 3z^2 + z + 1 $&.2970\\
\hline
\end{tabular}
\end{table}
\bigskip

We also have the following more specific arithmetic data about these groups.
\begin{theorem} All the groups $\G_{3,i}$ are arithmetic with the exception of $\G_{3,12}$ which is $(2,3,7)$ triangle group,  and $\G_{3,13}$ a nearly arithmetic web group.
\end{theorem}~\begin{table}
\centering
  \caption{Groups $\G_{3,i}$}
\begin{tabular}{|c|c|c|c|c|c|c|c|}
\hline
  $i$ & discr. $k\G$  & Ram. & Covol. & $i$ & discr. $k\G$  & Ram. & Covol.\\
\hline
\hline
1 &   -- &-- & 0 S4 & 2 &  --&  -- & 0 A5\\  
\hline
3&  -275 & $\emptyset$ &  .0390& 4 &  -283 & $\emptyset$ & .0408\\
\hline
5 &  -44 & ${\cal P}_2$ & .0661 & 6 &  -23 & ${\cal P}_5$ & .0785\\
\hline
7 & -23 &${\cal P}_5$ & .0785 & 8 & -563& $\emptyset$ & .1274\\
\hline
9 & -3& $\emptyset$ &.0845 & 10 &  -31 &${\cal P}_3 $ & .0659\\
\hline
11 & -491& $\emptyset$ &.1028 & 12 &  --& --& $\infty$ Fuch.\\
\hline
13 & ? & ?  & $\infty$ Web. & 14 &   -76 & ${\cal P}_2 $&.1654\\
\hline
\end{tabular}
\end{table}
 \begin{table}
\centering
  \caption{Groups $\G_{4,i}$}
\begin{tabular}{|c|c|c|c|}
\hline
  $i$ & $\gamma$ value & minimal polynomial & $\delta(f,g)$ \\
\hline
\hline
1& -1 & $z + 1 $& 0 \\  \hline 
2 &-.5 + .8660i &$ z^2 + z + 1 $&.4157 \\  \hline 
3 &-.1225 + .7448i & $ z^3 + 2z^2 + z + 1 $&.4269 \\  \hline 
4 &-1 + i &$z^2 + 2z + 2 $&.4406 \\  \hline 
5 &-.6588 + 1.1615i &$z^3 + 3z^2 + 4z + 3 $&.5049 \\  \hline 
6 &.2327 + .7925i &$z^3 + z^2 + 1 $&.5225 \\  \hline 
7& -.2281 + 1.1151i & $z^3 + 2z^2 + 2z + 2 $&.5297 \\  \hline 
8 &.4196 + .6062i &$z^3 + z^2 - z + 1 $&.5297 \\  \hline 
9 &i &$z^2 + 1$& .5306 \\  \hline 
10 &.6180 &$z^2 + z - 1 $&.5306 \\  \hline 
11 &-1 + 1.2720i &$z^4 + 4z^3 + 7z^2 + 6z + 1 $&.5306 \\  \hline 
12 &-.4063 + 1.1961i &$z^4 + 3z^3 + 4z^2 + 4z + 1 $&.5345 \\  \hline 
13 &.7881 + .4013i &$z^4 + z^3 - 2z^2 + 1 $& .6130 \\  \hline 
\hline
\end{tabular}
\end{table}
\bigskip

   \begin{table}
\centering
  \caption{Groups $\G_{5,i}$}
\begin{tabular}{|c|c|c|c|}
\hline
  $i$ & $\gamma$ value & minimal polynomial & $\delta(f,g)$ \\
\hline
\hline
1&$ -.3819 $&$z - \beta - 1 $ & 0 \\  \hline
2 &$-.6909 + .7228i$&$ z^2 - \beta z + 1 $&.4568\\  \hline
3 &$.1180 + .6066i 4$&$z(z - \beta - 1)^2 - \beta - 1 $&.5306\\  \hline
4&$ -.1909 + .9815i 4$&$z^2 - (\beta + 1)z + 1 $&.6097\\  \hline
5 &$.61804$&$ z - \beta - 2$& .6268\\  \hline
6 &$.2527 + .8507i $&$z(z - \beta - 1)^2 + 1$& .6514\\  \hline
7 &$-.6909 + 1.2339i $&$z^2 - \beta z + 2 $&.6717\\  \hline
8&$ -.3819 + 1.2720i $&$z^3 - (2\beta + 1)z^2 + (\beta^2 + \beta + 2)z - 2\beta - 1 $&.6949\\  \hline
9&$ .1180 + 1.1696i $&$z^3 - (2\beta + 2)z^2 + (\beta^2 + 2\beta + 2)z - \beta $&.7195\\  \hline
10& $-.0817 + 1.2880i $&$z^4 - (2\beta + 1)z^3 + (\beta^2 + \beta + 2)z^2 - 2\beta z + 1 $&.7273\\  \hline
11&$ .6180 + .7861i$& $z^3 - (2\beta + 3)z^2 + (\beta^2 + 3\beta + 2)z + 1$& .7323\\  \hline
12& $ .8776 + .5825i $& $z(z - \beta)(z - \beta - 2)^2 + 1 $& .7725 \\  \hline 
\end{tabular}
\end{table}
 
 \newpage
 
 \subsection{Other commutator planes}
 
 Before moving on we would like to point out that many of the other spaces ${\cal D}_{p,q}$ are not as symmetric as the spaces ${\cal D}_{p,2}$.  Illustrated below is the space ${\cal D}_{3,4}$.  There are no obvious symmetries of this space other than complex conjugation.
 
  \scalebox{0.5}{\includegraphics[viewport= -30 600 400 790]{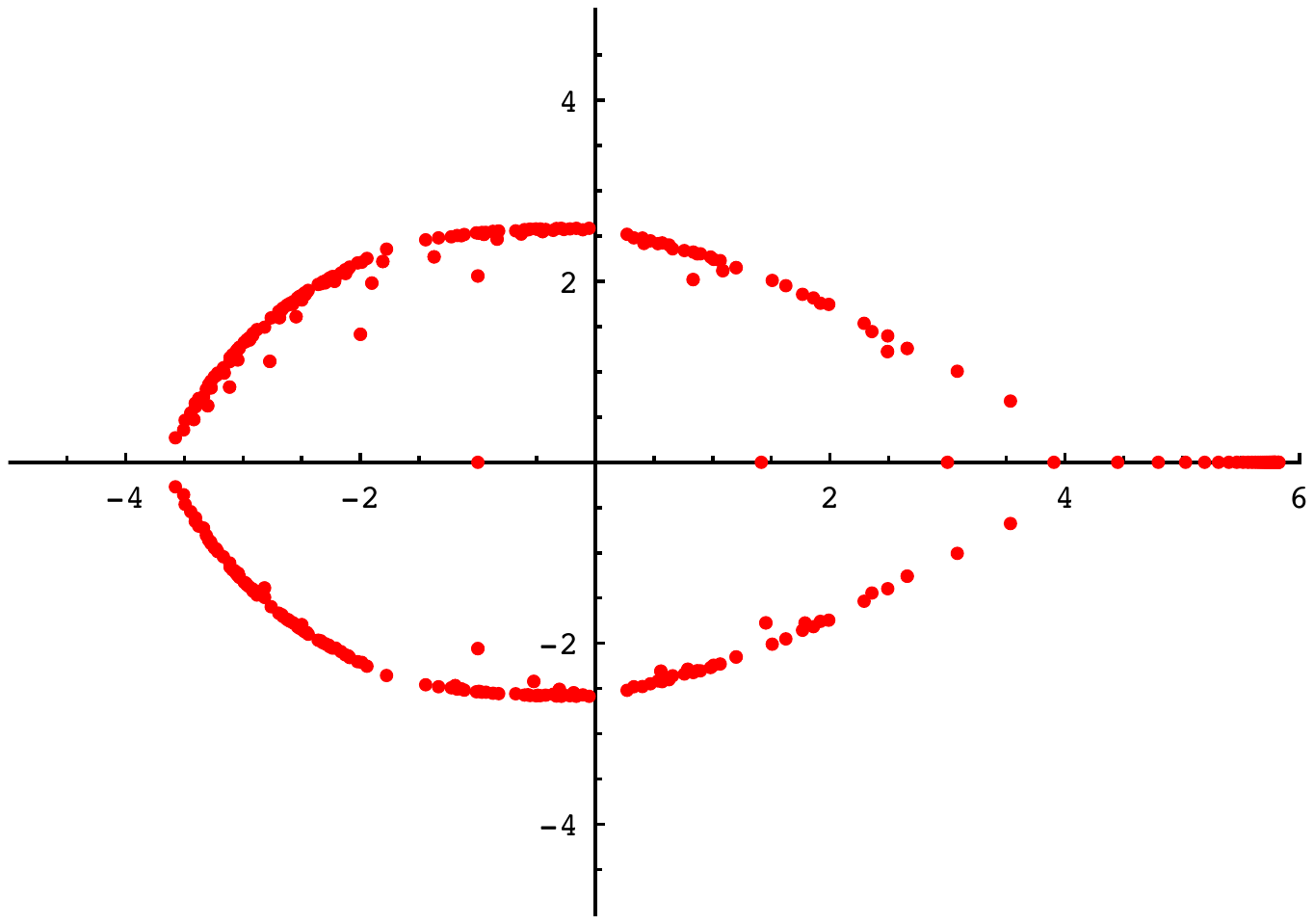}}
  
  \begin{center}  The $(3,4)$ commutator plane.
  \end{center}
 This picture is plotted data from $(3,0)$ and then $(4,0)$ orbifold Dehn surgery on two bridge links. From the matrix representation we calculate the $\gamma$ values.
 
  One can always find the extremal points on the real line by considering certain triangle groups.  For instance the $(3,4,\infty)$ triangle group has $\gamma = 3+2\sqrt{2}$,  while $[f,g]$ parabolic,  gives $\gamma=-4$.

\subsection{Axial distance spectra}
Here is the sort of data we pick up from descriptions of the spaces ${\cal D}_{p,q}$.  There are more general results available for other $p$ and $q$ \cite{GM1} and sharpness is usually given by an example of a truncated tetrahedral reflection group. \cite{CMar}.

\begin{theorem} Let  $\Gamma=\langle f,g\rangle$ be discrete generated by elliptics of orders  $p$ and  $q$ and let with $\delta+ i\theta$  be the complex distance between the axes.  Then   for the various values of $p$ and $q$ indicated in the tables below,  either $\delta(f,g)\geq \delta$,  the last tabulated value,  or $\delta(f,g)$ is one of the tabulated values and $\G$ is arithmetic if so stated.
\end{theorem}
 \begin{eqnarray*}
 \hline\\
 {\bf (p=3,q=2)} 		&&		    {\bf (p=4,q=2)} \\
0.19707 + i 0.78539   \; \mbox{  arithmetic}\;    &&	0.41572 + i 0.59803 \; \mbox{  arithmetic}\; \\
0.21084 + i 0.33189   \; \mbox{  arithmetic}\;  &&	0.42698 + i 0.44303 \; \mbox{  arithmetic}\;\\
0.23371 + i 0.49318   \; \mbox{  arithmetic}\; 	 &&0.44068 + i 0.78539 \; \mbox{  arithmetic}\;\\
0.24486 + i 0.67233   \; \mbox{  arithmetic}\;  &&	0.50495 + i 0.67478 \; \mbox{  arithmetic}\;\\
0.24809 + i 0.40575   \; \mbox{  arithmetic}\;  && 	0.52254 + i 0.34470 \; \mbox{  arithmetic}\;\\
0.27407 + i 0.61657   \; \mbox{  arithmetic}\; && 	0.52979 + i 0.24899 \; \mbox{  arithmetic}\;\\
0.27465 + i 0.78539   \; \mbox{  arithmetic}\;  &&	0.52979 + i 0.53640 \; \mbox{  arithmetic}\;\\
0.27702 + i 0.56753   \; \mbox{  arithmetic}\;   &&	0.53063 			 \; \mbox{  arithmetic}\;\\
0.27884 + i 0.22832   \; \mbox{  arithmetic}\; &&	0.53063 + i 0.45227 \; \mbox{  arithmetic}\;\\
  {\delta > 0.28088 }			 &&   {\delta   >  0.53264}\\ \\
  \hline
  \\
  {\bf (p=3,q=3)}					 &&   {\bf (p=2,q=5)} \\
0.39415 + i 1.57079   \; \mbox{  arithmetic}\; &&  	0.4568 	\; \mbox{  arithmetic}\;\\
0.42168 + i 0.66379   \; \mbox{  arithmetic}\;	 &&0.5306	\; \mbox{  arithmetic}\;\\
0.46742 + i 0.98637   \; \mbox{  arithmetic}\; &&	0.6097	\; \mbox{  arithmetic}\;\\
0.48973 + i 1.34468   \; \mbox{  arithmetic}\;    &&	0.6268	\; \mbox{  arithmetic}\;\\
0.49619 + i 0.81150   \; \mbox{  arithmetic}\; &&	0.6514  	\; \mbox{  non-arithmetic}\\
0.54814 + i 1.23135   \; \mbox{  arithmetic}\;    && 	0.6717  	\; \mbox{  non-arithmetic} \\
0.54930 + i 1.57079   \; \mbox{  arithmetic}\;      &&	0.6949	\; \mbox{  arithmetic}\;\\
0.55404 + i 1.13507   \; \mbox{  arithmetic}\; &&	0.7195	\; \mbox{  arithmetic}\;\\
0.55769 + i 0.45665   \; \mbox{  arithmetic}\; &&	0.7273	\; \mbox{  arithmetic}\;\\
  {\delta > 0.56177 } 			 &&   {\delta  > 0.73}\\
  \end{eqnarray*}
 
 \subsection{Distances between spherical points}
  If there is a non-simple elliptic,  or equivalently when there is a
finite spherical subgroup, in our Kleinian group $\G$  the arguments trying to find covolume estimates are of a different nature
but still based around the
knowledge of the spectra of possible axial distances for elliptics of
orders $3,4$ and $5$.  A {\em spherical point} is a point stabilized
by a spherical triangle
subgroup of a Kleinian group - namely the tetrahedral, octahedral and icosahedral groups $A_4, S_4$ and $A_5$.  Geometric
position arguments based around the axes emanating from a spherical
point show the distances
between spherical points to be uniformly bounded below.  In
\cite{GMannals} we identify the initial part of this spectra of distances,
significantly extending
earlier work of Derevnin and Mednykh, \cite{DM}.  Again,  a crucial
point is establishing arithmeticity of the first few extrema.  This
enables us to use
arithmetic criteria to eliminate small configurations from
consideration. 

 Once this is done a ball of maximal radius $r$ about a
suitable spherical point will
be precisely invariant (and of volume $\pi(\sinh(2r)-2r)$) and provide volume bounds using the following obvious result.

\begin{theorem}  Let $\G$ be a Kleinian group and $x_0\in \IH^3$ with stabiliser $\G_{x_0}$.  Let
\[ r = \min \left\{\; \frac{1}{2} \rho_{\IH}(g(x_0),x_0):g\in \G\setminus\G_{x_0}\;\right\} . \]
Then
\begin{equation}
{\rm vol}_{\IH}(\IH^3/\G) > \pi (\sinh(2r)-2r)/|\G_{x_0}|
. \end{equation}
\end{theorem}
 And as an example,  we will see in a minute that the distance between two points fixed by the spherical triangle group $A_5$ is either one of $6$ listed values and the group is arithmetic,  or it exceeds $2.826$.  Thus we obtain
 \begin{theorem}   Let $\G$ be a Kleinian group with a subgroup isomorphic to $A_5$.  Then either $\G$ is arithmetic (and one of $6$ possibilities listed) or
 \begin{equation}
 {\rm vol}_{\IH}(\IH^3/\G)> \pi (\sinh(2.826)-2.826)/60 = 0.292\ldots
  . \end{equation}
 \end{theorem}
 This estimate is nearly an order of magnitude more than we need and can in fact be improved by using B\"or\"oczky's sphere packing arguments  \cite{Boroc}.   The closest two $A_5$ spherical points can be is about $1.38257\ldots$,  however this distance does not provide enough volume. This is how the use of arithmeticity allows us to push past extrema.  It is important to notice that we don't need to know what the minimal volume arithmetic group is since we have identified the exceptional candidates and we can use Borel's formula to identify the volume in each case.
 
 \medskip
 
 The way we find the data on the distance between spherical points is as follows.  The axes of elliptics emanating from a spherical point are in general position.  It follows that the distance between these points is also bounded below and it is an exercise in hyperbolic trigonometry to obtain these bounds.  This is illustrated below. 
 
 \medskip
 
  \rotatebox{-89.99}{\scalebox{0.8}{\includegraphics[viewport= 230 250 380 270]{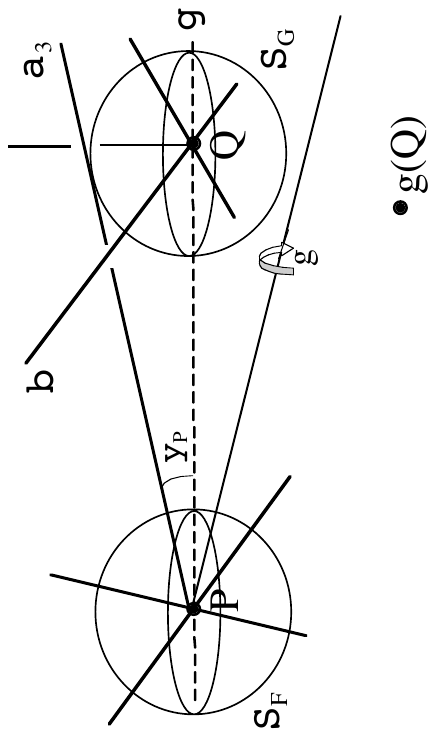}}}
  
  \medskip
  
  It becomes clear that in order to get points as close as possible the spherical points $P$ and $Q$ should lie on a common axis.  Geometrically this means that $P$ and $Q$ both lie on some
common elliptic axis of $\G_P$ and $\G_Q$.  Then as we push them closer another pair of axes should meet (since the distance is at least $\delta_{p,q}$ or $0$).  This motivates our study of $(p,q,r)$ Kleinian groups in the next section as these elliptic axes now form a hyperbolic triangle.    The corresponding
elliptics generate a $(p,q,r)$--Kleinian group.  Hence the distance $\rho(P,Q)$ is bounded below by the
corresponding possible edge-lengths of the associated hyperbolic triangle which we compute.  These $(p,q,r)$ groups are also needed in the study of the Margulis constant.

\medskip

We record the various possibilities for distances between spherical points in the following theorem.  In all instances, with the notable exception  of $S_4$, we shall see that these distances are the sharp bounds.  

\begin{theorem} Suppose that $P$ and $Q$ are spherical triangle points lying on a common axis 
$\eta$ of order $n$ and suppose that an axis of $P$ meets an axis of $Q$, possibly on the sphere at
infinity, other than the axis $\eta$.  If $P$ is isomorphic to $A_4$ and if $Q$ is isomorphic to $A_4$, $S_4$ or $A_5$, then
the distance $\rho(P,Q)$ is either one of the $5$ entries in the appropriate table below and
$\eta$ has the corresponding order $n$ or $\rho(P,Q)$ exceeds the $4^{th}$ tabulated value or the  $5^{th}$ tabulated value in case the orders of the  $4^{th}$ and $5^{th}$ entry are the same. 
{\rm
\begin{center}
\begin{tabular}{|c|c|}
\multicolumn{2}{c}{ {\boldmath $(A_4,A_4)$}} \vspace{.02in}\\
\hline
$n$ & $\rho(P,Q)$ \\
\hline
3 & 0.69314... \\
\hline 
3 & 0.76914... \\
\hline 
3  & 0.92905... \\
\hline 
3 & 1.0050... \\
\hline 
2  & 1.06128...\\ 
\hline
\end{tabular}
\hspace{.1in} 
\begin{tabular}{|c|c|}
\multicolumn{2}{c}{  {\boldmath $(A_4,S_4)$}} \vspace{.02in}\\
\hline
 $n$ & $\rho(P,Q)$ \\
\hline
3 & 1.01481... \\
\hline 
3 & 1.31696... \\
\hline 
3  & 1.43364... \\
\hline 
2 & 1.43796... \\
\hline
3 & 1.49279... \\
\hline 
\end{tabular}
\hspace{.1in}
\begin{tabular}{|c|c|}
\multicolumn{2}{c} { {\boldmath $(A_4,A_5)$}}\vspace{.02in} \\
\hline
$n$ & $\rho(P,Q)$ \\
\hline
3 & 1.22646... \\
\hline 
3 & 1.62669...\\
\hline 
2  & 1.76110...  \\
\hline 
3 & 1.87988...  \\
\hline
3 & 1.98339...  \\
\hline
\end{tabular}
\bigskip
\end{center}}

\noindent If $P$ is isomorphic to $S_4$ or $A_5$ and if $Q$ is isomorphic to $S_4$ or $A_5$, then the distance  
$\rho(P,Q)$ is either one of the first $7$ entries in the appropriate table below and $\eta$ has the
corresponding order $n$ or $\rho(P,Q)$ exceeds the $6^{th}$ tabulated value or the  $7^{th}$ tabulated value in case the orders of the  $6^{th}$ and $7^{th}$ entry are the same. 
\end{theorem}

\begin{center}
\begin{tabular}{|c|c|}
\multicolumn{2}{c}{ {\boldmath $(S_4,S_4)$}} \vspace{.02in}\\
\hline
 $n$ & $\rho(P,Q)$ \\
\hline
4 & 1.06128...  \\
\hline
4 & 1.12838... \\
\hline 
3,4 & 1.31696... \\
\hline
4 & 1.38433...  \\
\hline
4 & 1.48710...  \\
\hline 
3  & 1.56680... \\
\hline 
2 & 1.70004... \\
\hline 
\end{tabular}
\hspace{.1in}
\begin{tabular}{|c|c|}
\multicolumn{2}{c}{  {\boldmath $(S_4,A_5)$}} \vspace{.02in}\\
\hline
$n$ & $\rho(P,Q)$ \\
\hline
3 & 1.22646... \\
\hline 
3 & 1.98339... \\
\hline 
3 & 2.13275...  \\
\hline 
2 & 2.27311...  \\
\hline
3 & 2.34868...  \\
\hline 
2 & 2.35576...  \\
\hline
2 & 2.83641... \\
\hline
\end{tabular}
\hspace{.1in}
\begin{tabular}{|c|c|}
\multicolumn{2}{c}{  {\boldmath $(A_5,A_5)$}} \vspace{.02in}\\
\hline
 $n$ & $\rho(P,Q)$ \\
\hline
5 & 1.38257...  \\
\hline 
5 & 1.61692...  \\
\hline 
3 & 1.90285...  \\
\hline 
5 & 2.04442...  \\
\hline
5 & 2.16787...  \\
\hline
5 & 2.22404...  \\
\hline
2 & 2.82643... \\
\hline 
\end{tabular}
\end{center}

\bigskip

\noindent{\bf Remarks.} The last entry in each table is either the next possible distance or 
the smallest distance possible on a common axis with the given order.  
Most of the distances are achieved in the orientation-preserving subgroups of groups generated by
reflection in the faces of a hyperbolic tetrahedron.   

\bigskip

 There are of course additional complications in the calculations and sometimes we have to look elsewhere to pick up a little bit more volume,  but in the end we obtain the main
result of \cite{GMannals} using the volume estimates as above.

\begin{theorem}  Let $\G$ be a Kleinian group containing a finite spherical triangle group,  either tetrahedral, octahedral or icosahedral.  Then
\begin{equation}
{\rm vol}_{\IH}(\IH^3/\G) \geq 0.03905\ldots
\end{equation}
and this value is sharp and obtained in the minimal covolume arithmetic lattice $\G_0$ described above.
\end{theorem}

\bigskip

The above discussion shows us that there are one remaining case  to
deal with in order to identify the two smallest covolume Kleinian
groups $\Gamma$. That case is when   all the elliptics in $\Gamma$ are order $2$.
We can further assume that  two elements of order $2$ can only meet at
right angles and then the group they generate will be a Klein
$4$-group for otherwise we would pick up an element of order $3$ or more. Already in \cite{GM8} it
was shown that these Klein $4$--groups appear in many extremal
situations and one of the first results we establish in \cite{MMannals} is the following
universal constraint concerning discrete groups generated by two
loxodromics whose axes meet orthogonally.   This result significantly
refines a particular case of an earlier theorem \cite{GM3} concerning
groups generated by loxodromics with intersecting axes.

\begin{theorem}\label{orthogonal}
Let  $f$ and $g$ be loxodromic transformations generating a discrete
group such that the axes of $f$ and $g$ meet at right
angles.  Let
$\tau_f$ and $\tau_g$ be the respective translation lengths of $f$
and $g$.  Then
\begin{equation} \label{orthoineq1}
\max \{ \tau_f, \tau_g \} \geq \lambda_{\perp}
\end{equation}
where
\begin{equation}
   \lambda_{\perp} = \arccosh\left(\frac{\sqrt{3}+1}{2}\right)= 0.831446\ldots
. \end{equation}
Equality holds when $\langle f,g \rangle$ is the two-generator
arithmetic torsion-free lattice with presentation
\[ \Gamma = \langle f,g : fg^{-1}fgfgf^{-1}gfg =
gfg^2fgf^{-1}g^{2}f^{-1} = 1 \rangle . \]
This group is a four fold cover of (4,0) \& (2,0) Dehn surgery on the
2 bridge link complement $6^{2}_{2}$ of Rolfsen's tables, 
\cite{Rolfsen}. $\IH^3/\Gamma_0$  has volume
$1.01494160\ldots$,  Chern--Simons invariant $0$ and homology $\IZ_3+\IZ_6$.
\end{theorem}

   (We expect that Theorem \ref{orthogonal} represents the extreme case
independently of the angle at which the axes of loxodromics meet.)

\medskip

The proof of this theorem goes through a search of the two-dimensional parameter space using the good words (in fact about 22 of them) to eliminate small configurations.  

\medskip

   Groups with orthogonal loxodromics appear naturally as they ``wind
up'' Klein $4$--groups when projected to
the quotient orbifold.    
 It is clear how to use Theorem \ref{orthogonal} to get volume estimates as it implies a collar about one or the other loxodromics.  This is then the final piece of the puzzle in solving Siegel's problem.
     
\section{$(p,q,r)$--Kleinian groups}

 Here  we consider Kleinian groups generated by three rotations of orders $p,q$ and $r$ whose axes of rotation
are coplanar.  That is they lie in a common hyperbolic plane in $\IH^3$.  We call such a group a $(p,q,r)$--Kleinian group when the
generating elements have orders $p$, $q$ and $r$.  These groups are natural generalisations of triangle groups
and also of the index 2 orientation-preserving subgroups of the group generated by reflection in the faces of a
hyperbolic tetrahedron.  In particular the $(p,q,r)$--Kleinian groups properly contain both these classes. 
These groups also appear as extremals for  interesting problems in the theory of Kleinian groups.  Apart
from the triangle groups and reflection subgroups there are a number of interesting infinite families of web
groups.  These are the groups whose limit set is a Sierpinki gasket,  that is
$\oC$ minus an infinite family of disjoint disks so that what remains has no interior.  The picture below is from Curt McMullen's homepage.\\

\medskip

\scalebox{0.35}{\includegraphics[viewport= -150 400 580 790]{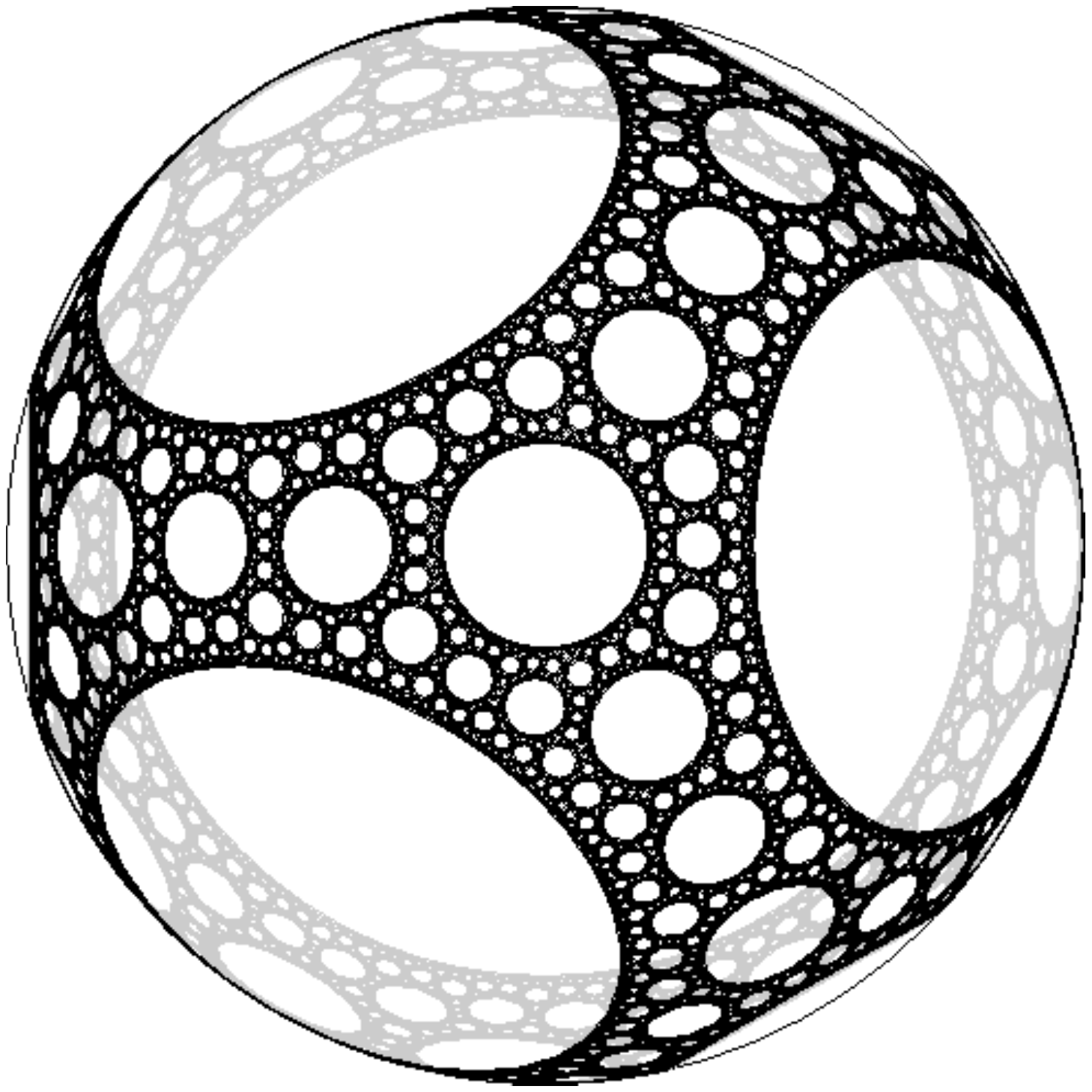}}

\bigskip

\bigskip

\begin{center} Sierpinsky limit set
\end{center} 
The stabilizers of
components of the ordinary set, the union of the disks removed, are triangle groups.  

When the axes of the rotations bound a hyperbolic triangle and the group is discrete, these groups will have a presentation
of the form
\begin{equation}
\langle a,b,c : a^p=b^q=c^r=(ab)^\ell=(ac)^m=(bc)^n \rangle
. \end{equation}
When such a group is discrete, the sum $\Sigma=\frac{1}{\ell}+\frac{1}{m}+\frac{1}{n}$ determines whether it is compact
($\Sigma>1$),  non-compact but of finite covolume ($\Sigma=1$), or a web group $(\Sigma<1)$.  If $p=q=r=2$ we in fact have a
$\IZ_2$--extension of a triangle group.

Recall that for  a Kleinian group  $\G=\langle g_1,g_2,g_3,\ldots \rangle $ we define the {\em Margulis
constant} for $\G$ to be
\begin{equation}
{\bf m}_\G = \inf_{x\in \IH} \sup_{i\geq 1} \{ \rho(g_i(x),x) \}
. \end{equation}

The Margulis constant for a particular group depends on the particular choice of generators. The Margulis constant for
Kleinian groups is
\begin{equation}
c_M= \inf_{\G} m_{\G}
\end{equation}
where the infimum is over all Kleinian groups $\G$ and all choices of generators. 

We point out the
the following obvious lemma
\begin{lemma}  If $\G=\langle g_1,g_2,g_3,\ldots \rangle $ and if
$\G^\prime=\langle h_1,h_2,h_3,\ldots \rangle $ with $h_i \in \{
g_1,g_2,g_3,\ldots \} $ for all
$i\geq 1$,  then
\begin{equation}
{\bf m}_{\G^\prime} \leq {\bf m}_\G
. \end{equation}
\end{lemma}

Thus in order to bound the Margulis constant for a given Kleinian group, we may choose a suitable subcollection of generators
and estimate the Margulis constant for the Kleinian group they generate. In
particular,  if there are two or three generators generating a nonelementary
subgroup,  then the Margulis constant is bounded below by the constant for such groups.  An elementary analysis reduces to the case that there are only two or three generators, \cite{Inada}.

\begin{theorem} Let $\G=\langle g_1,g_2,g_3,\ldots \rangle$ be a Kleinian group.  Then
there are indices $i_1<i_2$ such that $\G^\prime=\langle g_{i_1},g_{i_2} \rangle$ is a Kleinian group or there are indices
$i_1<i_2<i_3$ such that $\G^\prime=\langle g_{i_1},g_{i_2},g_{i_3} \rangle$ is a Kleinian group. 

The case of three generators is necessary only if $g_{i_1}$, $g_{i_2}$, $g_{i_3}$ are elliptic
elements of finite order $2\leq p\leq q \leq r \leq 6$ and $(p,q,r)$ is one of the following triples:
\medskip

$(2,2,k)$, $2 \leq k \leq 6$, \hskip15pt $(2,3,k)$, $3 \leq k \leq 6$, \hskip15pt $(2,k,k)$, $4 \leq k \leq 6$,
\medskip

$(3,3,k)$, $3 \leq k \leq 6$, \hskip15pt $(3,k,k)$, $4 \leq k \leq 6$, \hskip15pt $(k,k,k)$, $4 \leq p \leq 6$.  
\medskip

In all but the case $(2,2,k)$ the elliptic axes are coplanar and bound a hyperbolic triangle.
\end{theorem}
\noindent{\bf Proof:}  
If three generators are in fact necessary, then each pair of them must generate an elementary group. It is clear from the
classification of the elementary discrete groups, that no $g_{i_j}$ can be parabolic, loxodromic or elliptic of order at
least $7$.  If two elliptic elements generate an elementary group,  their axes must meet (possibly on the Riemann sphere) unless they
both have order $2$.  Suppose that $q \geq 2$.  Then the axes of the generators must pairwise intersect, again possibly on the
Riemann sphere, and hence they are coplanar and bound a hyperbolic triangle.  Each vertex of this triangle is stabilised by an
elementary group.  Using the classification of the elementary discrete groups once more, we obtain the stated restrictions on
$(p,q,r)$. \hfill $\Box$

\bigskip

The labeling $(p,q,r)$ and the hyperbolic triangle obtained when the axes are coplanar determine the Kleinian group
completely.  We thus want metric restrictions on this hyperbolic triangle.  However, first we show that the groups we are
interested in embed in reflection groups.

\subsection{The associated reflection group}

The following result helps in deciding when a possible triple $(p,q,r)$ and the associated
hyperbolic triangle arise from a discrete group.

\begin{theorem}  Let $f,g,h$ be elliptic M\"obius transformations whose axes are coplanar.  Then there are reflections
$\s_i$ in hyperplanes $\Pi_i$, $i=0,1,2,3$, such that
\begin{enumerate}
\item $\langle f,g,h \rangle$ is a subgroup of index 2 in $\langle \s_0,\s_1,\s_2,\s_3 \rangle$ and hence both groups are
simultaneously discrete or non--discrete.
\item Let $\theta_{i,j}$ denote the dihedral angle between $\Pi_i$ and $\Pi_j$,  $0\leq i \neq j \leq 3$ and let
$\theta_{i,j}= 0$ if they do not meet.  Then  $\theta_{0,1}=\pi/p$, $\theta_{0,2}=\pi/q$ and  $\theta_{0,3}=\pi/r$.  If
\[ \theta_{1,2} + \theta_{1,3} + \theta_{2,3} \geq \pi,
\] then the hyperplanes $\Pi_i$ bound a hyperbolic tetrahedron.
\end{enumerate}
\end{theorem}
{\bf Proof:}  We may assume that $f,g,h$ are primitive elliptic elements of order $p,q,r$ respectively.  Let $\Pi_0$ be the
plane in which the axes of
$f,g,h$ lie and let
$\sigma_0$ be reflection in $\Pi_0$.  $\Pi_0$ separates $\IH^3$ into two components.  Fix one such component $U$. Let
$\Pi_1$ be the plane having dihedral angle $\pi/p$ with $\Pi_0$, chosen so that the acute angle is formed in $U$.  Do the same
for $g$ and $h$.  Then $\s_0\s_1=f$ or $f^{-1}$; in the latter case we replace $f$ with $f^{-1}$.  We again do the same for 
$\s_0\s_2=g$ and $\s_0\s_3=h$.  Next we see that $\s_1\s_2=\s_1\s_0 \s_0\s_2 = f g^{-1}$. Arguing similarly we conclude
that $\langle f,g,h \rangle$ has index 2 in the reflection group and thus the first part of
the theorem is proved.  The second part is a well known fact from hyperbolic geometry, 
\cite{Rat}.

\bigskip

The Poincar\'e theorem,  carefully presented in \cite{Maskit},  can be used to obtain a presentation for the reflection group and then for the subgroup
$\langle f,g,h \rangle$ when it is discrete.  In fact the group is discrete if and only if the dihedral angles of
intersection of the hyperplanes are submultiples of $\pi$ together with a  number of exceptions, in fact infinitely many, but
all of a similar pattern.  

\subsection{Admissible Triangles and Discreteness}

We say that a hyperbolic triangle in $\IH^3$ is an {\em admissible (p,q,r)-triangle}
if the triangle
is formed from the axes of elliptic elements  $f,g,h$ of order $p,q,r$
respectively and that $\langle f,g,h\rangle$ is a discrete group.

\bigskip

We want to see what sort of admissible hyperbolic triangles our axes might bound.  There are two obvious conditions. The first
is that the sum of the (interior) angles of intersection of the axes is smaller than $\pi$ for we have a hyperbolic triangle of
positive area.  The second is that all the interior angles must be chosen from those of our list or in the subsequent
remark.

These two conditions readily reduce our list of possibilities to manageable proportions.  However, we seek to identify all
such groups, and so we must decide when a candidate generates a discrete group.

The third condition we apply concerns the sum of the dihedral angles at each vertex (opposite the angle of
intersection) formed in the associated reflection group.  Let $\Sigma$ be the sum of dihedral angles opposite each vertex. 
If $\Sigma \geq \pi$,  then the three hyperplanes (those other than the hyperplane containing our triangle) of the associated
reflection group meet at a finite point and together with the hyperplane containing our triangle bound  a finite volume
hyperbolic tetrahedron.  Our
$(p,q,r)$--group is of index 2 (if either group is discrete).  Notice that there are more reflection groups than those
appearing in the standard lists since we do not place the restriction that the dihedral angles are submultiples of $\pi$.  We
have tabulated the reflection groups appearing below. 

If $\Sigma <
\pi$,  then the three planes subtend a hyperbolic triangle on the sphere at infinity,  and therefore generate a reflection
group in the sides of a hyperbolic triangle.  This group must be discrete and this places restrictions on the angles.    Of
course if all angles are submultiples of $\pi$,  then the group is certainly discrete. However it is shown by Knapp
\cite{Kn} that in a discrete reflection group in the sides of a hyperbolic triangle,  if one angle is of the form
$2\pi/p$,
$p$ odd,  then the other two angles are of the form
$\pi/n,\pi/n$ or 
$\pi/2,\pi/p$ (with a further exception for the $(2,3,7)$--triangle group which will not concern us).  If the angle is of the
form $3k\pi/p$,  then the other two angles are $\pi/3,\pi/p$.  And finally if an angles is of the form $4\pi/p$,  then the
other two angles are $\pi/p,\pi/p$.  No other rational multiples of $\pi$ occur.  Knapp's theorem is more or less restated in
when we discuss admissible $(2,2,2)$--triangles.  This theorem gives us a powerful tool to limit the possible
dihedral angles.

\bigskip

There are two basic cases to consider.  The first case is where all the dihedral angles are submultiples of $\pi$ and the
second is where one or more such angles is not.  In the first case,  the group is discrete by the Poincar\'e theorem.  Thus we obtain necessary and sufficient conditions for discreteness.  It is a fairly routine matter to run
through the possible angles and find all the discrete groups.  We have indicated the Margulis
constant for these groups and also the associated tetrahedral reflection group when the group has finite covolume.  From the
tables of covolumes of these reflection groups,  one can obtain the covolume of our $(p,q,r)$--Kleinian group by doubling
this value.  We have also indicated the structure of the stabilizers of vertices of our hyperbolic triangle.  These are
important for our applications above in finding the distances between spherical points.   When the group is a web group,  we further indicate
the triangle group which is a component stabilizer.

From the tables we deduce the following theorem.

\begin{theorem}  Let $\Gamma$ be a $(p,q,r)$--Kleinian group,  with $3\leq p\leq q\leq r$.  Then $\Gamma$ is the index-two
orientation-preserving subgroup of the group of reflections in a hyperbolic tetrahedron.
\end{theorem}

Note that there are nonclassical reflection groups here.  By this we mean the angles at edges of the hyperbolic
tetrahedron may not be submultiples of $\pi$.

\subsection{Dihedral angles not submultiples of $\pi$.}

In this case we proceed as above to identify the triangles involved.  We  use Knapp's theorem to eliminate a number of
cases as well as the simple observation that if the dihedral sum $\Sigma\geq \pi$, then the stabiliser of this vertex must
be one of the elementary groups. Thus for instance there can be no angles of the form $m\pi/4$ and $n\pi/5$ meeting. 
If we have a candidate satisfying all these conditions we must still decide if it is discrete.  We do this in the following
way.  Let $\G$ be a $(p,q,r)$--Kleinian group and the dihedral angle between two hyperplanes $\Pi_1$ and $\Pi_2$ at some
vertex,  say $p-q$ with angle $\theta_{pq}$, is of the form
$k\pi/m$ (usually
$k=2$).  Let $\s_i$ denote the reflection in the hyperplane $\Pi_i$, $i=1,2$.  Then in $\langle \s_1,\s_2 \rangle$ we find
reflections in hyperplanes meeting $\Pi_1$ and $\Pi_2$ at all the diheral angles $m\pi/n$,  $m=1,2,\ldots,n$. At least one of
these hyperplanes intersects the interior of our triangle.  This intersection is the axis of an elliptic of some order, say
$s$, and thus we have found two new groups of the form $(p,r,s)$ and $(q,r,s)$ each of which must be discrete.  We will not
have $s=2$ unless the tetrahedron is symmetric. Notice that this process provides $k$ such groups,  all having angle of
intersection at the
$p-q$ vertex of the form
$\theta_{pq}/k$.  Notice the hyperbolic area of the new triangle has decreased and so this process of subdivision must stop
if the group is discrete and at one of the groups whose dihedral angles are all submultiples of
$\pi$.  In practise we find the subdivision process stops at the first stage and the group is either discrete as an index-two
subgroup of a group with dihedral angles all submultiples of
$\pi$,  or the angles of intersection of the axes $s$ and $r$ are not permissible in a discrete group.  Sometimes we must
subdivide 3 times.  This subdivision exhibits some interesting inclusions among the $(p,q,r)$--groups.

Incidentally this process identifies for us all discrete reflection groups in the faces of a hyperbolic tetrahedron such
that at one face all the dihedral angles are submultiples of $\pi$.   

The resulting discrete $(p,q,r)$--Kleinian groups are tabulated below.  
 
 \subsection{Admissible $(2,2,2)$-triangles}

A $(2,2,2)$--Kleinian group when restricted to the hyperbolic plane containing the axes is actually a reflection group in the
sides of a hyperbolic triangle,  and thus is an index-two supergroup of a Fuchsian triangle group.  The results of Knapp
\cite{Kn} can be used to determine all possible triples of angles.  Each such angle is a rational multiple of $\pi$,  and in
fact with few exceptions will be a submultiple of $\pi$.

\begin{theorem}  Let $\Gamma$ be a  $(2,2,2)$--Kleinian group. Then the angles of intersection of the axes are of the
following form.
\begin{enumerate}
\item $\pi/\ell$, $\pi/m$, $\pi/n$ with $\frac{1}{\ell}+\frac{1}{m}+\frac{1}{n} < 1$.
\item $2\pi/\ell$, $\pi/m$, $\pi/m$ with $\frac{1}{\ell}+\frac{1}{m} < \frac{1}{2}$.
\item $2\pi/\ell$, $\pi/2$, $\pi/\ell$ with $\ell \geq 7$.
\item $3\pi/\ell$, $\pi/3$, $\pi/\ell$ with $\ell \geq 7$.
\item $4\pi/\ell$, $\pi/\ell$, $\pi/\ell$ with $\ell \geq 7$.
\item $2\pi/7$, $\pi/3$, $\pi/7$.
\end{enumerate}
\end{theorem}

The reader may easily verify the geometrically evident fact that in a group generated
by three coplanar axes of order $2$ the Margulis constant is twice the radius of the
maximal inscribed disk of the triangle.  The following formula for the radius $r$ of the
inscribed disk in a hyperbolic triangle with angles $\alpha,\beta,\gamma$ is found in
\cite{Beardon}.  
\begin{equation}
\tanh^2(r)=
\frac{\cos^2\alpha+\cos^2\beta+\cos^2\gamma+2\cos\alpha\cos\beta\cos\gamma-1}{2(1+\cos\alpha)(1+\cos\beta)(1+\cos\gamma)}
. \end{equation}
 
\begin{theorem} Let $\G$ be a $(2,2,2)$--Kleinian group. Then
\begin{equation}
{\bf m}_\G \geq 
2\tanh^{-1}\sqrt{\frac{\cos^2\pi/3+\cos^2\pi/7-1}{2(1+\cos\pi/3)(1+\cos\pi/7)}}=0.2088\ldots
. \end{equation} 
The bound is sharp and uniquely achieved by the reflection group of the (2,3,7)
triangle.
\end{theorem} 

\subsection{The Margulis constant for groups generated by an admissible triangle}

In this section we indicate how the Margulis constant for each group in the  
tables is calculated.

We say that a group $\G$ is generated by an admissible triangle if $\G =
\langle f_1,f_2, f_3
\rangle$,  $f_i$ are all elliptic of order less than $6$ and the axes $\alpha_i$ of
$f_i$ form an admissible triangle.   The Margulis
constant for such a
group $\G=\langle f_1,f_2,f_3 \rangle$ is simply
\[ m_\G = \inf \{ \rho(x,f_i(x)) : x\in \IH,  i=1,2,3 \} \]
where $\rho$ is the hyperbolic metric of $\IH^3$.  A
moments thought should convince the reader that the infimum is actually a
minimum and that (by
convexity) the minimum is attained at a point  $x_0$ lying in the plane
spanned by the three
axes interior to the region bounded by the triangle and such that
\begin{equation}
\rho(x_0,f_1(x_0)) = \rho(x_0,f_2(x_0)) = \rho(x_0,f_3(x_0))
. \end{equation}
Since the plane spanned by the three axes is totally geodesic  we  may use
the metric
geometry in that plane.    Thus 
\begin{equation}
m_\G = \rho(x_0,f_i(x_0)) = 2 {\rm arcsinh}\left( \sin(\pi/n_i) \sinh(d_i)
\right)
\end{equation} 
where $f_i$ has order $n_i$ and the distance to the ``Margulis point'' is $d_i$.
Now after quite a bit of hyperbolic trigonometry we see
$\sinh^2\left(\frac{m_\G}{2}\right)$ is
equal to
\begin{equation}
\frac{\sin^2(\pi/n_2)(\cosh^2(\ell_2)-1)}
{1-\cosh^2(\ell_2)+\cot^2(\theta_{1,2})+\cot^2(\theta_{3,2})
+2\cot(\theta_{1,2})\cot(\theta_{3,2})\cosh(\ell_2)}
\end{equation}
where $\ell_i$ is the edge length of the triangle coinciding with the axis of $f_i$.  This gives $m_\G$ in terms of the three angles in view of our formulas above.

It is clear that the above formula remains valid if $\theta_2 = 0$,  since
$\cot(\theta_{1,2})$ and $\cot(\theta_{1,3})$ are not defined in terms of
$\theta_2$ and
$\cosh(\ell_2)$ is finite.   While if all the vertices are infinite we find,  after some simplification,
the following
symmetric formula
\[
\frac{1}{\sinh^2\left(\frac{m_\G}{2}\right)}
=
 \frac{1}{\sin(\pi/n_1)\sin(\pi/n_2)}+\frac{1}{\sin(\pi/n_1)\sin(\pi
/n_3)}
+\frac{1}{\sin(\pi/n_2)\sin(\pi/n_3)}
. \]

\subsection{Margulis constant for the $(2,2,n)$ case}

We can also deal with the case that there are two elliptic axes of
order $2$ meeting a
third axis of order $n\geq 3$ not only in the case of a $(2,2,n)$--Kleinian group,  but also in the case that the three axes
may not be coplanar.  In this latter case we are
unfortunately unable to give explicit computations of the constant. Instead
we will just give
a lower bound, which is sharp in some instances. 
{\scriptsize
\begin{center}
\begin{tabular}{|c|c|c|c|c|c|}
\multicolumn{6}{c}{\bf Margulis Constants for
(2,2,n)--nontriangles}\\
\hline
Order& Vertex & $\rho$ & $\theta$ & $\psi$ & Constant \\
\hline
$3$ & $A_4$ & $0.6931$ & $0.95531$ & $\pi/6$ & $.230989$ \\
\hline
$3$ & $A_5$ & $1.736$ & $0.36486$ & $\pi/3$ & $.276814$ \\
\hline
$3$ & $S_4$ & $1.059$ & $0.61547$ & $.8496$ & $.23167$ \\
\hline
$3$ & $A_5$ & $1.736$ & $1.20593$ & $0$ & $.302572$ \\
\hline
$4$ & $S_4$ & $1.059$ & $\pi/4$ & $.62894$ & $.28693$ \\
\hline
$5$ & $A_5$ & $1.382$ & $.55357$ & $.85474$ & $.26088$ \\
\hline
$5$ & $A_5$ & $1.382$ & $1.01772$ & $0$ & $420531$ \\
\hline
$6$ & $\Delta_{(2,3,6)}$ & $\infty$ & $0$ & $0$ & $.594241$ \\
\hline
\end{tabular}
\end{center}}

\subsection{Three elliptics of order 2}

In this section we show that a Kleinian group generated by two loxodromic
transformations can be identified as a subgroup of index at most two in a Kleinian
group generated by three elliptics of order 2 and that this supergroup has a smaller
Margulis constant.

\begin{theorem}  Let $\G=\langle f,g \rangle$ be a Kleinian group generated by two
loxodromic transformations.  Then there are three elliptic transformations
$\phi_1,\phi_2,\phi_3$,  each of order two,  generating a Kleinian
group $\G^\prime =
\langle
\phi_1,\phi_2,\phi_3 \rangle$ which contains $\G$ with index at most two and 
\begin{equation}
m_{\G}^{\prime} < m_\G
. \end{equation}
\end{theorem}
{\bf Proof.}  Let $\alpha=ax(g)$ and $\beta=ax(g)$ and let $\gamma$ be the hyperbolic
line perpendicular to $\alpha$ and $\beta$.  Such a line $\gamma$ exists because $\G$ is
discrete and so
$\alpha\cap\beta$ is either empty or a finite point of hyperbolic space.  Let $\phi_1$
be the half turn whose axis is $\gamma$.  There is a half turn $\phi_2$ whose axis is
perpendicular to $\alpha$ and such that $\phi_2\phi_1=f$.  Similarly there is $\phi_3$
with $\phi_3\phi_1=g$.  It is easy to see  that $\langle
\phi_1,\phi_2,\phi_3 \rangle$ contain $\G$ with index at most two and therefore
$\G^\prime$ is Kleinian. 

Next,  the point  $x_0\in \IH^3$ such that $m_\G=\rho(x_0,f(x_0))=\rho(x_0,g(x_0))$ lies
on the line $\gamma$ between the two axes.  To see this note that the set of all points
$y$ such that for a given constant $\delta$ we have $\rho(y,f(y))=\delta$ forms a
hyperbolic cylinder about $\alpha$,  a fact easily observed if one conjugates $f$ so
that its fixed points are $0,\infty\in \oC$.  These cylinders bound geodesically
convex regions and the two cylinders $\{y:\rho(y,f(y))=m_\G\}$ and
$\{y:\rho(y,g(y))=m_\G\}$ meet at a unique point of $\gamma$.  Now as $x_0\in \gamma$
we have
\begin{eqnarray}
\rho(x_0,\phi_2(x_0)) & = & 0\\
\rho(x_0,\phi_2(x_0)) & = & \rho(x_0,\phi_2\phi_1(x_0)) = \rho(x_0,f(x_0)) = m_\G \\
\rho(x_0,\phi_3(x_0)) & = & \rho(x_0,\phi_3\phi_1(x_0)) = \rho(x_0,g(x_0)) = m_\G 
\end{eqnarray}
so of course $m_{\G}^{\prime} \leq m_\G$.  That strict inequality holds follows since
for a Kleinian group generated by three elliptics of order 2 the Margulis constant is
the diameter of the smallest hyperbolic ball which is tangent to all three axes. $\Box$

\bigskip

In fact Lemma 2.7 of \cite{GM1} can be used to find the Margulis constant for any group
generated by a pair of loxodromic transformations in terms of the traces of the
elements and their commutators.  However the solution involves solving highly nonlinear
equations and is only practical numerically.  Inequalities for Kleinian groups, such
as J\o rgensen's inequality, can be used together with these formulae to give lower
bounds for the Margulis
constant in the general case.  However at present our methods are far from sharp and we
shall return to this problem elsewhere.

Further note that as limiting cases of the above it is easily seen that any two
generator Kleinian group lies with index at most two in a Kleinian group generated by
three elliptics of order 2.

\subsection{An Example}

In this section we give an example of a Kleinian group generated by three elliptic
elements of order 2 with a Margulis constant which we conjecture to be smallest
possible among all Kleinian groups.  The group in question arises from an arithmetic
Kleinian group (a $\IZ_2$--extension of the $3-5-3$ Coxeter group, $\Gamma_1$ in our list).  This group has
the smallest co-volume among all arithmetic groups and is conjectured to be of minimal
co-volume among all Kleinian  groups.  Unfortunately the calculation is rather messy
involving the numerical solution of a number of simultaneous nonlinear equations.  We
give the example and leave the computational details to the reader.

Consider the hyperbolic tetrahedron with vertices $A,B,C,D$ and dihedral angles at the
edges as below.
\begin{itemize}
\item dihedral angle at $\overline{AB}$ is $\pi/3$,
\item dihedral angle at $\overline{AC}$ is $\pi/2$,
\item dihedral angle at $\overline{AD}$ is $\pi/2$,
\item dihedral angle at $\overline{BC}$ is $\pi/5$,
\item dihedral angle at $\overline{BD}$ is $\pi/2$,
\item dihedral angle at $\overline{CD}$ is $\pi/3$,
\end{itemize}
Let $\ell$ be the hyperbolic line which is perpendicular to  $\overline{AD}$ and 
$\overline{BC}$.  We
note first that the angle
$\angle ADB =\pi/3$ and that $\ell$ and
$\overline{AD}$ meet at right angles.  Also the angle between the plane spanned by
$AD$ and $BD$ and the plane spanned by $AD$ and $\ell$ is $\pi/4$.

The group we want is the group generated by the three half turns
$\phi_1,\phi_2,\phi_3$ whose axes are  $\overline{AD}$,   $\overline{BD}$ and $\ell$
respectively.  That this group is discrete is easily seen from the $\IZ_2$ symmetry of
the Coxeter tetrahedron about
$\ell$ and the fact that the reflection group in the faces of this tetrahedron is
discrete.

We want the radius of the unique ball tangent to the axes of $\phi_i$, $i=1,2,3$. To
find this we subdivide the tetrahedron about the center of this ball and use elementary
hyperbolic trigonometry and symmetry to produce a number of different equations for
$r$ which we then solve numerically.   We find

 \subsection{The $(p,q,r)$ groups.}
 
 We now go about tabulating all these groups.  We include all the data necessary to computer the Margulis constant for each group.  The first three  angles given are those the edges of the triangle meet and therefore must be those occurring in our list of angles possibly occurring in the elementary groups.  We then indicate  the vertex stabilisers and an associated tetrahedral reflection group.
 
 \subsubsection{Dihedral angles submultiples of $\pi$}
{\scriptsize
\begin{center}
\begin{tabular}{|c|c|c|c|c|c|c|c|c|}
\multicolumn{9}{c}{\bf  $(p,q,r)$--triangles ($p\geq 3$).  Dihedral angles submultiples of $\pi$ }\\
\hline
$\#$ & $m_\G$ & $(p,q,r)$ & Angle & Angle & Angle & $\Sigma$ dihedral & Discrete & Vertex Structure\\
\hline
$1$ & $0.9624$ & $(3,3,3)$ & $0$ & $0$ & $0$ & $\pi$ & $\Gamma_{31}$ & $\Delta_{(3,3,3)},\Delta_{(3,3,3)},\Delta_{(3,3,3)}$\\ 
\hline
$2$ & $0.7564$ & $(3,3,3)$ & $0$ & $0$ & $1.230 $ & $7\pi/6$ & $\Gamma_{32}$ & $ A_4,\Delta_{(3,3,3)},\Delta_{(3,3,3)}$\\
\hline
$3$ & $0.4435$ & $(3,3,3)$ & $0$ & $1.230 $ & $1.230$ & $4\pi/3$ & $\Gamma_{27}$ & $A_4,A_4,\Delta_{(3,3,3)}$\\
\hline
$4$ & $0.6329$ & $(3,3,4)$ & $0$ & $.955$ & $.955$ & $4\pi/3$ &  $\Gamma_{28}$ & $\Delta_{(3,3,3)},S_4,S_4$\\
\hline
$5$ & $0.7362$ & $(3,3,5)$ & $0$ & $.652$ & $.652 $ & $4\pi/5$ & $\Gamma_{29}$ & $\Delta_{(3,3,3)},A_5,A_5$\\
\hline
$6$ & $0.3796$ & $(3,3,5)$ & $1.230$ & $.652$ & $.652$ & $3\pi/2$ & $\Gamma_9$ & $ A_4,A_5,A_5 $\\
\hline
$7$ & $0.7983$ & $(3,3,6)$ & $0$ & $0$ & $0$ & $4\pi/3$ & $\Gamma_{30}$ &
$\Delta_{(3,3,3)},\Delta_{(2,3,6)},\Delta_{(2,3,6)}$\\ 
\hline
$8$ & $0.6349$ & $(3,3,6)$ & $1.230$ & $0$ & $0$ & $3\pi/2$ & $\Gamma_{24}$ & $ A_4,\Delta_{(2,3,6)},\Delta_{(2,3,6)}$ \\
\hline
$9$ & $0.7953$ & $(3,4,4)$ & $.955$ & $.955$ & $0$ & $3\pi/2$ & $\Gamma_{25}$ & 
$ \Delta_{(2,4,4)},\Delta_{(2,4,4)},\Delta_{(2,4,4)}$
\\
\hline
$10$ & $0.5777$ & $(4,4,4)$ & $0$ & $0$ & $0$ & $3\pi/2$ & $\Gamma_{26}$ & $
\Delta_{(2,4,4)},\Delta_{(2,4,4)},\Delta_{(2,4,4)}$ \\
\hline
\end{tabular}
\end{center}}

{\scriptsize
\begin{center}
\begin{tabular}{|c|c|c|c|c|c|c|c|}
\multicolumn{8}{c}{\bf  (2,3,3)--triangles. Dihedral angles submultiples of $\pi$ }\\
\hline
$\#$ & $m_\G$ &$3-3$&$2-3$ &$2-3$ & dihedral angles & discrete & Vertex Structure\\
\hline
$1$ & 1.0050 & $ 0 $&$ 0 $&$ 0 $& $\pi/3 ,\pi/6 , \pi/6 $ & $\Delta_{(3,3,6)}$ &
$\Delta_{(3,3,3)},\Delta_{(2,3,6)},\Delta_{(2,3,6)}$
\\
\hline
$2$ & 0.9841 & $ 0 $&$ .364 $&$ 0 $& $\pi/3 ,\pi/5 , \pi/6 $ &$\Delta_{(3,5,6)}$ &$\Delta_{(3,3,3)},A_5,\Delta_{(2,3,6)}$ \\
\hline
$3$ & 0.9465 & $ 0 $&$ .615 $&$ 0 $& $\pi/3 ,\pi/4 , \pi/6 $ & $\Delta_{(3,4,6)}$&$\Delta_{(3,3,3)},S_4,\Delta_{(2,3,6)}$ \\
\hline
$ 4 $ & 0.8689 & $ 0 $&$ .955 $&$ 0 $& $\pi/3 ,\pi/3 , \pi/6 $ & $\Delta_{(3,3,6)}$ &$\Delta_{(3,3,3)},A_4,\Delta_{(2,3,6)}$
\\
\hline
$ 5 $ & 0.6687 & $ 0 $&$ \pi/2 $&$ 0 $&  $\pi/3 ,\pi/2 , \pi/6 $ & $\Gamma_{30}$ &$\Delta_{(3,3,3)},D_3,\Delta_{(2,3,6)}$ \\
\hline
$6 $ & 0.9624 & $ 0 $&$ .364 $&$ .364 $&  $\pi/3 ,\pi/5 , \pi/5 $ &$\Delta_{(3,5,5)}$ &$\Delta_{(3,3,3)},A_5,A_5$ \\
\hline
$ 7 $ & 0.9233 & $ 0 $&$ .615 $&$ .364 $& $\pi/3 ,\pi/4 , \pi/5 $ &$\Delta_{(3,4,5)}$ & $\Delta_{(3,3,3)},S_4,A_5$ \\
\hline
$ 8 $ & 0.8423 & $ 0 $&$ .955 $&$ .364 $& $\pi/3 ,\pi/3 , \pi/5 $ & $\Delta_{(3,3,5)}$& $\Delta_{(3,3,3)},A_4,A_5$ \\
\hline
$ 9 $ & 0.6309 & $ 0 $&$ \pi/2 $&$ .364 $&  $\pi/3 ,\pi/2 , \pi/5 $ & $\Gamma{29}$ &$\Delta_{(3,3,3)},D_3,A_5$ \\
\hline
$ 10 $ & 0.8813 & $ 0 $&$ .615 $&$ .615 $& $\pi/3 ,\pi/4 , \pi/4 $ & $\Delta_{(3,4,4)}$ & $\Delta_{(3,3,3)},S_4,S_4$ \\
\hline
$ 11 $ & 0.7941 & $ 0 $&$ .955 $&$ .615 $&  $\pi/3 ,\pi/3 , \pi/4 $ &$\Delta_{(3,3,4)}$ & $\Delta_{(3,3,3)},A_4,S_4$ \\
\hline
$ 12 $ & 0.5622 & $ 0 $&$ \pi/2 $&$ .615 $&  $\pi/3 ,\pi/2 , \pi/4 $ & $\Gamma_{28}$ & $\Delta_{(3,3,3)},D_3,S_4$ \\
\hline
$ 13 $ & 0.6931 & $ 0 $&$ .955 $&$ .955 $&  $\pi/3 ,\pi/3 , \pi/3 $ & $\Gamma_{32}$ & $\Delta_{(3,3,3)},A_4,A_4$ \\
\hline
$ 14 $ & 0.4170 & $ 0 $&$ \pi/2 $&$ .955 $& $\pi/3 ,\pi/2 , \pi/3 $ & $\Gamma_{27}$ &  $\Delta_{(3,3,3)},D_3,A_4$ \\
\hline
$ 15 $ & 0.7872 & $ 1.230 $&$ 0 $&$ 0 $& $\pi/2 ,\pi/6 , \pi/6 $ & $\Delta_{(2,6,6)}$&$A_4,\Delta_{(2,3,6)},\Delta_{(2,3,6)}$
\\
\hline
$ 16 $ & 0.7673 & $ 1.230 $&$ .364 $&$ 0 $& $\pi/2 ,\pi/5 , \pi/6 $ & $\Delta_{(2,5,6)}$&$A_4,A_5,\Delta_{(2,3,6)}$ \\
\hline
$ 17 $ & 0.7123 & $ 1.230 $&$ .615 $&$ 0 $& $\pi/2 ,\pi/4 , \pi/6 $ & $\Delta_{(2,4,6)}$&$A_4,S_4,\Delta_{(2,3,6)}$ \\
\hline
$ 18 $ & 0.5961 & $ 1.230 $&$ .955 $&$ 0 $& $\pi/2 ,\pi/3 , \pi/6 $ & $\Gamma_{10}$ &$A_4,A_4,\Delta_{(2,3,6)}$ \\
\hline
$ 19 $ & 0.2671 & $ 1.230 $&$ \pi/2 $&$ 0 $&  $\pi/2 ,\pi/2 , \pi/6 $ & $\Gamma_{17}$ &$A_4,D_3,\Delta_{(2,3,6)}$ \\
\hline
$ 20 $ & 0.7344 & $ 1.230 $&$ .364 $&$ .364 $& $\pi/2 ,\pi/5 , \pi/5 $ &$\Delta_{(2,5,5)}$ & $A_4,A_5,A_5$ \\
\hline
$ 21 $ & 0.6743 & $  1.230 $&$ .615 $&$ .364 $& $\pi/2 ,\pi/4 , \pi/5 $ & $\Delta_{(2,4,5)}$& $A_4,S_4,A_5$ \\
\hline
$ 22 $ & 0.5438 & $ 1.230 $&$ .955 $&$ .364 $&  $\pi/2 ,\pi/3 , \pi/5 $ & $\Gamma_{5}$ &$A_4,A_4,A_5$ \\
\hline
$ 23 $ & 0.6034 & $ 1.230 $&$ .615 $&$ .615 $ &  $\pi/2 ,\pi/4 , \pi/4 $ & $\Gamma_{14}$ &$A_4,S_4,S_4$ \\
\hline
$ 24 $ & 0.4397 & $ 1.230 $&$ .955 $&$ .615 $& $\pi/2 ,\pi/3 , \pi/4 $ & $\Gamma_{4}$ &  $A_4,A_4,S_4$ \\
\hline
\end{tabular}
\end{center}}
 
{\scriptsize
\begin{center}
\begin{tabular}{|c|c|c|c|c|c|c|c|}
\multicolumn{8}{c}{\bf   (2,3,4)--triangles. Dihedral angles submultiples of $\pi$ }\\
\hline
$ \# $ &$m_\G$ &$3-4$&$2-4$ &$2-3$ & dihedral angles & discrete & Vertex Structure\\
\hline
$ 1 $ & 0.8040 & $ .955 $&$ 0 $&$ 0 $& $\pi/2, \pi/4 , \pi/6 $& $\Delta_{(2,4,6)}$ & $S_4,\Delta_{(2,3,6)},\Delta_{(2,3,6)}$
\\
\hline
$ 2 $ & 0.7955 & $ .955 $&$ 0 $&$ .364 $& $ \pi/2, \pi/4 , \pi/5 $& $\Delta_{(2,4,5)}$ &$S_4,A_5,\Delta_{(2,4,4)}$ \\
\hline
$ 3 $ & 0.7478 & $ .955 $&$ 0 $&$ .615 $& $ \pi/2, \pi/4 , \pi/4 $& $\Gamma_{15}$ &$S_4,S_4,\Delta_{(2,4,4)}$ \\
\hline
$ 4 $ & 0.6494 & $ .955 $&$ 0 $&$ .955 $& $ \pi/2, \pi/4 , \pi/3 $& $\Gamma_{14}$ &$S_4,A_4,\Delta_{(2,4,4)}$ \\
\hline
$ 5 $ & 0.3887 & $ .955 $&$ 0 $&$ \pi/2 $& $ \pi/2, \pi/4 , \pi/2 $& $\Gamma_{19}$ &$S_4,D_3,\Delta_{(2,4,4)}$ \\
\hline
$ 6 $ & 0.7138 & $ .955 $&$ \pi/4 $&$ 0 $& $\pi/2, \pi/3 , \pi/6 $& $\Gamma_{11}$ &$S_4,S_4,\Delta_{(2,3,6)}$ \\
\hline
$ 7 $ & 0.6792 & $ .955 $&$ \pi/4 $&$ .364 $& $\pi/2, \pi/3 , \pi/5 $& $\Gamma_7$ &$S_4,S_4,A_5$ \\
\hline
$ 8 $ & 0.6155 & $ .955 $&$ \pi/4 $&$ .615 $& $ \pi/2, \pi/3 , \pi/4 $& $\Gamma_6$ &$S_4,S_4,S_4$ \\
\hline
$ 9 $ & 0.4749 & $ .955 $&$ \pi/4 $&$ .955 $& $\pi/2, \pi/3 , \pi/3 $& $\Gamma_4$ & $S_4,S_4,A_4$ \\
\hline
$ 10 $ & 0.4098 & $ .955 $&$ \pi/2 $&$ 0 $& $\pi/2, \pi/2 , \pi/6 $& $\Gamma_{20}$ &$S_4,D_4,\Delta_{(2,3,6)}$ \\
\hline
$ 11 $ & 0.3231 & $ .955 $&$ \pi/2 $&$ 0.364 $& $\pi/2, \pi/2 , \pi/5 $& $\Gamma_2$ &$S_4,D_4,A_5$ \\
\hline
\end{tabular}
\end{center}}

{\scriptsize 
\begin{center}
\begin{tabular}{|c|c|c|c|c|c|c|c|}
\multicolumn{8}{c}{\bf  (2,3,5)--triangles. Dihedral angles submultiples of $\pi$ }\\
\hline
$ \# $ & $m_\G$ &$3-5$&$2-5$ &$2-3$ & dihedral angles & discrete & Vertex Structure\\
\hline
$ 1 $ & 0.7855 & $ .652 $&$ .553 $&$ 0 $& $\pi/2 ,\pi/3 ,\pi/6 $&  $\Gamma_{12}$ & $A_5,A_5,\Delta_{(2,3,6)}$ \\
\hline
$ 2 $ & 0.7589 & $ .653 $&$ .553 $&$ .364 $& $\pi/2 ,\pi/3 ,\pi/5 $&  $\Gamma_8$ &$A_5,A_5,A_5$ \\
\hline
$ 3$ & 0.7115 & $ .652 $&$ .553 $&$ .615 $& $\pi/2 ,\pi/3 ,\pi/4 $&  $\Gamma_7$ &$A_5,A_5,S_4$ \\
\hline
$ 4 $ & 0.6147 & $ .652 $&$ .553 $&$ .955 $& $\pi/2 ,\pi/3 ,\pi/3 $& $\Gamma_5$  &$A_5,A_5,A_4$ \\
\hline
$ 5 $ & 0.3541 & $ .652 $&$ .553 $&$ \pi/2 $&  $\pi/2 ,\pi/3 ,\pi/2 $& $\Gamma_1$  &$A_5,A_5,D_3$ \\
\hline
$ 6 $ & 0.5014 & $ .652 $&$ \pi/2 $&$ 0 $&  $\pi/2 ,\pi/2 ,\pi/6 $& $\Gamma_{21}$  &$A_5,D_5,\Delta_{(2,3,6)}$ \\
\hline
$ 7 $ & 0.4481 & $ .652 $&$ \pi/2 $&$ .364 $&  $\pi/2 ,\pi/2 ,\pi/5 $& $\Gamma_3$  &$A_5,D_5,A_5$ \\
\hline
$ 8 $ & 0.3456 & $ .652 $&$ \pi/2 $&$ .615 $&  $\pi/2 ,\pi/2 ,\pi/4 $& $\Gamma_2$  &$A_5,D_5,S_4$ \\
\hline
\end{tabular}
\end{center}}

{\scriptsize 
\begin{center}
\begin{tabular}{|c|c|c|c|c|c|c|c|}
\multicolumn{8}{c}{\bf   (2,3,6)--triangles. Dihedral angles submultiples of $\pi$ }\\
\hline
$ \# $ & $m_\G$ &$2-3$&$2-6$ &$3-6$ & dihedral angles & discrete  & Vertex Structure\\
\hline
$ 1 $ & 0.8314 & $ 0 $&$ 0 $&$ 0 $& $\pi/6 ,\pi/3 ,\pi/2 $& $\Gamma_{13}$ & 
$\Delta_{(2,3,6)},\Delta_{(2,3,6)},\Delta_{(2,3,6)}$ \\
\hline
$ 2 $ & 0.8152 & $ .364 $&$ 0 $&$ 0 $& $\pi/5 ,\pi/3 ,\pi/2 $&  $\Gamma_{12}$ &$A_5,\Delta_{(2,3,6)},\Delta_{(2,3,6)}$ \\
\hline
$ 3 $ & 0.7860 & $ .615 $&$ 0 $&$ 0 $& $\pi/4 ,\pi/3 ,\pi/2 $& $\Gamma_{11}$ &$S_4,\Delta_{(2,3,6)},\Delta_{(2,3,6)}$ \\
\hline
$ 4 $ & 0.7249 & $ .955 $&$ 0 $&$ 0 $& $\pi/3 ,\pi/3 ,\pi/2 $& $\Gamma_{10}$ &$A_4,\Delta_{(2,3,6)},\Delta_{(2,3,6)}$ \\
\hline
$ 5 $ & 0.5636 & $ \pi/2 $&$ 0 $&$ 0 $&  $\pi/2 ,\pi/3 ,\pi/2 $& $\Gamma_{18}$ & $D_3,\Delta_{(2,3,6)},\Delta_{(2,3,6)}$ \\
\hline
$ 6 $ & 0.5120 & $ 0 $&$ \pi/2 $&$ 0 $&  $\pi/6 ,\pi/2 ,\pi/2 $& $\Gamma_{22}$ & $D_6,\Delta_{(2,3,6)},\Delta_{(2,3,6)}$ \\
\hline
$ 7 $ & 0.6572 & $ .364 $&$ \pi/2 $&$ 0 $&  $\pi/5 ,\pi/2 ,\pi/2 $& $\Gamma_{21}$ & $D_6,\Delta_{(2,3,6)},A_5$ \\
\hline
$ 8 $ & 0.5879 & $ .615 $&$ \pi/2 $&$ 0 $&  $\pi/4 ,\pi/2 ,\pi/2 $& $\Gamma_{20}$ & $D_6,\Delta_{(2,3,6)},S_4$ \\
\hline
$ 9 $ & 0.4393 & $ .955 $&$ \pi/2 $&$ 0 $&  $\pi/3 ,\pi/2 ,\pi/2 $& $\Gamma_{17}$ & $D_6,\Delta_{(2,3,6)},A_4$ \\
\hline
\end{tabular}
\end{center}}

{\scriptsize
\begin{center}
\begin{tabular}{|c|c|c|c|c|c|c|c|}
\multicolumn{8}{c}{\bf (2,4,4)--triangles. Dihedral angles submultiples of $\pi$ }\\
\hline
$ \# $ &$m_\G$&$2-4$&$2-4$ &$4-4$ & dihedral angles & discrete & Vertex Structure\\
\hline
$ 1 $ & 0.8813 & $ 0 $&$ 0 $&$ 0 $& $\pi/4 ,\pi/4 ,\pi/2 $ & $\Gamma(16)$ &
$\Delta_{(2,4,4)},\Delta_{(2,4,4)},\Delta_{(2,4,4)}$ \\
\hline
$ 2 $ & 0.5777 & $ \pi/2 $&$ 0 $&$ 0 $& $\pi/3 ,\pi/4 ,\pi/2 $ &  $\Gamma_{23}$ & $D_4,\Delta_{(2,4,4)},\Delta_{(2,4,4)}$ \\
\hline
$ 3 $ & 0.7953 & $ \pi/4 $&$ 0 $&$ 0 $& $\pi/3 ,\pi/4 ,\pi/2 $ &  $\Gamma_{15}$ & $S_4,\Delta_{(2,4,4)},\Delta_{(2,4,4)}$ \\
\hline
$ 4 $ & 0.4232 & $ \pi/2 $&$ \pi/4 $&$ 0 $& $\pi/3 ,\pi/3 ,\pi/2 $ &  $\Gamma_{19}$  & $D_4, S_4, \Delta_{(2,4,4)}$ \\
\hline
$ 5 $ & 0.6931 & $ \pi/4 $&$ \pi/4 $&$ 0 $& $\pi/3 ,\pi/3 ,\pi/2 $ &  $\Gamma_{14}$  & $S_4,S_4,\Delta_{(2,4,4)}$ \\
\hline
\end{tabular}
\end{center}}

{\scriptsize
\begin{center}
\begin{tabular}{|c|c|c|c|c|c|c|c|}
\multicolumn{8}{c}{\bf  (2,2,3)--triangles. Dihedral angles submultiples of $\pi$ }\\
\hline
$\#$ & $m_\G$ &$2-3$&$2-3$ &$2-2$ & dihedral angles & discrete & Vertex Structure\\
\hline
$1$ & 0.7095 & $ 0 $&$ 0 $&$ \pi/n\;\;(2\leq n\leq \infty) $&$\pi/6,\pi/6,\pi/n $&
$\Delta_{(n,6,6)}$&$\Delta_{(2,3,6)},\Delta_{(2,3,6)},D_n$\\
\hline
$2$ & 0.6694 & $ 0 $&$ .364 $&$ \pi/n \;\; (2\leq n\leq \infty)  $&  $\pi/6 ,\pi/5 ,\pi/n $& $\Delta_{(5,6,n)}$
&$\Delta_{(2,3,6)},A_5,D_n$ \\
\hline
$3$ & 0.5966 & $ 0 $&$ .615 $&$ \pi/n \;\; (2\leq n\leq \infty)  $&  $\pi/6 ,\pi/4 ,\pi/n $&$\Delta_{(4,6,n)}$
&$\Delta_{(2,3,6)},S_4,D_n$ \\
\hline
$4$ & 0.4429 & $ 0 $&$ .955 $&$ \pi/2 $&  $\pi/6 ,\pi/3 ,\pi/2 $& $\Gamma_{24}$
&$\Delta_{(2,3,6)},A_4,D_2$ \\
\hline
$5$ & 0.6941 & $ 0 $&$ .955 $&$ \pi/n \;\; (3\leq n\leq \infty)  $&  $\pi/6 ,\pi/3 ,\pi/n $& $\Delta_{(3,n,6)}$
&$\Delta_{(2,3,6)},A_4,D_n$ \\
\hline
$6$ & 0.3908 & $ 0 $&$ \pi/2 $&$ \pi/3 $&   $\pi/6 ,\pi/2 ,\pi/3 $& $\Gamma_{18}$ &$\Delta_{(2,3,6)},D_3,D_3$ \\
\hline
$7$ & 0.5274 & $ 0 $&$ \pi/2 $&$ \pi/n \;\; (4\leq n\leq \infty)  $&   $\pi/6 ,\pi/2 ,\pi/n $& $\Delta_{(2,6,n)}$
&$\Delta_{(2,3,6)},D_3,D_n$ \\
\hline
$8$ & 0.6146 & $ .364 $&$ .364 $&$ \pi/n \;\; (2\leq n\leq \infty)  $ &$\pi/5 ,\pi/5 ,\pi/n $& $\Delta_{(5,5,n)}$ &
$A_5,A_5,D_n$
\\
\hline
$9$ & 0.5349 & $ .364 $&$ .615 $&$ \pi/n \;\; (2\leq n\leq \infty)  $& $\pi/5 ,\pi/4 ,\pi/n $& $\Delta_{(4,5,n)}$ &
$A_5,S_4,D_n$ \\
\hline
$10$ & 0.3490 & $ .364 $&$ .955 $&$ \pi/2  $& $\pi/5 ,\pi/3 ,\pi/n $& $\Gamma_9$ & $A_5,A_4,D_2$ \\
\hline
$11$ & 0.6496 & $ .364 $&$ .955 $&$ \pi/n
\;\; (3\leq n\leq \infty)  $&
$\pi/5 ,\pi/3 ,\pi/n $&
$\Delta_{(3,5,n)}$ & $A_5,A_4,D_n$
\\
\hline
$12$ & 0.2768 & $ .364 $&$ \pi/2 $&$ \pi/3 $&  $\pi/5 ,\pi/2 ,\pi/3 $& $\Gamma_1$& $A_5,D_3,D_3$ \\
\hline
$13$ & 0.4610 & $ .364 $&$ \pi/2 $&$ \pi/n
\;\; (4\leq n\leq \infty)  $& 
$\pi/5 ,\pi/2 ,\pi/n $&
$\Delta_{(2,5,n)}$& $A_5,D_3,D_n$
\\
\hline
$14$ & 0.4299 & $ .615 $&$ .615 $&$ \pi/2
$&
$\pi/4 ,\pi/4 ,\pi/2 $&
$\Gamma_{25}$ & $S_4,S_4,D_2$ \\
\hline
$15$ & 0.7048 & $ .615 $&$ .615 $&$ \pi/n
\;\; (3\leq n\leq \infty)  $&
$\pi/4 ,\pi/4 ,\pi/n $& 
$\Delta_{(4,4,n)}$ & $S_4,S_4,D_n$
\\
\hline
$16$ & 0.5659 & $ .615 $&$ .955 $&$ \pi/n
\;\; (3\leq n\leq \infty)  $&
$\pi/4 ,\pi/3 ,\pi/n $&
$\Delta_{(3,4,n)}$ & $S_4,A_4,D_n$
\\
\hline
$17$ & 0.3157 & $ .615 $&$ \pi/2 $&$ \pi/4
$&  $\pi/4 ,\pi/2 ,\pi/4
$&$\Gamma_{19}$ & $S_4,D_3,D_4$ \\
\hline
$18$ & 0.4305 & $ .615 $&$ \pi/2 $&$ \pi/n
\;\; (5\leq n\leq \infty) $& 
$\pi/4 ,\pi/2 ,\pi/n
$&$\Delta_{(2,4,n)}$ &
$S_4,D_3,D_n$ \\
\hline
$19$ & 0.3517 & $ .955 $&$ .955 $&$ \pi/3
$&
$\pi/3 ,\pi/3 ,\pi/3 $& $\Gamma_{27} $ & $A_4,A_4,D_3$ \\
\hline
$20$ & 0.5372 & $ .955 $&$ .955 $&$ \pi/n
\;\; (4\leq n\leq \infty)  $&
$\pi/3 ,\pi/3 ,\pi/n $&
$\Delta_{(3,3,n)}$ & $A_4,A_4,D_n$
\\
\hline
$21$ & 0.2309& $ .955 $&$ \pi/2 $&$ \pi/6
$&  $\pi/3 ,\pi/2 ,\pi/6 $& $\Gamma_{17}$ & $A_4,D_3,D_6$ \\
\hline
$22$ & 0.30345& $ .955 $&$ \pi/2 $&$ \pi/n
\;\; (7\leq n\leq \infty)  $& 
$\pi/3 ,\pi/2 ,\pi/n $&
$\Delta_{(2,3,n)}$ & $A_4,D_3,D_n$
\\
\hline
\end{tabular}
\end{center}}

{\scriptsize
\begin{center}
\begin{tabular}{|c|c|c|c|c|c|c|c|}
\multicolumn{8}{c}{\bf  (2,2,4)--triangles. Dihedral angles submultiples of $\pi$ }\\
\hline
$\# $ & $m_\G$ &$2-4$&$2-4$ &$2-2$ & dihedral angles & discrete & Vertex Structure\\
\hline
$1$ & 0.5777 & $ 0 $&$ 0 $&$ \pi/2 $& $\pi/4 ,\pi/4 ,\pi/2 $ & $\Gamma_{26}$ & $\Delta_{(2,4,4)},\Delta_{(2,4,4)},D_2$ 
\\
\hline
$2$ & 0.7337 & $ 0 $&$ 0 $&$ \pi/n \;\; (3\leq
n
\leq \infty) $& $\pi/4 ,\pi/4 ,\pi/n $ &
$\Delta_{(4,4,n)}$ &
$\Delta_{(2,4,4)},\Delta_{(2,4,4)},D_n$ 
\\
\hline
$3$ & 0.4738 & $ 0 $&$ \pi/4 $&$ \pi/2 $& $\pi/4 ,\pi/3 ,\pi/2 $& $\Gamma_{25}$ &
$\Delta_{(2,4,4)},S_4,D_2$  \\
\hline
$4$ & 0.6931 & $ 0 $&$ \pi/4 $&$ \pi/n \;\; (3\leq n \leq \infty) $& $\pi/4 ,\pi/3 ,\pi/n $& $\Delta_{(3,4,n)}$ &
$\Delta_{(2,4,4)},S_4,D_n$  \\
\hline
$5$ & 0.3574 & $ 0 $&$ \pi/2 $&$ \pi/3 $& $\pi/4 ,\pi/2 ,\pi/3 $& $\Gamma_{19}$  & $\Delta_{(2,4,4)},D_4,D_3$  \\
\hline
$6$ & 0.4877 & $ 0 $&$ \pi/2 $&$ \pi/4 $& $\pi/4 ,\pi/2 ,\pi/4 $& $\Gamma_{23}$  & $\Delta_{(2,4,4)},D_4,D_4$  \\
\hline
$7$ & 0.5504 & $ 0 $&$ \pi/2 $&$ \pi/n \;\;  (5\leq n \leq \infty) $& $\pi/4 ,\pi/2 ,\pi/n $& $\Delta_{(2,4,n)}$  &
$\Delta_{(2,4,4)},D_4,D_n$  \\
\hline
$8$ & 0.5283 & $ \pi/4 $&$ \pi/4 $&$ \pi/3 $& $\pi/3 ,\pi/3 ,\pi/3$ & $\Gamma_{28}$ & $S_4,S_4,D_3$ 
\\
\hline
$9$ & 0.6441 & $ \pi/4 $&$ \pi/4 $&$ \pi/n \;\;
(4\leq n \leq \infty) $& $\pi/3 ,\pi/3
,\pi/n$ & $\Delta_{(3,3,n)}$ &
$S_4,S_4,D_n$ 
\\
\hline
$10$ & 0.2874 & $ \pi/4 $&$ \pi/2 $&$ \pi/5 $& $\pi/3 ,\pi/2 ,\pi/5 $& $\Gamma_2$ & $S_4,D_4,D_5$  \\
\hline
$11$ & 0.3640 & $ \pi/4 $&$ \pi/2 $&$ \pi/6 $&
$\pi/3 ,\pi/2 ,\pi/6 $& $\Gamma_{20}$ &
$S_4,D_4,D_6$  \\
\hline
$12$ & 0.4055 & $ \pi/4 $&$ \pi/2 $&$ \pi/n
\;\;  (7\leq n \leq \infty) $& $\pi/3
,\pi/2 ,\pi/n $&$\Delta_{(2,3,n)}$ &
$S_4,D_4,D_n$ 
\\
\hline
\end{tabular}
\end{center}}

{\scriptsize
\begin{center}
\begin{tabular}{|c|c|c|c|c|c|c|c|}
\multicolumn{8}{c}{\bf  (2,2,5)--triangles. Dihedral angles submultiples of $\pi$ }\\
\hline
$\#$ & $m_\G$ &$2-5$&$2-5$ &$2-2$ & dihedral angles & discrete & Vertex Structure\\
\hline
$1$ & 0.4064 & $ .55357 $&$ .55357 $&$ \pi/2 $& $\pi/3 ,\pi/3 ,\pi/2 $ & $\Gamma_9$ &  $A_5,A_5,D_2$ \\
\hline
$2$ & 0.6300 & $ .55357 $&$ .55357 $&$ \pi/3 $&
$\pi/3 ,\pi/3 ,\pi/3 $ & $\Gamma_{29}$ & 
$A_5,A_5,D_3$ \\
\hline
$3$ & 0.7157 & $ .55357 $&$ .55357 $&$ \pi/n,
\;\; (4\leq n \leq \infty) $& $\pi/3
,\pi/3 ,\pi/n $ & $\Delta(3,3,n)$  & 
$A_5,A_5,D_n$ \\
\hline
$4$ & 0.3126 & $ .55357 $&$ \pi/2 $&$ \pi/4 $& $\pi/3 ,\pi/2 ,\pi/4 $ & $\Gamma_2$ & $A_5,D_5,D_4$ \\
\hline
$5$ & 0.4042 & $ .55357 $&$ \pi/2 $&$ \pi/5 $&
$\pi/3 ,\pi/2 ,\pi/5 $ & $\Gamma_3$ &
$A_5,D_5,D_5$ \\
\hline
$6$ & 0.4515 & $ .55357 $&$ \pi/2 $&$ \pi/6 $&
$\pi/3 ,\pi/2 ,\pi/6 $ & $\Gamma_{21}$ &
$A_5,D_5,D_6$ \\
\hline
$7$ & 0.4798 & $ .55357 $&$ \pi/2 $&$ \pi/n,
\;\; (7\leq n \leq \infty) $& $\pi/3
,\pi/2 ,\pi/n $ & $\Delta(2,3,n)$ &
$A_5,D_5,D_n$
\\
\hline
\end{tabular}
\end{center}}

{\scriptsize
\begin{center}
\begin{tabular}{|c|c|c|c|c|c|c|c|}
\multicolumn{8}{c}{\bf  (2,2,6)--triangles. Dihedral angles submultiples of $\pi$}\\
\hline
$\#$ & $m_\G$ &$2-6$&$2-6$ &$2-2$ & dihedral angles & discrete & Vertex Structure\\
\hline
$1$ & 0.5942 & $ 0 $&$ 0 $&$ \pi/2 $ & $\pi/3 ,\pi/3 ,\pi/2 $ & $\Gamma_{24}$ & $\Delta_{(2,3,6)},\Delta_{(2,3,6)},D_2$ \\
\hline
$2$ & 0.7389 & $ 0 $&$ 0 $&$ \pi/3 $ & $\pi/3 ,\pi/3 ,\pi/3 $ & $\Gamma_{30}$ & $\Delta_{(2,3,6)},\Delta_{(2,3,6)},D_3$ \\
\hline
$3$ & 0.7935 & $ 0 $&$ 0 $&$ \pi/n,\; 4\leq n\leq \infty $  & $\pi/3 ,\pi/3 ,\pi/n $& $\Delta_{(3,3,n)}$
& $\Delta_{(2,3,6)},\Delta_{(2,3,6)},D_n$ \\
\hline
$4$ & 0.3003 & $ 0 $&$ \pi/2 $&$ \pi/3 $&   $\pi/3 ,\pi/2 ,\pi/3 $ & $\Gamma_{17}$ & $\Delta_{(2,3,6)},D_6,D_3 $
\\
\hline
$5$ & 0.4187 & $ 0 $&$ \pi/2 $&$ \pi/4 $&   $\pi/3 ,\pi/2 ,\pi/4 $ & $\Gamma_{20}$ & $\Delta_{(2,3,6)},D_6,D_4 $
\\
\hline
$6$ & 0.4781 & $ 0 $&$ \pi/2 $&$ \pi/5 $&   $\pi/3 ,\pi/2 ,\pi/5 $ & $\Gamma_{21}$ & $\Delta_{(2,3,6)},D_6,D_5 $
\\
\hline
$7$ & 0.5120 & $ 0 $&$ \pi/2 $&$ \pi/6 $&   $\pi/3 ,\pi/2 ,\pi/6 $ & $\Gamma_{22}$ & $\Delta_{(2,3,6)},D_6,D_6 $
\\
\hline
$8$ & 0.5330 & $ 0 $&$ \pi/2 $&$ \pi/n, \; 7\leq n\leq \infty $&  $\pi/3 ,\pi/2 ,\pi/n $ & $\Delta_{(2,3,n)}$ &
$\Delta_{(2,3,6)},D_6,D_n
$ \\
\hline
\end{tabular}
\end{center}}

\subsubsection{Dihedral angles not submultiples of $\pi$}

In this section we record those groups which arise from reflection groups whose dihedral angles are not submultiples of
$\pi$.   
{\scriptsize
\begin{center}
\begin{tabular}{|c|c|c|c|c|c|c|c|}
\multicolumn{8}{c}{\bf  $(p,q,r)$--triangles ($p\geq 3$). Dihedral angles not submultiples of $\pi$ }\\
\hline
$\#$ & $(p,q,r)$ & Angle & Angle & Angle & $\Sigma$ dihedral & group & Vertex Structure\\
\hline
$1^*$ & $(3,3,3)$ & $0.7297$ & $0.7297$ & $0.7297$ & $6\pi/5$ & $6 \times (2,2,3)_{12}$  &
$A_5,A_5,A_5$\\ 
\hline
$2^*$ & $(3,3,3)$ & $0$ & $0$ & $0.7297 $ & $16\pi/15$ & $2 \times (2,3,3)_9$ & $ A_5,\Delta_{(3,3,3)},\Delta_{(3,3,3)}$\\
\hline
$3^*$ & $(3,3,3)$ & $0$ & $0$ & $1.9106$ & $4\pi/3$ & $2 \times (2,3,3)_{14}$ & $A_4,\Delta_{(3,3,3)},\Delta_{(3,3,3)}$\\
\hline
$4^*$ & $(3,3,4)$ & $0.7297$ & $.955$ & $.955$ & $16\pi/15$ &  $2 \times (2,3,4)_{11}$ & $\Delta_{(3,3,3)},S_4,S_4$\\
\hline
$5^*$ & $(3,3,5)$ & $0.7297$ & $.652$ & $.652 $ & $7\pi/5$ & $2 \times (2,3,5)_{7}$ & $A_5,A_5,A_5$\\
\hline
$6^*$ & $(3,3,6)$ & $0.7297$ & $0$ & $0$ & $7\pi/5$ & $2\times (2,3,6)_7$ &
$A_5,\Delta_{(2,3,6)},\Delta_{(2,3,6)}$\\ 
\hline
$7^*$ & $(3,3,6)$ & $1.9106$ & $0$ & $0$ & $5\pi/3$ & $2 \times (2,3,6)_9$ & $ A_4,\Delta_{(2,3,6)},\Delta_{(2,3,6)}$ \\
\hline
$8^*$ & $(3,5,5)$ & $1.10719$ & $.652$ & $.652$ & $5\pi/3$ & $2\times (2,3,5)_2$ & 
$A_5,A_5,A_5$\\
\hline
$9^*$ & $(3,6,6)$ & $0$ & $0$ & $0$ & $5\pi/3$ & $2\times (2,3,6)_5$ & 
$A_5,A_5,A_5$\\
\hline
$10^*$ & $(6,6,6)$ & $0$ & $0$ & $0$ & $2\pi$ & $6 \times (2,2,6)_2$ & 
$A_5,A_5,A_5$\\
\hline
\end{tabular}
\end{center}}

{\scriptsize
\begin{center}
\begin{tabular}{|c|c|c|c|c|c|c|}
\multicolumn{7}{c}{\bf  (2,3,3)--triangles. Dihedral angles not submultiples of $\pi$ }\\
\hline
$\#$  &$3-3$&$2-3$ &$2-3$ & dihedral angles & group & Vertex Structure\\
\hline
$1^*$ &  $ 0.7297 $&$ 0 $&$ 0 $& $2\pi/5 ,\pi/6 , \pi/6 $ & $2 \times (2,2,3)_2$ &
$A_5,\Delta_{(2,3,6)},\Delta_{(2,3,6)}$
\\
\hline
$2^*$ &  $ 0.7297 $&$ 0.364 $&$ 0.364 $& $2\pi/5 ,\pi/5 , \pi/5 $ &$2 \times (2,2,3)_8$ &$A_5,A_5,A_5$
\\
\hline
$3^*$ &  $ 0.7297 $&$ 0.615 $&$ 0.615 $& $2\pi/5 ,\pi/4 , \pi/4 $ & $2 \times (2,2,3)_9,$&$A_5,S_4,S_4$ \\
\hline
$ 4^* $ & $ 0.7297 $&$ 0.955 $&$ 0.955 $& $2\pi/5 ,\pi/3 , \pi/3 $ & $2\times (2,2,3)_{10}$ &$A_4,A_4,A_5$
\\
\hline
$ 5^* $ & $ 1.9106 $&$ 0 $&$ 0 $&  $2\pi/3 ,\pi/6 , \pi/6 $ & $2\times (2,2,3)_4$ &$A_4,\Delta_{(2,3,6)},\Delta_{(2,3,6)} $ \\
\hline
$6^* $ &  $ 1.9106 $&$ .364 $&$ .364 $&  $2\pi/3 ,\pi/5 , \pi/5 $ &$2 \times (2,2,3)_{10}$ &$A_4,A_5,A_5$ \\
\hline
\end{tabular}
\end{center}}

{\scriptsize
\begin{center}
\begin{tabular}{|c|c|c|c|c|c|c|}
\multicolumn{7}{c}{\bf (2,p,q)--triangles. Dihedral angles not submultiples of $\pi$ }\\
\hline
$\#$  &$4-4$&$2-4$ &$2-4$ & dihedral angles & group & Vertex Structure\\
\hline
$1^*$ &  $ \pi/2 $&$ 0 $&$ 0 $& $2\pi/3 ,\pi/4 , \pi/4 $ & $2 \times (2,2,4)_3$ &
$S_4,\Delta_{(2,4,4)},\Delta_{(2,4,4)}$
\\
\hline
$\#$  &$5-5$&$2-5$ &$2-5$ &   &  &  \\
\hline
$1^*$ &  $ 1.107149 $&$ 0.5535 $&$ 0.5535 $& $2\pi/3 ,\pi/3 , \pi/3 $ &$2 \times (2,2,5)_1$ &$A_5,A_5,A_5$
\\
\hline
$\#$  &$6-6$&$2-6$ &$2-6$ &   &  &  \\
\hline
$1^*$ &  $ 0 $&$ 0 $&$ 0 $& $2\pi/3 ,\pi/3 , \pi/3 $ & $2 \times (2,2,6)_1,$&$\Delta(2,3,6),\Delta_{(2,3,6)},\Delta_{(2,3,6)}$
\\
\hline
$2^*$ &  $ 0 $&$ 0 $&$ \pi/2 $& $2\pi/3 ,\pi/3 , \pi/2 $ & $3 \times
(2,2,6)_2,$&$\Delta(2,3,6),\Delta_{(2,3,6)},D_6$
\\
\hline
\end{tabular}
\end{center}}

{\scriptsize
\begin{center}
\begin{tabular}{|c|c|c|c|c|c|c|}
\multicolumn{7}{c}{\bf  (2,2,3)--triangles. Dihedral angles not submultiples of $\pi$ }\\
\hline
$\#$  &$2-3$&$2-3$ &$2-2$ & dihedral angles & group & Vertex Structure\\
\hline
$1^*$ &  $ 0 $&$ 0 $&$ 2\pi/3 $&$\pi/6,\pi/6,2\pi/3 $& $2 \times (2,2,3)_6 $&$\Delta_{(2,3,6)},\Delta_{(2,3,6)},D_3$\\
\hline
$2^*$ &  $ 0 $&$ 0 $&$ 2\pi/n\;\;(5\leq n\leq \infty) $&$\pi/6,\pi/6,2\pi/n $&
$\Delta_{(2,6,n)}$&$\Delta_{(2,3,6)},\Delta_{(2,3,6)},D_n$\\
\hline
$3^*$ &  $ .364 $&$ .364 $&$ 2\pi/3 $ &$\pi/5 ,\pi/5 ,2\pi/3 $& $2 \times (2,2,3)_{12} $ & $A_5,A_5,D_3$ \\
\hline
$4^*$ &  $ .364 $&$ .364 $&$ 2\pi/n \;\; (2\leq n\leq \infty)  $ &$\pi/5 ,\pi/5 ,2\pi/n $& $\Delta_{(5,5,n)}$ & $A_5,A_5,D_n$
\\
\hline
$5^*$ &  $ .615 $&$ .615 $&$ 2\pi/n \;\; (5\leq n\leq \infty)  $& $\pi/4 ,\pi/4 ,2\pi/n $&  $\Delta_{(2,4,n)}$ & $S_4,S_4,D_n$
\\
\hline
$6^*$ &  $ .955 $&$ .955 $&$ 2\pi/n \;\; (7 \leq n\leq \infty)  $& $\pi/3 ,\pi/3 ,2\pi/n $& $\Delta_{(2,3,n)}$ &
$A_4,A_4,D_n$
\\
\hline
\end{tabular}
\end{center}}

{\scriptsize
\begin{center}
\begin{tabular}{|c|c|c|c|c|c|c|}
\multicolumn{7}{c}{\bf   (2,2,4)--triangles. Dihedral angles not submultiples of $\pi$ }\\
\hline
$\# $ &$2-4$&$2-4$ &$2-2$ & dihedral angles & group & Vertex Structure\\
\hline
$1^*$ & $ 0 $&$ 0 $&$ 2\pi/3 $& $\pi/4 ,\pi/4 ,2\pi/3 $ &
$2 \times (2,2,4)_5 $ & $\Delta_{(2,4,4)},\Delta_{(2,4,4)},D_3$ \\
\hline
$2^*$ & $ 0 $&$ 0 $&$ 2\pi/n \;\; (5\leq n \leq \infty) $& $\pi/4 ,\pi/4 ,2\pi/n $ &
$\Delta_{(2,4,n)}$ & $\Delta_{(2,4,4)},\Delta_{(2,4,4)},D_n$ \\
\hline
$3^*$ & $ \pi/4 $&$ \pi/4 $&$ 2\pi/5  $& $\pi/3 ,\pi/3 ,2\pi/5$ & $2 \times (2,2,4)_{10} $ &
$S_4,S_4,D_5$  \\
\hline
$4^*$ & $ \pi/4 $&$ \pi/4 $&$ 2\pi/n \;\; (7\leq n \leq \infty) $& $\pi/3 ,\pi/3 ,2\pi/n$ & $\Delta_{(2,3,n)}$ &
$S_4,S_4,D_n$  \\
\hline
\end{tabular}
\end{center}}

{\scriptsize
\begin{center}
\begin{tabular}{|c|c|c|c|c|c|c|}
\multicolumn{7}{c}{\bf   (2,2,5)--triangles. Dihedral angles not submultiples of $\pi$ }\\
\hline
$\#$ &$2-5$&$2-5$ &$2-2$ & dihedral angles & group & Vertex Structure\\
\hline
$1^*$ & $ .55357 $&$ .55357 $&$ 2\pi/5, $& $\pi/3
,\pi/3 ,2\pi/5 $ & $2 \times (2,2,5)_5 $  & 
$A_5,A_5,D_5$ \\
\hline
$2^*$ & $ .55357 $&$ .55357 $&$ 2\pi/n,
\;\; (7\leq n \leq \infty) $& $\pi/3
,\pi/3 ,2\pi/n $ & $\Delta(2,3,n)$  & 
$A_5,A_5,D_n$ \\
\hline
\end{tabular}
\end{center}}

{\scriptsize
\begin{center}
\begin{tabular}{|c|c|c|c|c|c|c|}
\multicolumn{7}{c}{\bf  (2,2,6)--triangles. Dihedral angles not submultiples of $\pi$ }\\
\hline
$\#$  &$2-6$&$2-6$ &$2-2$ & dihedral angles & group & Vertex Structure\\
\hline
$1^*$  & $ 0 $&$ 0 $&$ 2\pi/3 $ & $\pi/3 ,\pi/3 ,2\pi/3 $ & $2\times (2,2,6)_4$ & $\Delta_{(2,3,6)},\Delta_{(2,3,6)},D_3$ \\
\hline
$2^*$  & $ 0 $&$ 0 $&$ 2\pi/5 $ & $\pi/3 ,\pi/3 ,2\pi/5 $ & $2\times (2,2,6)_6$ & $\Delta_{(2,3,6)},\Delta_{(2,3,6)},D_5$ \\
\hline
$3^*$  & $ 0 $&$ 0 $&$ 2\pi/n,\; 7\leq n\leq \infty $  & $\pi/3 ,\pi/3 ,2\pi/n $& $\Delta_{(2,3,n)}$
& $\Delta_{(2,3,6)},\Delta_{(2,3,6)},D_n$ \\
\hline
\end{tabular}
\end{center}}

\subsubsection{The reflection groups}

Here, for reference, we collect together the tetrahedral reflection groups which occur and
tabulate them so that they can be identified in our previous tables.  First contains the classical reflection groups with all
dihedral angles submultiples of $\pi$. 
{\scriptsize
\begin{center}
\begin{tabular}{|c|c|c|c|c|c|c|c|}
\multicolumn{8}{c}{\bf  Reflection groups. Angles submultiples of $\pi$}\\
\hline
$\G_i$ &$\angle\overline{AB}$ &$\angle\overline{BC}$ &$\angle\overline{CD}$ &$\angle\overline{AC}$&$\angle\overline{AD}$&
$\angle\overline{BD}$ & covolume\\
\hline
1 & $ \pi/3 $ & $ \pi/2 $ & $ \pi/3 $ & $\pi/2 $ & $\pi/5 $ & $\pi/2  $  & 0.039050 \\
\hline
2 & $ \pi/5 $ & $ \pi/2 $ & $ \pi/4 $ & $\pi/2 $ & $\pi/3 $ & $\pi/2  $  & 0.035885 \\
\hline
3 & $ \pi/5 $ & $ \pi/2 $ & $ \pi/5 $ & $\pi/2 $ & $\pi/3 $ & $\pi/2  $  & 0.093326 \\
\hline
4 & $ \pi/3 $ & $ \pi/3 $ & $ \pi/3 $ & $\pi/2 $ & $\pi/4 $ & $\pi/2  $ & 0.085770 \\
\hline
5 & $ \pi/3 $ & $ \pi/3 $ & $ \pi/3 $ & $\pi/2 $ & $\pi/5 $ & $\pi/2  $  & 0.205289 \\
\hline
6 & $ \pi/3 $ & $ \pi/4 $ & $ \pi/3 $ & $\pi/2 $ & $\pi/4 $ & $\pi/2  $ & 0.222229 \\
\hline
7 & $ \pi/3 $ & $ \pi/5 $ & $ \pi/3 $ & $\pi/2 $ & $\pi/4 $ & $\pi/2  $  & 0.358653\\
\hline
8 & $ \pi/3 $ & $ \pi/5 $ & $ \pi/3 $ & $\pi/2 $ & $\pi/5 $ & $\pi/2  $  & 0.502131\\
\hline
9 & $ \pi/2 $ & $ \pi/3 $ & $ \pi/3 $ & $\pi/5 $ & $\pi/3 $ & $\pi/2  $  & 0.071770\\
\hline
10 & $ \pi/3 $ & $ \pi/3 $ & $ \pi/6 $ & $\pi/2 $ & $\pi/3 $ & $\pi/2  $ & 0.364107 \\
\hline
11 & $ \pi/4 $ & $ \pi/3 $ & $ \pi/6 $ & $\pi/2 $ & $\pi/3 $ & $\pi/2  $  & 0.525840 \\
\hline
12 & $ \pi/5 $ & $ \pi/3 $ & $ \pi/6 $ & $\pi/2 $ & $\pi/3 $ & $\pi/2  $ & 0.672986  \\
\hline
13 & $ \pi/6 $ & $ \pi/3 $ & $ \pi/6 $ & $\pi/2 $ & $\pi/3 $ & $\pi/2  $ & 0.845785 \\
\hline
14 & $ \pi/3 $ & $ \pi/3 $ & $ \pi/4 $ & $\pi/2 $ & $\pi/4 $ & $\pi/2  $ & 0.305322\\
\hline
15 & $ \pi/3 $ & $ \pi/4 $ & $ \pi/4 $ & $\pi/2 $ & $\pi/4 $ & $\pi/2  $ & 0.556282 \\
\hline
16 & $ \pi/4 $ & $ \pi/4 $ & $ \pi/4 $ & $\pi/2 $ & $\pi/4 $ & $\pi/2  $ & 0.915965 \\
\hline
17 & $ \pi/6 $ & $ \pi/2 $ & $ \pi/3 $ & $\pi/2 $ & $\pi/3 $ & $\pi/2  $ & 0.042289 \\
\hline
18 & $ \pi/3 $ & $ \pi/2 $ & $ \pi/3 $ & $\pi/2 $ & $\pi/6 $ & $\pi/2  $  & 0.169157\\
\hline
19 & $ \pi/4 $ & $ \pi/2 $ & $ \pi/3 $ & $\pi/2 $ & $\pi/4 $ & $\pi/2  $  & 0.076330\\
\hline
20 & $ \pi/6 $ & $ \pi/2 $ & $ \pi/4 $ & $\pi/2 $ & $\pi/3 $ & $\pi/2  $  & 0.105723\\
\hline
21 & $ \pi/6 $ & $ \pi/2 $ & $ \pi/5 $ & $\pi/2 $ & $\pi/3 $ & $\pi/2  $  & 0.171502\\
\hline
22 & $ \pi/6 $ & $ \pi/2 $ & $ \pi/6 $ & $\pi/2 $ & $\pi/3 $ & $\pi/2  $ & 0.253735 \\
\hline
23 & $ \pi/4 $ & $ \pi/2 $ & $ \pi/4 $ & $\pi/2 $ & $\pi/4 $ & $\pi/2  $ & 0.228991\\
\hline
24 & $ \pi/2 $ & $ \pi/2 $ & $ \pi/3 $ & $\pi/6 $ & $\pi/3 $ & $\pi/2  $ & 0.211446 \\
\hline
25 & $ \pi/2 $ & $ \pi/2 $ & $ \pi/3 $ & $\pi/4 $ & $\pi/4 $ & $\pi/2  $ & 0.152661 \\
\hline
26 & $ \pi/2 $ & $ \pi/2 $ & $ \pi/4 $ & $\pi/4 $ & $\pi/4 $ & $\pi/2  $  & 0.457983\\
\hline
27 & $ \pi/3 $ & $ \pi/2 $ & $ \pi/3 $ & $\pi/3 $ & $\pi/3 $ & $\pi/2  $  & 0.084578\\
\hline
28 & $ \pi/3 $ & $ \pi/2 $ & $ \pi/4 $ & $\pi/3 $ & $\pi/3 $ & $\pi/2  $  & 0.211446\\
\hline
29 & $ \pi/3 $ & $ \pi/2 $ & $ \pi/5 $ & $\pi/3 $ & $\pi/3 $ & $\pi/2  $ & 0.343003\\
\hline
30 & $ \pi/3 $ & $ \pi/2 $ & $ \pi/6 $ & $\pi/3 $ & $\pi/3 $ & $\pi/2  $ & 0.507471\\
\hline
31 & $ \pi/3 $ & $ \pi/3 $ & $ \pi/3 $ & $\pi/3 $ & $\pi/3 $ & $\pi/3  $  & 1.014941\\
\hline
32 & $ \pi/3 $ & $ \pi/3 $ & $ \pi/3 $ & $\pi/3 $ & $\pi/3 $ & $\pi/2  $  & 0.422892\\
\hline
\end{tabular}
\end{center}}

{\scriptsize
\begin{center}
\begin{tabular}{|c|c|c|c|c|c|c|}
\multicolumn{7}{c}{\bf  Reflection groups. Non--submultiples of $\pi$}\\
\hline
$\G_{i}^{*}$ &$\angle\overline{AB}$ &$\angle\overline{AC}$ &$\angle\overline{AD}$ &$\angle\overline{BC}$&$\angle\overline{BD}$&
$\angle\overline{CD}$ \\
\hline
$1^*$ & $ 2\pi/5 $ & $ 2\pi/5 $ & $ 2\pi/5 $ & $\pi/3 $ & $\pi/3 $ & $\pi/3  $  \\
\hline
$2^*$ & $ 2\pi/5 $ & $ \pi/3 $ & $ \pi/3 $ & $\pi/3 $ & $\pi/3 $ & $\pi/3  $  \\
\hline
$3^*$ & $ 2\pi/5 $ & $ \pi/3 $ & $ \pi/3 $ & $\pi/3 $ & $\pi/3 $ & $\pi/4  $  \\
\hline
$4^*$ & $ 2\pi/5 $ & $ \pi/2 $ & $ \pi/2 $ & $\pi/3 $ & $\pi/3 $ & $\pi/5  $  \\
\hline
$5^*$ & $ 2\pi/5 $ & $ \pi/2 $ & $ \pi/2 $ & $\pi/3 $ & $\pi/3 $ & $\pi/6  $  \\
\hline
$6^*$ & $ 2\pi/5 $ & $ \pi/6 $ & $ \pi/6 $ & $\pi/3 $ & $\pi/3 $ & $\pi/2  $  \\
\hline
$7^*$ & $ 2\pi/5 $ & $ \pi/5 $ & $ \pi/5 $ & $\pi/3 $ & $\pi/3 $ & $\pi/2  $  \\
\hline
$8^*$ & $ 2\pi/5 $ & $ \pi/4 $ & $ \pi/4 $ & $\pi/3 $ & $\pi/3 $ & $\pi/2  $  \\
\hline
$9^*$ & $ 2\pi/5 $ & $ \pi/3 $ & $ \pi/3 $ & $\pi/3 $ & $\pi/3 $ & $\pi/2  $  \\
\hline
$10^*$ & $ 2\pi/5 $ & $ \pi/3 $ & $ \pi/3 $ & $\pi/2 $ & $\pi/2 $ & $\pi/4  $  \\
\hline
$11^*$ & $ 2\pi/3 $ & $ \pi/3 $ & $ \pi/3 $ & $\pi/3 $ & $\pi/3 $ & $\pi/3  $  \\
\hline
$12^*$ & $ 2\pi/3 $ & $ \pi/2 $ & $ \pi/2 $ & $\pi/3 $ & $\pi/3 $ & $\pi/6  $  \\
\hline
$13^*$ & $ 2\pi/3 $ & $ \pi/2 $ & $ \pi/2 $ & $\pi/5 $ & $\pi/5 $ & $\pi/3  $  \\
\hline
$14^*$ & $ 2\pi/3 $ & $ \pi/2 $ & $ \pi/2 $ & $\pi/6 $ & $\pi/6 $ & $\pi/3  $  \\
\hline
$15^*$ & $ 2\pi/3 $ & $ 2\pi/3 $ & $ 2\pi/3 $ & $\pi/6 $ & $\pi/6 $ & $\pi/6  $  \\
\hline
$16^*$ & $ 2\pi/3 $ & $ \pi/6 $ & $ \pi/6 $ & $\pi/3 $ & $\pi/3 $ & $\pi/2  $  \\
\hline
$17^*$ & $ 2\pi/3 $ & $ \pi/5 $ & $ \pi/5 $ & $\pi/3 $ & $\pi/3 $ & $\pi/2  $  \\
\hline
$18^*$ & $ 2\pi/3 $ & $ \pi/4 $ & $ \pi/4 $ & $\pi/4 $ & $\pi/4 $ & $\pi/2  $  \\
\hline
$19^*$ & $ 2\pi/3 $ & $ \pi/4 $ & $ \pi/4 $ & $\pi/2 $ & $\pi/2 $ & $\pi/4  $  \\
\hline
\end{tabular}
\end{center}}

\end{document}